\documentclass[10pt]{article}
\usepackage{amsmath}
\usepackage{amsthm,amsfonts,amssymb,epsfig,graphics}
\usepackage[latin1]{inputenc}
\usepackage[francais,english]{babel}
\usepackage{pst-all}
\def\ie{\emph{i.e.} }
\def\cf{\emph{cf.} }

\def\inpa{\mathrm{IP}}

\def\indi{\mathrm{Indice}}
\def\adhe{\mathrm{Adh}}
\def\inte{\mathrm{Int}}

\newcommand{\bbR}{{\mathbb{R}}}
\newcommand{\bbC}{{\mathbb{C}}}
\newcommand{\bbZ}{{\mathbb{Z}}}

\newcommand{\bbS}{{\mathbb{S}}}
\def\cF{{\cal F}}
\def\cA{{\cal A}}
\def\cG{{\cal G}}
\def\cP{{\cal P}}

\def\fixe{\mathrm{Fixe}}

\def\woi{W_{0 \rightarrow \infty}}
\def\wio{W_{\infty \rightarrow 0}}
\def\woo{W_{0 \rightarrow 0}}
\def\wii{W_{\infty \rightarrow \infty}}

\def\wao{W_{\alpha \rightarrow \omega}}

\def\wt{\widetilde}
\def\t{\tilde }
\def\diam{\mathrm{Module}}

\def\?{$^{***}$\marginpar{?}}

\def\cA{{\cal A}}

\newtheorem{theo}{Th\'eor\`eme}

\newtheorem*{ques*}{Question}
\newtheorem*{prop*}{Proposition}
\newtheorem*{conj*}{Conjecture}
\newtheorem*{theo*}{Th\'eor\`eme}
\newtheorem{coro}{Corollaire}[section]

\newtheorem{affi}[coro]{Affirmation}
\newtheorem{affi*}{Affirmation}
\newtheorem{prop}[coro]{Proposition}
\newtheorem{lemm}[coro]{Lemme}
\def\?{\footnote{?}}

\newenvironment{demo}[1][]{\addvspace{8mm} \emph{Preuve #1
    ---~~}}{~~~$\square$\bigskip}

\newlength{\espaceavantspecialthm}
\newlength{\espaceapresspecialthm}
\newenvironment{defi}[1][]{\refstepcounter{coro} %ref=> ca labelise bien!
\vskip \espaceavantspecialthm \noindent \textbf{D\'efinition~\thecoro
#1.} }%
{\vskip \espaceapresspecialthm}

\newenvironment{rema}[1][]{\refstepcounter{coro} %ref=> ca labelise bien!
\vskip \espaceavantspecialthm \noindent \textbf{Remarque~\thecoro
#1.} }%
{\vskip \espaceapresspecialthm}

\author{Fr\'ed\'eric Le Roux \\
Universit\'e Paris Sud \\
Bat. 430 \\
91405 Orsay Cedex   FRANCE \\
e-mail : leroux@topo.math.u-psud.fr}
\title{Un indice qui affine l'indice de Poincar\'e-Lefschetz pour les hom\'eomorphismes de surfaces}

\def\bbRo{\bbR^2 \setminus\{0\}}

\def\h01{{\cal H}_{0 \neq 1}}

\def\la{\leftarrow}
\def\ra{\rightarrow}
\def\ua{\uparrow}
\def\da{\downarrow}

\begin{document}
\sloppy

%************************************
%************************************
%************************************
%\pagebreak
%%%%%%%%%%%%%%%%%%%%%%%%%%%%%%%
%%%%%%%%%%%%%%%%%%%%%%%%%%%%%%%%

%
%RESTE : 

%- Barcelonais

%- r\'ealisabilit\'e et figure (secteur indiff\'erent notamment) ;

%- VERIF : compatibilit\'e composantes

%- VERIF : ordre cyclique, def de $O(F,G)$

%RESTE :

%A) Cons\'equences dynamiques : en local ??

%B) Appendice graphe : preuves.

%1) Figures

%2) structures dynamiques sur un exemple (fin)

%3) utilisation de l'ordre (cyclique) est conforme \`a la def ?

%4)  plan de la partie ``preuves des lemmes'', p21...

%5) Intro de ``cons\'equences dynamiques''. Lien questions Asterisque.

\clearpage
\pagenumbering{arabic}
\maketitle

%\pagebreak
%\setcounter{\page}{0}
%%%%%%%%%%%
%%%%%%%%%%%

%\maketitle
\selectlanguage{francais}
\begin{abstract}
Nous \'etudions la dynamique des hom\'eomorphismes de surfaces autour des points fixes isol\'es dont l'indice de Poincar\'e-Lefschetz est diff\'erent de $1$. Nous d\'efinissons un invariant de conjugaison, mot cyclique sur l'alphabet $\{\ua, \ra , \da , \la\}$, qui affine l'indice de Poincar\'e-Lefschetz. Grossi\`erement, il s'agit de d\'ecomposer canoniquement la dynamique en un  nombre fini de secteurs hyperboliques,  elliptiques ou indiff\'erents, chaque type de secteur contribuant \`a l'indice de Poincar\'e-Lefschetz par un terme valant respectivement $-1/2$,  $+1/2$ ou $0$. 
Le mot cyclique permet de lire l'existence de structures dynamiques canoniques. 
%en particulier des versions d\'eg\'en\'er\'ees \emph{canoniques} des vari\'et\'es stables et instables de la dynamique diff\'erentiable hyperbolique, ou des p\'etales de la dynamique holomorphe parabolique.
\end{abstract}

\selectlanguage{english}
\begin{abstract}
We study the dynamics of surface homeomorphisms around isolated fixed points whose Poincar\'e-Lefschetz index is not equal to $1$. We construct a new conjugacy invariant, which is a cyclic word on the alphabet $\{\ua, \ra , \da , \la\}$. This invariant is a refinement of the P.-L. index. It can be seen as a canonical decomposition of the dynamics into a finite number of sectors of hyperbolic, elliptic or indifferent type. The contribution of each type of sector to the P.-L. index is respectively 
$-1/2$,  $+1/2$ and $0$. The construction of the invariant implies the existence of some canonical dynamical structures.
\end{abstract}
\selectlanguage{francais}
\bigskip	
\bigskip

\textbf{Mots-clés} 
{Homéomorphisme, surface, point fixe, Brouwer, indice, composantes de
Reeb.}
\bigskip	

\bigskip	

\textbf{Classification AMS (2000)} 
{37E30, 37B30}.

\clearpage
\setcounter{tocdepth}{1}
\tableofcontents
\clearpage

%\bigskip
%\bigskip

\part{Introduction et r\'esultats}
\label{p.resultats}

%%%%%%%%%%%%%%%%%%%%%%%%%%%%
\section{Introduction}

%%%%%%%%%%%%%%%%%%%%%%%
\subsection{R\'esultat principal}
Nous nous proposons de construire un invariant pour un certain type de points fixes des  syst\`emes dynamiques topologiques en dimension $2$. Plus pr\'ecis\'ement, nous consid\'erons
 les donn\'ees suivantes~:
\begin{itemize}
\item une surface orientable $S$, pas n\'ecessairement compacte~;
\item un point $x_{0}$ sur la surface $S$~;
\item un hom\'eomorphisme $h$ de la surface $S$ dans elle-m\^eme, qui pr\'eserve l'orientation, et qui fixe le point $x_{0}$. 
\end{itemize}
Nous faisons deux hypoth\`eses importantes :
\begin{itemize}
\item le point $x_{0}$ est isol\'e parmi les points fixes de $h$~;
\item l'indice de Poincar\'e-Lefschetz de ce point fixe est diff\'erent de $1$.
\end{itemize}

Le principal r\'esultat du texte aura la forme suivante.
\pagebreak[3]
\begin{theo*}
Il existe un invariant de conjugaison qui associe \`a ces donn\'ees~:
\begin{itemize}
\item un entier $d \geq 2$ (appel\'e le \emph{module} de $h$ en $x_{0}$)~;
\item un mot cyclique $\mathbf{M}(h)$ sur l'alphabet $\{\la, \ra, \ua,\da  \}$.
\end{itemize}
Le mot $\mathbf{M}(h)$ est compos\'e de $2 \times d$ lettres, qui sont alternativement des fl\`eches horizontales et verticales.
\end{theo*}
\`A vrai dire, il faut exclure de l'\'enonc\'e le cas o\`u $S$ est la sph\`ere et $x_{0}$ l'unique point fixe de $h$ (ce qu'on appelle un \emph{hom\'eomorphisme de Brouwer}).
D'autre part, la notion de module se g\'en\'eralise quand l'indice est \'egal \`a $1$, et le th\'eor\`eme se g\'en\'eralise en indice $1$ lorsque le module est sup\'erieur ou \'egal \`a $4$.

\begin{figure}[htbp]
\begin{center}
\psset{xunit=1mm,yunit=1mm,runit=1mm}
\psset{linewidth=0.3,dotsep=1,hatchwidth=0.3,hatchsep=1.5,shadowsize=1}
\psset{dotsize=0.7 2.5,dotscale=1 1,fillcolor=black}
\psset{arrowsize=1 2,arrowlength=1,arrowinset=0.25,tbarsize=0.7 5,bracketlength=0.15,rbracketlength=0.15}
\begin{pspicture}(0,0)(151.25,41.25)
\psline(15,37.5)(15,12.5)
\psline(2.5,25)(27.5,25)
\psline(32.5,25)(57.5,25)
\psline(62.5,37.5)(87.5,12.5)
\psline(62.5,12.5)(87.5,37.5)
\psline(75,37.5)(75,12.5)
\psline(62.5,25)(87.5,25)
\psline(92.5,25)(117.5,25)
\psline(122.5,25)(147.5,25)
\psline(135,37.5)(135,12.5)
\psbezier(2.5,27.5)(10,27.5)(12.5,30)(12.5,37.5)
\psbezier(2.5,26.25)(13.75,26.25)(13.75,26.25)(13.75,37.5)
\psbezier(2.5,28.75)(6.25,28.75)(11.25,33.75)(11.25,37.5)
\psbezier(12.5,12.5)(12.5,20)(10,22.5)(2.5,22.5)
\psbezier(13.75,12.5)(13.75,23.75)(13.75,23.75)(2.5,23.75)
\psbezier(11.25,12.5)(11.25,16.25)(6.25,21.25)(2.5,21.25)
\psbezier(27.5,22.5)(20,22.5)(17.5,20)(17.5,12.5)
\psbezier(27.5,23.75)(16.25,23.75)(16.25,23.75)(16.25,12.5)
\psbezier(27.5,21.25)(23.75,21.25)(18.75,16.25)(18.75,12.5)
\psbezier(17.5,37.5)(17.5,30)(20,27.5)(27.5,27.5)
\psbezier(16.25,37.5)(16.25,26.25)(16.25,26.25)(27.5,26.25)
\psbezier(18.75,37.5)(18.75,33.75)(23.75,28.75)(27.5,28.75)
\psline{->}(15,37.5)(15,30)
\psline{->}(15,25)(10,25)
\psline{->}(15,12.5)(15,20)
\psline{->}(15,25)(20,25)
\psbezier(45,25)(45,25)(45,25)(45,25)
\psdots[](45,25)
(45,25)
\psline(32.5,27.5)(57.5,27.5)
\psline(32.5,22.5)(57.5,22.5)
\psline(32.5,20)(57.5,20)
\psline(32.5,17.5)(57.5,17.5)
\psline(32.5,30)(57.5,30)
\psline(32.5,32.5)(57.5,32.5)
\psline{->}(32.5,27.5)(46.25,27.5)
\psline{->}(32.5,30)(46.25,30)
\psline{->}(32.5,25)(40,25)
\psline{->}(32.5,32.5)(46.25,32.5)
\psline{->}(32.5,22.5)(46.25,22.5)
\psline{->}(32.5,20)(46.25,20)
\psline{->}(32.5,17.5)(46.25,17.5)
\psline{->}(45,25)(51.25,25)
\psline{->}(62.5,37.5)(70,30)
\psline{->}(75,37.5)(75,30)
\psline{->}(87.5,37.5)(80,30)
\psline{->}(87.5,25)(80,25)
\psline{->}(87.5,12.5)(80,20)
\psline{->}(75,12.5)(75,20)
\psline{->}(62.5,12.5)(70,20)
\psline{->}(62.5,25)(70,25)
\psline{->}(92.5,25)(97.5,25)
\psline{->}(105,25)(115,25)
\psbezier(105,25)(93.75,25)(93.75,13.75)(105,13.75)
\psbezier(105,25)(97.12,25)(97.12,17.5)(105,17.5)
\psbezier(105,25)(100.5,25)(100.5,20)(105,20)
\psbezier{->}(105,25)(117.5,25)(117.5,36.25)(105,36.25)
\psbezier{->}(105,25)(113.75,25)(113.75,32.5)(105,32.5)
\psbezier{->}(105,25)(110,25)(110,30)(105,30)
\psbezier{->}(105,25)(117.5,25)(117.5,13.75)(105,13.75)
\psbezier{->}(105,25)(113.75,25)(113.75,17.5)(105,17.5)
\psbezier{->}(105,25)(110,25)(110,20)(105,20)
\psbezier(106.25,25)(93.75,25)(93.75,36.25)(106.25,36.25)
\psbezier(106.25,25)(97.5,25)(97.5,32.5)(106.25,32.5)
\psbezier(106.25,25)(101.25,25)(101.25,30)(106.25,30)
\psbezier{->}(135,25)(135,41.25)(128.75,41.25)(123.75,36.25)
\psbezier{->}(135,25)(135,36.25)(130,37.5)(126.25,33.75)
\psbezier{->}(135,25)(135,31.25)(132.5,35)(128.75,31.25)
\psbezier(135,25)(118.75,25)(118.75,31.25)(123.75,36.25)
\psbezier(135,25)(123.75,25)(122.5,30)(126.25,33.75)
\psbezier(135,25)(128.75,25)(125,27.5)(128.75,31.25)
\psbezier{->}(135,25)(135,41.25)(141.25,41.25)(146.25,36.25)
\psbezier{->}(135,25)(135,36.25)(140,37.5)(143.75,33.75)
\psbezier{->}(135,25)(135,31.25)(137.5,35)(141.25,31.25)
\psbezier(135,25)(151.25,25)(151.25,31.25)(146.25,36.25)
\psbezier(135,25)(146.25,25)(147.5,30)(143.75,33.75)
\psbezier(135,25)(141.25,25)(145,27.5)(141.25,31.25)
\psbezier{->}(135,25)(135,8.75)(141.25,8.75)(146.25,13.75)
\psbezier{->}(135,25)(135,13.75)(140,12.5)(143.75,16.25)
\psbezier{->}(135,25)(135,18.75)(137.5,15)(141.25,18.75)
\psbezier(135,25)(151.25,25)(151.25,18.75)(146.25,13.75)
\psbezier(135,25)(146.25,25)(147.5,20)(143.75,16.25)
\psbezier(135,25)(141.25,25)(145,22.5)(141.25,18.75)
\psbezier{->}(135,25)(135,8.75)(128.75,8.75)(123.75,13.75)
\psbezier{->}(135,25)(135,13.75)(130,12.5)(126.25,16.25)
\psbezier{->}(135,25)(135,18.75)(132.5,15)(128.75,18.75)
\psbezier(135,25)(118.75,25)(118.75,18.75)(123.75,13.75)
\psbezier(135,25)(123.75,25)(122.5,20)(126.25,16.25)
\psbezier(135,25)(128.75,25)(125,22.5)(128.75,18.75)
\rput(15,5){$-1$}
\rput(45,5){$0$}
\rput(75,5){$1$}
\rput(105,5){$2$}
\rput(135,5){$3$}
\end{pspicture}
\caption{Diff\'erentes valeurs de l'indice}
\label{f.indice-PL}
\end{center}
\end{figure}
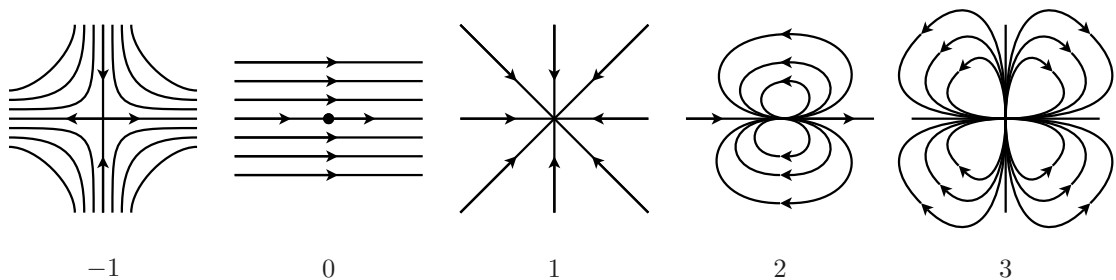

Rappelons la d\'efinition informelle de l'indice de Poincar\'e-Lefschetz du point fixe $x_{0}$. Il s'agit d'un invariant local, on peut se ramener au cas o\`u la surface $S$ est le plan. Puisque le point fixe est suppos\'e isol\'e, nous pouvons choisir un disque topologique $D$, contenant le point fixe dans son int\'erieur, et suffisamment petit pour ne contenir aucun autre point fixe de $h$. Pour chaque point $x$ du disque $D$,  tra\c cons le vecteur allant du point $x$ \`a son image $h(x)$~: on obtient ainsi un champ de vecteurs sur $D$, qui ne s'annule qu'au point $x_{0}$.
Nous pouvons maintenant \'evaluer l'indice de ce champ de vecteurs sur le bord de $D$~: 
lorsque le point $x$ parcourt la circonf\'erence de $D$, le vecteur reliant $x$ \`a $h(x)$ effectue un certain nombre de tours, qui donne l'indice du champ de vecteurs. On peut montrer que cet indice ne d\'epend pas du choix du disque $D$~: on a ainsi d\'efini l'indice de Poincar\'e-Lefschetz du point fixe $x_{0}$ pour l'hom\'eomorphisme $h$.
Les exemples  les plus simples sont dessin\'es sur la figure~\ref{f.indice-PL}.

\begin{figure}[htbp]
\begin{center}
\def\JPicScale{0.8}
\ifx\JPicScale\undefined\def\JPicScale{1}\fi
\psset{unit=\JPicScale mm}
\psset{linewidth=0.3,dotsep=1,hatchwidth=0.3,hatchsep=1.5,shadowsize=1}
\psset{dotsize=0.7 2.5,dotscale=1 1,fillcolor=black}
\psset{arrowsize=1 2,arrowlength=1,arrowinset=0.25,tbarsize=0.7 5,bracketlength=0.15,rbracketlength=0.15}
\begin{pspicture}(0,0)(36.96,31.25)
\psbezier(18.21,10)(18.21,10)(18.21,10)(18.21,10)
\rput{0}(18.21,15){\psellipse[linewidth=0.2,linestyle=dashed,dash=1 1](0,0)(15,15)}
\psbezier(18.21,15)(18.21,15)(18.21,15)(18.21,15)
\psdots[](18.21,15)
(18.21,15)
\psline{->}(18.21,0)(18.21,4.29)
\psline{->}(18.21,30)(13.57,30)
\psline{->}(32.32,9.82)(36.96,7.86)
\psline{->}(26.51,2.5)(30.17,5.36)
\psline{->}(9.46,27.14)(6.78,31.25)
\psline{->}(26.25,27.59)(23.57,23.39)
\psline{->}(32.32,20)(30.53,24.64)
\psline{->}(3.92,19.91)(5.53,24.64)
\psline{->}(4.01,9.91)(0,7.5)
\psline{->}(9.19,2.95)(5.71,5.71)
\rput(18.21,10){$0$}
\rput(18.75,11.25){}
\end{pspicture}
\caption{Lien avec l'indice de Poincar\'e-Lefschetz pour le mot $\mathbf{M}=(\ua \ra \da \ra \ua \ra \da \la \da \la )$ (pr\'ecis\'e en section~\ref{ss.exemple})}
\label{f.mot-cyclique}
\end{center}
\end{figure}
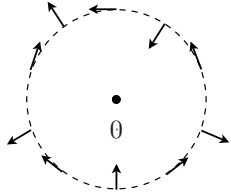

L'invariant de conjugaison que nous voulons d\'ecrire sera plus fin que l'indice de Poincar\'e-Lefschetz.  Pour pr\'eciser ceci, il faut penser au mot cyclique $\mathbf{M}(h)$ de la mani\`ere suivante. On dessine un cercle autour du point fixe $x_{0}$, et on \'ecrit les lettres du mot $\mathbf{M}(h)$ le long de ce cercle, comme sur la figure~\ref{f.mot-cyclique} (il faut s'imaginer \`a l'ext\'erieur du cercle, tourn\'e vers le point fixe,  \'ecrivant le mot $\mathbf{M}(h)$ en faisant le tour du cercle dans le sens trigonom\'etrique, en comman\c cant par la fl\`eche du bas). On a maintenant le dessin d'un champ de vecteurs le long du cercle. L'indice de Poincar\'e-Lefschetz est alors \'egal \`a l'indice de ce champ de vecteurs symbolique, c'est-\`a-dire au nombre de tours effectu\'es par les fl\`eches-lettres du mot $\mathbf{M}(h)$ lorsqu'on a \'ecrit ce mot autour du cercle.
Ainsi, le nouvel invariant donne une mani\`ere canonique de d\'ecomposer le nombre de tours \'evalu\'e par l'indice de Poincar\'e-Lefschetz en une somme de quarts-de-tour dans les sens positifs ou n\'egatifs, somme dont les termes sont cycliquement ordonn\'es.

%%%%%%%%%%%%%%%%%%%%%%%
\subsection{Esquisse de d\'efinition}
Esquissons la d\'efinition de l'invariant. Pour cela, on se place dans le cadre d'un hom\'eomorphisme du plan qui poss\`ede un unique point fixe : nous verrons qu'on peut toujours se ramener \`a cette situation id\'eale. 
Avant tout, rappelons que lorsque l'indice de Poincar\'e-Lefschetz du point fixe est suppos\'e diff\'erent de $1$, la th\'eorie de Brouwer donne une premi\`ere description de la dynamique~: toutes les orbites, dans le pass\'e ou le futur, tendent vers le point fixe ou vers l'infini. En particulier, il n'y a que quatre types d'orbites $(h^n(x))_{n \in \bbZ}$ possibles.
\'Etant donn\'ee une courbe ferm\'ee $\gamma$ qui entoure le point fixe, on peut d\'ecomposer $\gamma$ en un nombre fini de morceaux (connexes) \emph{libres}, c'est-\`a-dire qui sont chacun  disjoints de leur image par $h$. Lorsque le nombre de morceaux est le plus petit possible, on dit qu'on a une \emph{d\'ecomposition minimale} de $\gamma$, et 
on d\'efinit la \emph{$h$-longueur} de la courbe $\gamma$ comme \'etant pr\'ecis\'ement le nombre de morceaux d'une d\'ecomposition minimale. On peut alors d\'efinir le \emph{module} de $h$~: c'est la plus petite $h$-longueur parmi toutes les courbes ferm\'ees entourant le point fixe. Les courbes r\'ealisant ce minimum sont appel\'ees \emph{courbes g\'eod\'esiques}.

Passons \`a la d\'efinition du mot cyclique $\mathbf{M}(h)$. On construit d'abord les fl\`eches ``horizontales'' ($\ra$ et $\la$). On consid\`ere une d\'ecomposition minimale d'une courbe g\'eod\'esique $\gamma$, et on va  associer une  fl\`eche \`a chacun des \emph{sommets}, c'est-\`a-dire des extr\'emit\'es des morceaux de la d\'ecomposition. Appelons $\gamma_{1}$ et $\gamma_{2}$  les deux morceaux situ\'es de part et d'autre du sommet $x_{1}$. Par minimalit\'e, l'image de $\gamma_{1}$ par $h$ rencontre $\gamma_{2}$, ou l'image de $\gamma_{2}$  rencontre $\gamma_{1}$. D'autre part,  un lemme de Franks dit que \emph{ces deux possibilit\'es ne peuvent pas arriver simultan\'ement lorsque l'indice de Poincar\'e-Lefschetz est diff\'erent de $1$}.
La fl\`eche horizontale associ\'ee au sommet $x_{1}$ d\'ecrit alors le sens dans lequel s'effectue la transition.
Il reste \`a construire les fl\`eches ``verticales'' ($\ua$ et  $\da$), qui vont venir s'intercaler entre deux fl\`eches horizontales successives. En reprenant les notations du paragraphe pr\'ec\'edent, on va associer une fl\`eche verticale \`a chaque morceau de la d\'ecomposition~;
au final, les trois fl\`eches associ\'ees au sommet $x_{1}$ et aux deux morceaux qui lui sont adjacents contiendront l'information sur la dynamique de $x_{1}$. Pour fixer les id\'ees, supposons que la transition au sommet $x_{1}$ s'effectue de $\gamma_{1}$ vers $\gamma_{2}$. Alors la fl\`eche associ\'ee \`a $\gamma_{1}$ d\'ecrit 
le pass\'e du point $x_{1}$ (par convention, elle vaut $\da$ si l'orbite de $x_{1}$ vient de l'unique point fixe, et $\ua$ s'il vient de l'infini)~; tandis que la fl\`eche associ\'ee \`a $\gamma_{2}$ d\'ecrit le futur de $x_{1}$. Ainsi, par exemple, le mot $(\da \ra \ua)$ symbolise une orbite qui vient du point fixe et retourne vers le point fixe apr\`es avoir ``effectu\'e un demi-tour sur sa gauche''.

Prouver le th\'eor\`eme \'enonc\'e ci-dessus revient maintenant \`a montrer que le mot $\mathbf{M}(h)$ ainsi construit ne d\'epend pas de la courbe $\gamma$ utilis\'ee dans la construction (parmi toutes les courbes r\'ealisant le module de $h$). Ceci sera fait \textit{via} la construction de structures dynamiques canoniques qui g\'en\'eralisent les secteurs  des points hyperboliques de type selle ou des points paraboliques en dynamique holomorphe. Par analogie au cas des feuilletages (ou des flots), ces structures sont baptis\'ees \emph{composantes de Reeb}. Notons que dans la section~\ref{s.definitions}, le point de vue est renvers\'e : nous d\'efinissons le mot $\mathbf{M}(h)$ \`a partir de la dynamique dans les bords des composantes de Reeb canoniques, puis nous montrerons l'\'equivalence avec la d\'efinition esquiss\'ee ci-dessus en termes de courbes g\'eod\'esiques.

%%%%%%%%%%%%%%%%%%%%%%%
\subsection{Structures dynamiques, sur un exemple}
\label{ss.exemple}
L'invariant permet de lire certaines caract\'eristiques de la dynamique de $h$,
en pr\'ecisant la mani\`ere dont les quatre types d'orbites sont organis\'ees autour du point fixe. \`A chaque mot cyclique $\textbf{M}$ correspond un mod\`ele dynamique $h_{\textbf{M}}$, qui est en quelque sorte la dynamique la plus simple possible parmi les hom\'eomorphismes $h$ tels que $\textbf{M}(h) = \textbf{M}$. On essaie alors de retrouver pour l'hom\'eomorphisme g\'en\'eral $h$ des propri\'et\'es dynamiques du mod\`ele $h_{\textbf{M}(h)}$.
 Les r\'esultats pr\'ecis seront \'enonc\'es dans la partie~\ref{p.consequences-dynamiques} ; nous nous contentons ici d'expliquer certaines cons\'equences dynamiques sur un exemple.
 On se restreint encore aux  hom\'eomorphismes du plan qui fixent uniquement le point $0$. Supposons que l'indice par quarts-de-tour de $h$ soit 
$$\mathbf{M}=(\ua \ra \da \ra \ua \ra \da \la \da \la ).$$
 Nous commen\c cons par dessiner la dynamique du mod\`ele $h_{\textbf{M}}$ correspondant (figure~\ref{f.exemple-modele}).
Le module est \'egal \`a $5$, il y a
donc $5$ secteurs  $S_{1}, \dots, S_{5}$,
qui correspondent respectivement aux mots partiels $(\ua \ra \da)$, $(\da \ra \ua)$, $(\ua \ra \da)$, 
$(\da \la \da)$, $(\da \la \ua)$.
 Les secteurs $1$, $3$ et $5$ sont hyperboliques, le  deuxi\`eme est
elliptique, le quatri\`eme est indiff\'erent.\footnote{La construction pr\'ecise d'un mod\`ele de secteur  indiff\'erent sera donn\'ee \`a la section~\ref{ss.realisabilite}.}
 Pour retrouver l'indice de Poincar\'e-Lefschetz, 
on additionne les contributions de chacun des secteurs, et on ajoute $1$. Les secteurs hyperboliques comptent pour $-1/2$, les secteurs elliptiques pour $+1/2$, les secteurs indiff\'erents comptent pour du beurre. Par cons\'equent l'indice de Poincar\'e-Lefschetz
du point fixe est ici \'egal \`a $(-1/2+1/2-1/2+0-1/2) +1 = 0$.

\begin{figure}[htbp]
\begin{center}
%Num\'eroter les secteurs.
\def\JPicScale{1}
\ifx\JPicScale\undefined\def\JPicScale{1}\fi
\psset{unit=\JPicScale mm}
\psset{linewidth=0.3,dotsep=1,hatchwidth=0.3,hatchsep=1.5,shadowsize=1}
\psset{dotsize=0.7 2.5,dotscale=1 1,fillcolor=black}
\psset{arrowsize=1 2,arrowlength=1,arrowinset=0.25,tbarsize=0.7 5,bracketlength=0.15,rbracketlength=0.15}
\begin{pspicture}(0,0)(45,40)
\psbezier(28.11,17.96)(28.11,17.96)(28.11,17.96)(28.11,17.96)
\rput{0}(22.7,19.86){\psellipse[linewidth=0.2,linestyle=dashed,dash=1 1](0,0)(20.04,20.04)}
\psbezier(22.79,19.97)(22.79,19.97)(22.79,19.97)(22.79,19.97)
\psdots[](22.79,19.97)
(22.79,19.97)
\psline{->}(22.89,-0.18)(22.82,2.53)
\psline{->}(23.03,39.9)(19.98,40)
\psline{->}(41.55,12.94)(45,11.4)
\psline{->}(33.88,3.16)(37.04,5.39)
\psline{->}(11.31,36.08)(9.32,38.85)
\psline{->}(35.65,35.11)(33.62,32.42)
\psline{->}(42.12,24.54)(41.44,28.56)
\psline{->}(3.99,26.42)(5.33,29.85)
\psline{->}(3.55,15.09)(0,14.05)
\psline{->}(10.96,3.76)(8.32,5.68)
\psline{<-}(34.65,15.43)(22.79,19.91)
\psline{->}(22.79,19.91)(15.14,30.13)
\psline{<-}(9.9,16.39)(22.79,19.91)
\psline(22.79,19.91)(22.64,6.62)
\psline(22.79,19.91)(31.64,30.06)
\pscustom[]{\psbezier{-}(21.49,6.66)(21.6,12.9)(21.6,16.04)(20.7,16.77)
\psbezier{->}(20.7,16.77)(19.8,17.5)(17.41,17.36)(10.61,15.21)
}
\pscustom[]{\psbezier{-}(24.07,6.6)(24.24,12.3)(23.66,16.31)(25.06,16.96)
\psbezier{->}(25.06,16.96)(26.48,17.6)(29.35,15.87)(34.11,14.33)
}
\pscustom[]{\psbezier{<-}(16.33,30.86)(18.85,27.16)(21.63,24.15)(23.12,23.97)
\psbezier{-}(23.12,23.97)(24.63,23.8)(26.97,26.97)(30.43,30.64)
}
\pscustom[linestyle=dashed,dash=1 1]{\psbezier{-}(22.34,20.02)(19.19,19.58)(14.18,17.37)(13.23,18.98)
\psbezier(13.23,18.98)(12.3,20.58)(16.91,19.73)(17.68,20.98)
\psbezier{->}(17.68,20.98)(18.45,22.25)(19.11,24.11)(14.73,28.84)
}
\psbezier{->}(22.79,19.91)(28.19,18)(36.85,15.39)(34.72,22.61)
\psbezier(34.7,22.72)(32.01,29.49)(26.44,23.77)(22.95,19.98)
\psbezier{->}(22.83,20.03)(26.54,18.72)(32.48,16.92)(31,21.85)
\psbezier(31,21.93)(29.14,26.55)(25.33,22.66)(22.93,20.08)
\psbezier{->}(18.83,7.92)(19.69,12.31)(17.96,14.67)(12.89,12.91)
\psbezier{<-}(19.24,31.33)(21.83,27.47)(24.8,27.4)(28.14,31.5)
\psbezier{->}(26.86,7.78)(26.5,12.75)(28.12,13.97)(32.38,12.5)
\psline{->}(22.69,11.73)(22.71,13.26)
\rput(32.42,6.81){$S_1$}
\rput(31.11,5.59){}
\rput(15.27,6.43){$S_5$}
\rput(39.32,23.02){$S_2$}
\rput(9.86,23.87){$S_4$}
\rput(24.12,34.74){$S_3$}
\psline{->}(31.7,30.14)(28.01,25.84)
\end{pspicture}
\caption{}
\label{f.exemple-modele}
\end{center}
\end{figure}

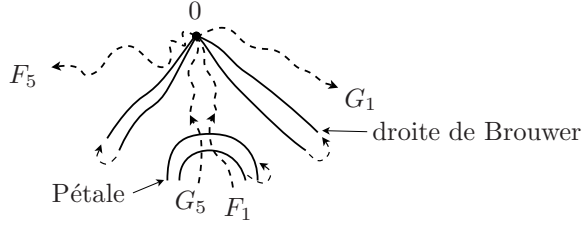
\begin{figure}[htbp]
\begin{center}
\def\JPicScale{0.8}
\ifx\JPicScale\undefined\def\JPicScale{1}\fi
\psset{unit=\JPicScale mm}
\psset{linewidth=0.3,dotsep=1,hatchwidth=0.3,hatchsep=1.5,shadowsize=1}
\psset{dotsize=0.7 2.5,dotscale=1 1,fillcolor=black}
\psset{arrowsize=1 2,arrowlength=1,arrowinset=0.25,tbarsize=0.7 5,bracketlength=0.15,rbracketlength=0.15}
\begin{pspicture}(0,0)(65,35)
\rput(41.49,12.92){}
\psdots[](30.89,30.71)
(30.89,30.71)
\pscustom[linestyle=dashed,dash=1 1]{\psbezier{-}(30.71,30.8)(27.68,32.32)(28.66,30.89)(29.01,29.46)
\psbezier(29.01,29.46)(29.37,28.04)(25.8,31.7)(26.51,28.75)
\psbezier(26.51,28.75)(27.23,25.8)(25.44,24.55)(23.3,26.43)
\psbezier(23.3,26.43)(21.16,28.3)(20.53,29.73)(18.12,28.75)
\psbezier(18.12,28.75)(15.71,27.77)(15.89,25.27)(13.48,24.38)
\psbezier{->}(13.48,24.38)(11.07,23.48)(13.12,26.34)(6.78,25.54)
}
\pscustom[linestyle=dashed,dash=1 1]{\psbezier(30.71,30.8)(32.5,26.79)(28.57,26.16)(29.82,23.57)
\psbezier(29.82,23.57)(31.07,20.98)(30.62,21.88)(30,18.84)
\psbezier(30,18.84)(29.37,15.8)(31.43,16.79)(31.69,13.13)
\psbezier(31.69,13.13)(31.96,9.46)(31.25,6.25)(31.25,6.25)
}
\pscustom[linestyle=dashed,dash=1 1]{\psbezier(30.89,30.63)(33.66,26.79)(30.98,27.14)(33.12,24.29)
\psbezier(33.12,24.29)(35.26,21.43)(33.3,20.54)(34.01,18.93)
\psbezier(34.01,18.93)(34.73,17.32)(32.76,17.95)(33.12,14.11)
\psbezier(33.12,14.11)(33.48,10.27)(35.35,10.63)(36.07,8.93)
\psbezier(36.07,8.93)(36.78,7.23)(36.96,5.54)(36.96,5.54)
}
\pscustom[linestyle=dashed,dash=1 1]{\psbezier{-}(30.89,30.54)(33.84,32.14)(35.8,33.04)(37.85,30.63)
\psbezier(37.85,30.63)(39.91,28.21)(41.78,28.75)(43.39,28.48)
\psbezier(43.39,28.48)(45,28.21)(46.07,26.96)(47.23,25.89)
\psbezier(47.23,25.89)(48.39,24.82)(49.1,25.36)(50.8,23.93)
\psbezier{->}(50.8,23.93)(52.5,22.5)(52.05,23.04)(54.46,22.14)
}
\pscustom[]{\psbezier(30.89,30.71)(25.98,21.88)(26.07,20.71)(23.07,18.71)
\psbezier(23.07,18.71)(20.07,16.71)(18.07,11.71)(18.07,11.71)
}
\pscustom[]{\psbezier(30.71,30.8)(26.43,24.91)(25.07,22.71)(22.07,20.71)
\psbezier(22.07,20.71)(19.07,18.71)(16.07,13.71)(16.07,13.71)
}
\psline{->}(30.26,16.7)(29.82,17.59)
\psline{->}(33.57,17.14)(33.97,17.9)
\pscustom[]{\psbezier(31.07,30.71)(36.07,28.71)(36.07,26.71)(39.07,24.71)
\psbezier(39.07,24.71)(42.07,22.71)(48.07,17.71)(51.07,14.71)
}
\pscustom[]{\psbezier(31.07,30.71)(33.07,28.71)(34.07,25.71)(38.07,21.71)
\psbezier(38.07,21.71)(42.07,17.71)(46.07,13.71)(49.07,11.71)
}
\pscustom[]{\psbezier(28.07,6.71)(28.07,10.71)(31.07,11.71)(33.07,11.71)
\psbezier(33.07,11.71)(35.07,11.71)(38.07,10.71)(39.07,6.71)
}
\pscustom[]{\psbezier(26.07,6.71)(26.07,10.71)(27.07,12.71)(29.07,13.71)
\psbezier(29.07,13.71)(31.07,14.71)(35.07,14.71)(37.07,13.71)
\psbezier(37.07,13.71)(39.07,12.71)(41.07,10.71)(41.07,6.71)
}
\pscustom[linewidth=0.2,linestyle=dashed,dash=1 1]{\psbezier{-}(39.37,7.59)(40.76,6.55)(42.14,5.63)(42.94,6.88)
\psbezier{->}(42.94,6.88)(43.75,8.13)(42.76,8.75)(41.34,9.64)
}
\pscustom[linewidth=0.2,linestyle=dashed,dash=1 1]{\psbezier{-}(49.64,11.07)(51.28,9.36)(51.78,10)(52.59,10.71)
\psbezier{->}(52.59,10.71)(53.39,11.43)(52.76,12.59)(51.78,13.93)
}
\pscustom[linewidth=0.2,linestyle=dashed,dash=1 1]{\psbezier{-}(17.59,10.63)(16.6,9.64)(16.25,9.02)(15.07,9.71)
\psbezier{->}(15.07,9.71)(13.89,10.41)(14.46,10.98)(15.07,12.71)
}
\rput(30.71,35){$0$}
\rput(31.25,1.61){}
\rput(37.85,1.96){$F_1$}
\rput(58.21,20.18){$G_1$}
\rput(30,3.21){$G_5$}
\rput(2.14,24.11){$F_5$}
\rput(36.25,3.57){}
\psline[linewidth=0.15]{->}(59.28,14.82)(51.78,15)
\psline[linewidth=0.15]{->}(20.35,4.64)(25.53,6.79)
\rput[l](60.18,15){droite de Brouwer}
\rput[r](18.93,4.64){P\'etale}
\rput[l](65,16.43){}
\end{pspicture}
\caption{Quelques propri\'et\'es dynamiques h\'erit\'ees du mod\`ele sur les secteurs $S_{1}$ et $S_{5}$}
\label{f.exemple-vrai}
\end{center}
\end{figure}

Pour imaginer un hom\'eomorphisme g\'en\'eral $h$ ayant le m\^eme indice par quarts-de-tour, il faut complexifier le mod\`ele $h_{\textbf{M}}$ de la fa\c con suivante. On commence par \'epaissir chaque s\'eparatrice en un secteur de dynamique parabolique (orbites ``parall\`eles'').
On complexifie la topologie des bords des secteurs : par exemple, on peut obtenir des bords qui ne sont plus localement connexes. Enfin, on peut complexifier \`a l'infini la dynamique au voisinage de $0$ et de l'infini (voir par exemple les figures 10, 13 et 14 du texte~\cite{0}).
L'hom\'eomorphisme g\'en\'eral $h$ h\'erite d'un certain nombre de propri\'et\'es dynamiques du mod\`ele $h_{\textbf{M}}$. Voici la plus importante (figure~\ref{f.exemple-vrai}).
%La propri\'et\'e dynamique essentielle du mod\`ele qui se transmet \`a $h$ sont les
%suivantes.
 \`A chaque secteur $S_{i}$ du mod\`ele correspond canoniquement un ouvert $O(F_{i}, G_{i})$ invariant par $h$, dont les bords $F_{i}$ et $G_{i}$ forment une \emph{composante de Reeb} (voir la d\'efinition~\ref{d.reeb} pour les d\'etails).
 Les orbites de presque tous les points d'un de ces bords ont une dynamique identique \`a celle du bord correspondant pour le mod\`ele : ainsi, presque tous les points de $F_{1}$ vont de l'infini \`a $0$, et presque tous les points de $G_{1}$ vont de $0$ \`a l'infini ; ceci g\'en\'eralise l'existence de s\'eparatrices stables ou instables (proposition~\ref{p.dynamique-bord}). Entre les ensembles $F_{i}$ et $G_{i}$, la dynamique va globalement dans le m\^eme sens orthoradial que la dynamique du mod\`ele dans le secteur $S_{i}$ : par exemple, il existe une demi-droite topologique qui joint $0$ \`a l'infini, situ\'ee entre $F_{1}$ et $G_{1}$, et qui est pouss\'ee par $h$ vers $G_{1}$ (c'est ce qu'on appelle une \emph{droite de Brouwer}). De plus, tout ouvert rencontrant $F_{1}$ a ses it\'er\'es positifs  qui s'accumulent sur $G_{1}$.

Voici quelques propri\'et\'es suppl\'ementaires.
Premi\`erement, dans chaque composante de Reeb, il existe un ouvert, d\'efini canoniquement, dont les points ont tous le m\^eme comportement dynamique, et ce comportement se lit sur le mot $\textbf{M}$~: par exemple, l'ouvert canonique entre $F_{2}$ et $G_{2}$ est constitu\'e de points dont les orbites vont toutes de $0$ \`a $0$, comme dans le secteur $S_{2}$ du mod\`ele. On en d\'eduira en particulier que si $h$ pr\'eserve la mesure de Lebesgue, alors le mod\`ele $h_{\textbf{M(h)}}$ n'a que des secteurs hyperboliques ; le mot $\textbf{M(h)}$ est alors obtenu par r\'ep\'etition du mot \'el\'ementaire $(\ua \la \da \ra)$, et on retrouvera naturellement un th\'eor\`eme de Slaminka : l'indice de Poincar\'e-Lefschetz d'un hom\'eomorphisme conservatif $h$ est inf\'erieur ou \'egal \`a $1$. 
Deuxi\`emement, on peut lire sur le mot $\textbf{M(h)}$ l'existence de \emph{p\'etales}. L'un d'eux, associ\'e au sous-mot $(\ra \da \la)$ de $\textbf{M(h)}$, est repr\'esent\'e sur la figure~\ref{f.exemple-vrai} ; il s'agit d'un p\'etale r\'epulsif bas\'e en l'infini.
Cette propri\'et\'e  g\'en\'eralise un \'enonc\'e du texte~\cite{0} (version topologique du th\'eor\`eme de la fleur de Leau-Fatou, p9).
Enfin, chaque secteur indiff\'erent  contient au moins une orbite dont la dynamique est oppos\'ee \`a celle des bords : ainsi, dans le secteur correspondant \`a  $S_{4}$, il existe une orbite qui va de l'infini \`a $0$.

%-- dans tout secteur indiff\'erent, il existe au moins une orbite dont la dynamique est oppos\'ee \`a la dynamique aux bords : par exemple, $h$ admet, dans le secteur  entre $F_{4}$ et $G_{4}$, une orbite qui va de l'infini \`a $0$.

%-- supposons que l'indice par quart-de-tour indique l'existence de deux secteurs hyperboliques ou deux secteurs elliptiques adjacents. Alors il existe des p\'etales pour $h$ : par exemple, les secteurs 1 et 5 \'etant hyperboliques, on peut trouver une droite topologique qui va de l'infini \`a l'infini, qui est incluse dans le secteur d\'elimit\'e par $F_{5}$ et $G_{1}$, et qui borde un attracteur pour $h$. De plus, il existe une droite de ce type qui passe   arbitrairement pr\`es de $0$, mais il en existe une autre qui est contenue dans un voisinage arbitrairement petit de l'infini.

%
%-- ouvert de points qui sont OK  ?

%%%%%%%%%%%%%%%%%%%%%%%
\subsection{Invariant global \emph{vs} invariant local}
\begin{figure}[htbp]
\begin{center}
\def\JPicScale{0.8}
\ifx\JPicScale\undefined\def\JPicScale{1}\fi
\psset{unit=\JPicScale mm}
\psset{linewidth=0.3,dotsep=1,hatchwidth=0.3,hatchsep=1.5,shadowsize=1}
\psset{dotsize=0.7 2.5,dotscale=1 1,fillcolor=black}
\psset{arrowsize=1 2,arrowlength=1,arrowinset=0.25,tbarsize=0.7 5,bracketlength=0.15,rbracketlength=0.15}
\begin{pspicture}(0,0)(120,40)
\psline(30,32.5)(30,7.5)
\psline(10,20)(50,20)
\psline{->}(30,40)(30,30)
\psline{->}(30,20)(20,20)
\psline{->}(30,-0)(30,10)
\psline{->}(30,20)(40,20)
\psline(100,32.5)(100,7.5)
\psline(87.5,20)(120,20)
\psbezier(87.5,22.5)(95,22.5)(97.5,25)(97.5,32.5)
\psbezier(87.5,21.25)(98.75,21.25)(98.75,21.25)(98.75,32.5)
\psbezier(87.5,23.75)(91.25,23.75)(96.25,28.75)(96.25,32.5)
\psbezier(112.5,17.5)(105,17.5)(102.5,15)(102.5,7.5)
\psbezier(112.5,18.75)(101.25,18.75)(101.25,18.75)(101.25,7.5)
\psbezier(112.5,16.25)(108.75,16.25)(103.75,11.25)(103.75,7.5)
\psline{->}(100,40)(100,30)
\psline{->}(100,20)(90,20)
\psline{->}(100,6.25)(100,10)
\psline{->}(100,20)(110,20)
\psbezier{->}(100,20)(83.75,20)(83.75,13.75)(88.75,8.75)
\psbezier{->}(100,20)(88.75,20)(87.5,15)(91.25,11.25)
\psbezier{->}(100,20)(93.75,20)(90,17.5)(93.75,13.75)
\psbezier(100,20)(100,3.75)(93.75,3.75)(88.75,8.75)
\psbezier(100,20)(100,8.75)(95,7.5)(91.25,11.25)
\psbezier(100,20)(100,13.75)(97.5,10)(93.75,13.75)
\psbezier(42.5,17.5)(35,17.5)(32.5,15)(32.5,7.5)
\psbezier(42.5,18.75)(31.25,18.75)(31.25,18.75)(31.25,7.5)
\psbezier(42.5,16.25)(38.75,16.25)(33.75,11.25)(33.75,7.5)
\psline(33.75,7.5)(33.75,-0)
\psline(32.5,7.5)(32.5,-0)
\psline(31.25,7.5)(31.25,-0)
\psline(42.5,18.75)(50,18.75)
\psline(42.5,17.5)(50,17.5)
\psline(42.5,16.25)(50,16.25)
\psbezier(36.25,-0)(36.25,9.62)(40.38,13.75)(50,13.75)
\psbezier(41.25,-0)(41.25,6.12)(43.88,8.75)(50,8.75)
\psbezier(17.5,22.5)(25,22.5)(27.5,25)(27.5,32.5)
\psbezier(17.5,21.25)(28.75,21.25)(28.75,21.25)(28.75,32.5)
\psbezier(17.5,23.75)(21.25,23.75)(26.25,28.75)(26.25,32.5)
\psline(26.25,32.5)(26.25,40)
\psline(27.5,32.5)(27.5,40)
\psline(28.75,32.5)(28.75,40)
\psline(17.5,21.25)(10,21.25)
\psline(17.5,22.5)(10,22.5)
\psline(17.5,23.75)(10,23.75)
\psbezier(23.75,40)(23.75,30.37)(19.62,26.25)(10,26.25)
\psbezier(18.75,40)(18.75,33.87)(16.12,31.25)(10,31.25)
\psbezier(32.5,32.5)(32.5,25)(35,22.5)(42.5,22.5)
\psbezier(31.25,32.5)(31.25,21.25)(31.25,21.25)(42.5,21.25)
\psbezier(33.75,32.5)(33.75,28.75)(38.75,23.75)(42.5,23.75)
\psline(42.5,23.75)(50,23.75)
\psline(42.5,22.5)(50,22.5)
\psline(42.5,21.25)(50,21.25)
\psline(31.25,32.5)(31.25,40)
\psline(32.5,32.5)(32.5,40)
\psline(33.75,32.5)(33.75,40)
\psbezier(50,26.25)(40.38,26.25)(36.25,30.37)(36.25,40)
\psbezier(50,31.25)(43.88,31.25)(41.25,33.87)(41.25,40)
\psbezier{->}(30,20)(13.75,20)(13.75,13.75)(18.75,8.75)
\psbezier{->}(30,20)(18.75,20)(17.5,15)(21.25,11.25)
\psbezier{->}(30,20)(23.75,20)(20,17.5)(23.75,13.75)
\psbezier(30,20)(30,3.75)(23.75,3.75)(18.75,8.75)
\psbezier(30,20)(30,8.75)(25,7.5)(21.25,11.25)
\psbezier(30,20)(30,13.75)(27.5,10)(23.75,13.75)
\psbezier(30,20)(30,-0)(22.5,-0)(16.25,6.25)
\psbezier(30,20)(30,-2.5)(21.25,-3.75)(13.75,3.75)
\psbezier{->}(30,20)(7.5,20)(6.25,10)(13.75,3.75)
\psbezier{<-}(16.25,6.25)(10,12.5)(10,20)(30,20)
\pscustom[]{\psbezier(87.5,23.75)(82.5,23.75)(80,21.25)(80,15)
\psbezier(80,15)(80,8.75)(80,5)(82.5,2.5)
\psbezier(82.5,2.5)(85,-0)(88.75,-0)(95,-0)
\psbezier(95,-0)(101.25,-0)(103.75,2.5)(103.75,7.5)
}
\pscustom[]{\psbezier(87.5,21.25)(83.75,21.25)(82.5,17.5)(82.5,15)
\psbezier(82.5,15)(82.5,12.5)(82.5,7.5)(85,5)
\psbezier(85,5)(87.5,2.5)(92.5,2.5)(95,2.5)
\psbezier(95,2.5)(97.5,2.5)(101.25,3.75)(101.25,7.5)
\psbezier(101.25,7.5)(101.25,11.25)(101.25,7.5)(101.25,7.5)
}
\pscustom[]{\psbezier(87.5,20)(83.75,20)(83.75,16.25)(83.75,13.75)
\psbezier(83.75,13.75)(83.75,11.25)(83.75,8.75)(86.25,6.25)
\psbezier(86.25,6.25)(88.75,3.75)(91.25,3.75)(93.75,3.75)
\psbezier(93.75,3.75)(96.25,3.75)(100,3.75)(100,7.5)
}
\pscustom[]{\psbezier(87.5,22.5)(82.5,22.5)(81.25,18.75)(81.25,13.75)
\psbezier(81.25,13.75)(81.25,8.75)(81.25,6.25)(83.75,3.75)
\psbezier(83.75,3.75)(86.25,1.25)(88.75,1.25)(93.75,1.25)
\psbezier(93.75,1.25)(98.75,1.25)(102.5,2.5)(102.5,7.5)
}
\psline(112.5,16.25)(120,16.25)
\psline(112.5,17.5)(120,17.5)
\psline(112.5,18.75)(120,18.75)
\psline(98.75,32.5)(98.75,40)
\psline(97.5,40)(97.5,32.5)
\psline(96.25,32.5)(96.25,40)
\pscustom[]{\psbezier(112.5,13.75)(106.25,13.75)(106.25,7.5)(106.25,5)
\psbezier(106.25,5)(106.25,2.5)(105,-0)(105,-0)
}
\psline(112.5,13.75)(120,13.75)
\psbezier(120,10)(112.5,10)(110,7.5)(110,-0)
\psbezier(115,-0)(115,3.75)(116.25,5)(120,5)
\pscustom[]{\psbezier(93.75,32.5)(93.75,26.25)(87.5,26.25)(85,26.25)
\psbezier(85,26.25)(82.5,26.25)(80,25)(80,25)
}
\psline(93.75,32.5)(93.75,40)
\psbezier(90,40)(90,32.5)(87.5,30)(80,30)
\psbezier(80,35)(83.75,35)(85,36.25)(85,40)
\psbezier(101.25,32.5)(101.25,25)(103.75,22.5)(111.25,22.5)
\psbezier(100,32.5)(100,21.25)(100,21.25)(111.25,21.25)
\psbezier(102.5,32.5)(102.5,28.75)(107.5,23.75)(111.25,23.75)
\psline(111.25,23.75)(118.75,23.75)
\psline(111.25,22.5)(118.75,22.5)
\psline(111.25,21.25)(118.75,21.25)
\psline(100,32.5)(100,40)
\psline(101.25,32.5)(101.25,40)
\psline(102.5,32.5)(102.5,40)
\psbezier(118.75,26.25)(109.13,26.25)(105,30.38)(105,40)
\psbezier(118.75,31.25)(112.63,31.25)(110,33.87)(110,40)
\end{pspicture}
\caption{Un germe donn\'e peut s'\'etendre de plusieurs mani\`eres, donnant diff\'erents indices par quarts-de-tour : $(\ua \ra \da \la \ua \ra \da \ra)$ et $(\ua \ra \da \la)$ }
\label{f.pas-local}
\end{center}
\end{figure}
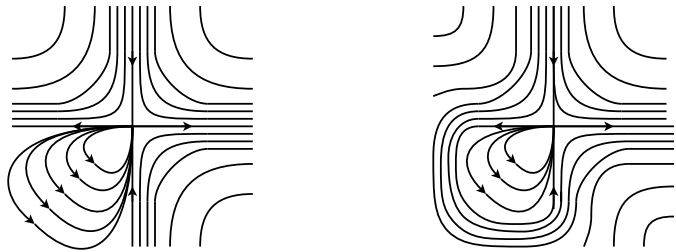

La figure~\ref{f.pas-local} montre deux hom\'eomorphismes du plan, fixant uniquement le point $0$, avec un indice \'egal \`a $0$, qui ont le m\^eme germe au voisinage du point fixe. Cependant, leurs  indices par quarts-de-tour diff\`erent. Ceci montre que cet invariant n'est pas un invariant de la dynamique locale.
Il est cependant possible qu'on puisse obtenir un invariant local similaire \`a l'indice par quart-de-tours~; en particulier, les r\'ecents travaux de 
F. R. Ruiz Del Portal et J. M. Salazar vont dans ce sens (voir~\cite{rdps}).

%%%%%%%%%%%%%%%%%%%%%%%
\subsection{Changement de cadre}
Dans cette section, nous expliquons  bri\`evement un proc\'ed\'e canonique qui produit, \`a partir des donn\'ees g\'en\'erales de l'introduction, un hom\'eomorphisme du plan dont $0$ est l'unique point fixe, avec un indice diff\'erent de $1$.  Ce proc\'ed\'e permet de se ramener du cadre g\'en\'eral \`a ce cadre plus simple : comme le proc\'ed\'e est canonique, toute construction d'un invariant de conjugaison dans le cadre simple s'\'etendra automatiquement au cadre g\'en\'eral.
Nous nous contentons de donner les \'enonc\'es, la preuve  \'etant identique \`a la preuve du th\'eor\`eme d'extension 2.1 du texte~\cite{0}.

%En fait, ce proc\'ed\'e est d\'ej\`a partiellement d\'ecrit au chapitre 2 du texte~\cite{0} (p 25), nous en rappellerons seulement les grandes lignes.

On se donne  un hom\'eomorphisme $h_0$  d'une surface orientable $S$, et  un point fixe isol\'e $x_{0}$ (sans condition d'indice  pour le moment).
Comme indiqu\'e plus haut,
on exclut le cas des hom\'eomorphismes de Brouwer, o\`u $S$ est la sph\`ere et $x_{0}$ l'unique point fixe de $h_0$.
 On consid\`ere la surface (non compacte) $O$ obtenue en enlevant \`a $S$ les points fixes de $h_{0}$. La premi\`ere proposition  est purement topologique, elle  donne une  correspondance entre le plan et la surface $O$.
\begin{prop}
\label{p.cadre1}
Il existe une unique application $p : \bbR^2 \rightarrow O \cup \{x_{0}\}$ telle que
\begin{enumerate}
\item $p(0) = x_{0}$ ;
\item  il existe un voisinage $U$ de $0$ dans le plan sur lequel la restriction de $p$ est un hom\'eomorphisme de $U$ vers un voisinage $V$ de $x_{0}$ dans la surface $S$ ;
\item la restriction de $p$ au plan trou\'e $\bbR^2 \setminus \{ 0\}$ est un rev\^etement
sur $O$.
\end{enumerate}
\end{prop}
L'unicit\'e s'entend \`a isomorphisme de rev\^etement pr\`es (si $p_{1}$ et $p_{2}$ sont deux telles applications, alors il existe un hom\'eomorphisme $g$ du plan, fixant $0$, pr\'eservant l'orientation, tel que $p_{2} = g \circ p_{1}$).
On peut ensuite relever la dynamique.
\begin{prop}
\label{p.cadre2}
Il existe un unique hom\'eomorphisme $h$ du plan, fixant $0$, qui v\'erifie
$$
p \circ h = h_{0} \circ p.
$$
En particulier, 
\begin{enumerate}
\item le point $0$ est l'unique point fixe de $h$ ;
\item l'hom\'eomorphisme $h$, au voisinage de $0$,  est localement conjugu\'e \`a l'hom\'eomorphisme $h_0$ au voisinage de $x_{0}$ ;
\item si $h_0$ pr\'eserve l'orientation, alors $h$ aussi ; 
\item l'indice  de Poincar\'e-Lefschetz de $0$ pour $h$ est \'egal \`a l'indice 
 de Poincar\'e-Lefschetz de $x_{0}$ pour $h_0$.
\end{enumerate}
\end{prop}
%Remarquons enfin que cette conjugaison locale pourrait permettre  de retrouver dans la dynamique de $h_0$ certaines des propri\'et\'es dynamiques de $h$ induites par la connaissance de l'indice par quarts-de-tour (existence d'ensembles stables et instables locaux, par exemple).\footnote{A ENLEVER si pas de r\'esultat local ???}

%%%%%%%%%%%%%%%%%%%%%%%%%%%%%%%%%%%
%%%%%%%%%%%%%%%%%%%%%%%%%%%%%%%%%
\section{D\'efinitions, r\'esultats}
\label{s.definitions}
Pour toute cette section, 
 on se donne un hom\'eomorphisme $h$ du plan,
pr\'eservant l'orientation, fixant le point $0$, sans autre point fixe, et  tel que l'indice de
Poincar\'e-Lefschetz du point fixe est diff\'erent de $1$. On utilisera les r\'esultats fondamentaux de la th\'eorie de Brouwer~: voir ~\cite{0}, section~3.2, et les r\'ef\'erences qui y sont donn\'ees. Tous les \'enonc\'es se g\'en\'eralisent en indice $1$ si le module est assez grand ; cette g\'en\'eralisation est expliqu\'ee dans l'appendice~\ref{as.indice-1}.

%%%%%%%%%%%%
\subsection{Module}
Comme expliqu\'e en introduction, on peut d\'efinir le module de $h$ par des consid\'erations sur les courbes qui entourent le point fixe. Nous reviendrons sur cette d\'efinition \`a la section~\ref{ss.definition-courbes}~; nous commen\c cons par donner une d\'efinition en termes de domaines de translation.
Un \emph{domaine de translation} est un ouvert du plan, connexe et
simplement connexe, invariant
par $h$, sur lequel la restriction de $h$ est conjugu\'ee \`a une translation.
Le \emph{th\'eor\`eme des translations planes} de Brouwer, adapt\'e \`a notre
cadre avec un point fixe, dit notamment que \emph{tout point du plan, autre que le point
fixe $0$, est dans un domaine de translation} (voir Slaminka \cite{slam}, Guillou ~\cite{guil},  ou encore~\cite{0}, th\'eor\`eme 3.79 p71).
\begin{defi}
\label{d.module}
Le \emph{module} de $h$ est le nombre minimal de domaines de
translation dont on a besoin pour faire le tour du point fixe~:
$$
\diam(h):=\min\{k \mid  \exists O_1, \dots O_k \mbox{ domaines de
translation}, O_1 \cup \cdots \cup O_k \mbox{ s\'epare } 0 \mbox{ et } \infty\}.
$$
\end{defi}

%%%%%%%%%%%%%%%%%%%%%%%%%%%%%%%%%%
\subsection{Composantes de Reeb}
\label{ss.composantes-de-reeb}
%%%%%%%%%%%%
\paragraph{D\'efinition des composantes}
On consid\`ere la sph\`ere $\bbS^2=\bbR^2 \cup \{\infty\}$, compactifi\'ee
d'Alexandroff du plan.
Pour toute partie $F$ de $\bbR^2 \setminus \{0\}$, on notera $\hat F
:= F \cup \{0, \infty\}$.
Soient $F$ et $G$ deux parties  ferm\'ees et  connexes dans $\bbR^2
 \setminus\{0\}$, disjointes. On suppose que $\hat F$ et $\hat G$ sont encore connexes (autrement dit, $F$ et $G$
rencontrent tout voisinage de $0$ et tout voisinage de $\infty$).
 Il existe alors exactement deux composantes connexes de
$\bbR^2 \setminus (F \cup G \cup \{0\})$ dont l'adh\'erence rencontre \`a la fois $F$
et $G$. On
les note $O(F,G)$ et $O(G,F)$, de fa\c{c}on \`a ce que l'ordre cyclique
autour de $0$ soit donn\'e par
$$
F < O(F,G) < G < O(G,F) < F.
$$
L'appendice~\ref{as.ordre-cyclique} pr\'ecise l'existence des ensembles $O(F,G)$, $O(G,F)$, ainsi que la d\'efinition de l'ordre cyclique.

\begin{figure}[htbp]
\begin{center}
\def\JPicScale{0.8}
\input{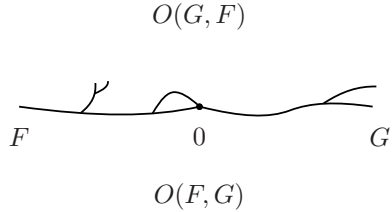}
\caption{Ordre cyclique autour du point fixe}
\label{default}
\end{center}
\end{figure}

\begin{defi}
Soit $U$ un ouvert du plan invariant ($h(U)=U$),
 $x$ et $y$ deux points du plan qui ne sont
pas dans la m\^eme orbite de $h$. On dit que le couple $(x,y)$ est
\emph{singulier via l'ouvert $U$} si pour tous $V_x, V_y$ voisinages
respectifs de $x$ et $y$, il existe un entier positif $n$ tel que 
$$
h^n(V_x) \cap V_y \cap U \neq \emptyset.
$$
%On notera $\sing_U(h)$ l'ensemble des couples singuliers \textit{via} $U$.
\end{defi}
Notons qu'un couple $(x,y)$ est singulier pour l'hom\'eomorphisme $h^{-1}$ si et seulement si $(y,x)$ est singulier pour $h$.

\begin{defi}\label{d.reeb}
Un couple $(F,G)$ de parties de $\bbR^2 \setminus \{0\}$ est une \emph{composante de Reeb}
pour $h$ si
\begin{itemize}
\item $F$ et $G$ sont  deux parties connexes et ferm\'ees  dans $\bbR^2 \setminus \{0\}$~;
\item $\hat F$ et $\hat G$ sont encore connexes~;
\item $F$ et $G$ sont invariants par $h$ (c'est-\`a-dire $h(F)=F$, $h(G)=G$)~;
\item il existe $\varepsilon=\pm 1$ tel que tout couple $(x,y) \in F \times G$ est singulier \emph{via} l'ouvert $O(F,G)$ pour l'hom\'eomorphisme $h^\varepsilon$.
% $F \times G$ est inclus dans $\sing_{O(F,G)}(h)$,  ou bien dans
% $\sing_{O(F,G)}(h^{-1})$. 
\end{itemize}
Les ensembles $F$ et $G$ seront appel\'es \emph{bords} de la composante
de Reeb. 
On dira que la dynamique de la composante de Reeb $(F,G)$  \emph{va de $F$ vers $G$} 
si $F \times G$  
est singulier pour $h$ ($\varepsilon=1$), et qu'elle  \emph{va de $G$ vers $F$} dans le
cas contraire ($\varepsilon=-1$).  
% $(F,G)$ est \emph{de type Ouest-Est} si  c'est $(F,G)$ qui 
%est inclus dans $\sing_{O(F,G)}$, et \emph{de type Est-Ouest} dans le
%cas contraire.  
\end{defi}

\begin{defi}
Une composante de Reeb $(F,G)$ est dite \emph{minimale}
si $F$ et $G$ sont chacun minimal, pour l'inclusion, dans la
collection des parties $K$ ferm\'ees dans $\bbR^2 \setminus \{0\}$,
connexes, telles que $\hat K$ soit connexe.
Autrement dit, $(F,G)$ est minimale si pour
toute composante de Reeb $(F',G')$ telle que $F' \subset F$ et $G'
\subset G$, on a $F'=F$ et $G'=G$.
\end{defi}

\begin{figure}[htbp]
\begin{center}
\def\JPicScale{0.6}
\ifx\JPicScale\undefined\def\JPicScale{1}\fi
\psset{unit=\JPicScale mm}
\psset{linewidth=0.3,dotsep=1,hatchwidth=0.3,hatchsep=1.5,shadowsize=1}
\psset{dotsize=0.7 2.5,dotscale=1 1,fillcolor=black}
\psset{arrowsize=1 2,arrowlength=1,arrowinset=0.25,tbarsize=0.7 5,bracketlength=0.15,rbracketlength=0.15}
\begin{pspicture}(0,0)(45,50)
\psdots[](22.5,45)
(22.5,45)
\pscustom[linewidth=0.5]{\psbezier(22.5,45)(12.5,35)(12.5,35)(12.5,30)
\psbezier(12.5,30)(12.5,25)(12.5,25)(7.5,20)
\psbezier(7.5,20)(2.5,15)(2.5,5)(2.5,5)
}
\pscustom[linewidth=0.5]{\psbezier(22.5,45)(32.5,35)(32.5,35)(32.5,30)
\psbezier(32.5,30)(32.5,25)(32.5,25)(37.5,20)
\psbezier(37.5,20)(42.5,15)(42.5,5)(42.5,5)
}
\psdots[linewidth=0.5](12.68,33.75)
(12.68,33.75)
\psdots[linewidth=0.5](40.54,15.18)
(40.54,15.18)
\rput(2.5,0){$F$}
\rput(22.5,0){$O(F,G)$}
\rput(42.5,0){$G$}
\rput(7.5,35){$x$}
\rput(45,17.5){$y$}
\rput(22.5,50){$0$}
\psdots[linewidth=0.5,linestyle=dashed,dash=1 1,dotsize=0.7 1.5](15,32.5)
(37.5,15)
\pscustom[linewidth=0.5,linestyle=dashed,dash=1 1]{\psbezier{-}(15,32.5)(15,25)(13.57,25.36)(12.5,20)
\psbezier(12.5,20)(11.43,14.64)(17.5,10)(20,12.5)
\psbezier(20,12.5)(22.5,15)(24.11,18.12)(27.5,20)
\psbezier{->}(27.5,20)(30.89,21.88)(33.93,18.66)(36.43,16.16)
}
\rput(23.3,20.45){$h^n$}
\end{pspicture}
\caption{Une composante de Reeb}
\label{   }
\end{center}
\end{figure}
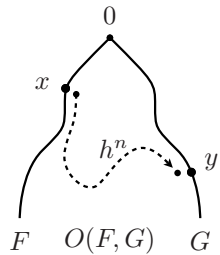

%%%%%%%%%%%%
\subsubsection*{Th\'eor\`eme d'existence de composantes de Reeb}

\begin{theo}
\label{t.compo-Reeb}
Soit $h$ un hom\'eomorphisme du plan,
pr\'eservant l'orientation, fixant le point $0$,  tel que l'indice de
Poincar\'e-Lefschetz du point fixe est diff\'erent de $1$. Soit $d$ le module de $h$.

Alors il existe exactement $d$ composantes de Reeb minimales pour $h$,
que l'on note $(F_1,G_1), \dots , (F_d,G_d)$. 
De plus, 
%les ouverts $O(F_i,G_i), i = 1 \dots d$ sont deux \`a deux
%disjoints, et
 on peut choisir la num\'erotation de fa\c{c}on \`a avoir pour ordre
cyclique
$$
F_1 < O(F_1,G_1) < G_1 \leq F_2 < \cdots < G_d \leq F_1.
$$
\end{theo}
Il est clair que la num\'erotation des composantes de Reeb minimales qui
v\'erifie la condition sur l'ordre cyclique est unique \`a permutation
circulaire pr\`es.\footnote{Ici encore, la notation $\leq$ est pr\'ecis\'ee dans l'appendice~\ref{as.ordre-cyclique}.}
Pour tout entier $k$, les bords $G_k$ et $F_{k+1}$ seront dit
\emph{adjacents}
(les indices $k$ et $k+1$ sont consid\'er\'es modulo $d$~; en
particulier, $G_0=G_d$).
Pour la fin de cette section, on adopte les hypoth\`eses et notations du th\'eor\`eme~\ref{t.compo-Reeb} (notamment pour les \'enonc\'es des propositions~\ref{p.dynamique-bord} et~\ref{p.bords-adjacents} ci-dessous).

%%%%%%%%%%%%
\paragraph{Dynamique dans les bords des composantes de Reeb}
La th\'eorie de Brouwer d\'ecrit le comportement des orbites de $h$ (voir par exemple la section 3.2 de ~\cite{0}, et plus particuli\`erement le corollaire 3.13 p37).
Soit $x$ un point du plan. Vue dans la sph\`ere $\bbS^2=\bbR^2 \cup \{\infty\}$, la suite
des it\'er\'es n\'egatifs de $x$ (son orbite pass\'ee) tend vers $0$ ou l'infini~: on
note $\alpha(x) \in \{0, \infty\}$ le point limite.
 On d\'efinit de m\^eme le point $\omega(x)$, limite de la suite des it\'er\'es positifs du point $x$. 
Soit maintenant $\alpha,\omega \in \{0, \infty\}$~;
on notera 
$$
\wao := \{x \in \bbR^2, \alpha(x)=\alpha \mbox{ et } \omega(x)=\omega \}.
$$
On d\'efinit ainsi quatre ensembles qui forment une partition du plan.

\begin{prop}
\label{p.dynamique-bord}
Soit $F$ un bord  d'une composante de Reeb $(F,G)$ ou $(G,F)$ minimale pour $h$.
Alors :
\begin{itemize}
\item  l'un des deux ensembles  $F  \cap \woi$, $F  \cap \wio$ est vide ;
\item  l'autre est un ouvert   dense de $F$, et il est   connexe ;
\item   les ensembles $F \cap \wii$ et $F \cap \woo$ sont ferm\'es.
\end{itemize}
\end{prop}
%Voici une cons\'equence imm\'ediate. De deux choses l'une~:
%\begin{enumerate}
%\item ou bien $F$ rencontre $\woi$~;
%\item ou bien $F$ rencontre $\wio$. 
%\end{enumerate}
\begin{defi}
On dira que $F$ est \emph{de type dynamique $0 \rightarrow \infty$} si il contient des points qui vont de $0$ \`a l'infini, et 
 \emph{de type dynamique $\infty \rightarrow 0$} dans le cas contraire.
\end{defi}

\begin{prop}
\label{p.bords-adjacents}
%On se place sous les hypoth\`eses du th\'eor\`eme~\ref{t.compo-Reeb}.
Pour tout entier $k$, les deux bords adjacents $G_k$ et $F_{k+1}$ sont du m\^eme type dynamique.
\end{prop}

%%%%%%%%%%%%%%%%%%%%%%%%%%%%%%%%%%
\subsection{D\'efinition dynamique de l'indice}
\label{ss.definition-indice}

%%%%%%%%%
%\subsubsection*{D\'efinition par la dynamique des composantes de Reeb}
Soit $d$ un entier positif.
On appellera \emph{mot cyclique index\'e de longueur $\ell$} une
application d\'efinie sur  $\bbZ$  \`a valeur dans l'alphabet
   $\{\leftarrow,\rightarrow,\uparrow,\downarrow\}$, que l'on notera
$M=(m_i)_{i \in \bbZ}$ et telle que, pour
tout~$i$, $m_{i+\ell}=m_i$. Deux mots cycliques index\'es qui sont
\'egaux \`a d\'ecalage d'indice pr\`es d\'efinissent un \emph{mot
cyclique}~: autrement dit, un mot cyclique $\mathbf{M}$  est une classe
d'\'equivalence de mots cycliques index\'es modulo la relation 
$$
(m_i)_{i \in \bbZ} \sim (m'_i)_{i \in \bbZ} \Leftrightarrow \exists
i_0, \forall i, m_{i+i_0}=m'_i.
$$ 
Un repr\'esentant $M$ de la classe d'\'equivalence du mot $\mathbf{M}$
est une \emph{indexation} de $\mathbf{M}$. 

\begin{defi}
\label{d.mot-cyclique}
Notons $d=\diam(h)$.
On appelle \emph{indice par quarts-de-tour de $h$}
le mot cyclique $\mathbf{M}(h)$ sur l'alphabet
$\{\leftarrow,\rightarrow,\uparrow,\downarrow\}$, de longueur $2d$,
dont une indexation $M=(m_i)_{i \in \bbZ}$ est
d\'efinie de la mani\`ere suivante.
Pour tout entier $k$, 
\begin{itemize}
\item la lettre $m_{2k}$ vaut
	\begin{itemize}
	 \item  $\rightarrow$ si la dynamique de la composante
	 $(F_{k},G_{k})$ va de $F_{k}$ vers $G_{k}$, 
	\item   $\leftarrow$ si la dynamique va de $G_{k}$ vers $F_{k}$.
	\end{itemize}
\item la lettre $m_{2k+1}$ vaut
	\begin{itemize}
	 \item  $\uparrow$ si $G_k$ et $F_{k+1}$ sont de type dynamique $\infty
	\rightarrow 0$,
	\item  $\downarrow$ si $G_k$ et $F_{k+1}$ sont de type dynamique $0
	\rightarrow \infty$.
	\end{itemize}
\end{itemize}
\end{defi}

%Deux mots cycliques sont consid\'er\'es comme \'egaux si on passe de l'un \`a
%l'autre par une permutation circulaire des indices.

Les r\'esultats \'enonc\'es dans la section~\ref{ss.composantes-de-reeb}
impliquent clairement le th\'eor\`eme suivant.
\begin{theo}
\label{t.invariant}
L'indice par quarts-de-tour est un invariant de conjugaison. Plus
pr\'ecis\'ement, soit $g$ un hom\'eomorphisme du plan, pr\'eservant
l'orientation, fixant $0$. Soit $h$  v\'erifiant les hypoth\`eses du th\'eor\`eme~\ref{t.compo-Reeb}.
Alors les mots cycliques $\mathbf{M}(h)$ et $\mathbf{M}(g h g^{-1})$
co\"{\i}ncident.
\end{theo}

%%%%%%%%%%%
\subsection{D\'efinition par les courbes g\'eod\'esiques}
\label{ss.definition-courbes}
%%%%%
\paragraph{Courbes}
Dans ce paragraphe,  $X$ d\'esigne $\bbR^2$ ou bien $\bbR^2\setminus \{0\}$.
Une \emph{courbe} de $X$ est une application continue de $[0,1]$ dans
$X$. Une courbe $\gamma$ est \emph{ferm\'ee} si
$\gamma(1)=\gamma(0)$. On pourra voir une courbe ferm\'ee comme une
application continue de $\bbR$ dans $X$, obtenue en prolongeant
$\gamma$ par p\'eriodicit\'e, ou bien comme une application du cercle
$\bbR/\bbZ$ dans $X$. Une courbe est un \emph{arc} si elle
est injective (autrement dit, sans point double). Une courbe ferm\'ee
est \emph{simple}  si le seul point double
est $\gamma(0)=\gamma(1)$~; on dira aussi que c'est une \emph{courbe de Jordan}.
Pour tous $t_1 < t_2$ dans $[0,1]$, on notera $[t_1  t_2]_\gamma$ (ou
$[\gamma(t_1)  \gamma(t_2)]_\gamma$) la
sous-courbe de $\gamma$ obtenue en restreignant $\gamma$ \`a
l'intervalle $[t_1 t_2]$ (et en reparam\'etrant pour obtenir une
application d\'efinie sur $[0,1]$). Si $\gamma$ est une courbe
ferm\'ee, on la voit comme une application $1$-p\'eriodique de $\bbR$, et
 ceci a encore un sens si $t_1 < t_2$ sont deux r\'eels
quelconque.

Soit $\gamma$ une courbe de Jordan dans $\bbR^2\setminus \{0\}$. Il
existe alors un unique entier relatif $n$ tel que l'application $\gamma :
\bbS^1 \ra \bbR^2\setminus \{0\}$ est homotope \`a
l'application $c_n : \theta \mapsto \exp (2i\pi n \theta)$ (on a
identifi\'e $\bbR^2$ au plan complexe)~; l'entier $n$ est appel\'e
\emph{degr\'e} de $\gamma$ (dans $\bbR^2 \setminus\{0\}$).

On consid\`ere maintenant  une courbe ferm\'ee  $\gamma$. Rappelons qu'un ensemble $E$ est \emph{libre pour l'hom\'eomorphisme $h$} si $h(E)$ est disjoint de $E$. 
\begin{defi}
Une \emph{d\'ecomposition} de $\gamma$ est la donn\'ee de r\'eels 
$t_1 < t_2 < \cdots < t_l < t_{l+1}=t_1+1$ tels que chacune des sous-courbes
$\gamma_i=[t_i t_{i+1}]_\gamma$, $i=1, \dots , l$ soit libre pour $h$.
 De mani\`ere
\'equivalente, on pourra se donner une d\'ecomposition en \'ecrivant $\gamma$ comme la
concat\'enation des sous-courbes $\gamma_i$, $\gamma=\gamma_1 * \cdots
* \gamma_l$. 
L'entier $l$ est appel\'e \emph{longueur} de la d\'ecomposition. Les
points $\gamma(t_i)$ en sont les \emph{sommets}.
\end{defi}
\begin{defi}
La \emph{$h$-longueur} de $\gamma$ est le minimum des longueurs des
d\'ecompositions de $\gamma$. Une d\'ecomposition r\'ealisant ce
minimum est appel\'ee d\'ecomposition \emph{minimale}. 
\end{defi}

%NOTATION (utilis\'ee dans le th\'eo)
%$$
%\gamma=x_0 \frac{\ \gamma_1 \ }{} x_2 \frac{\ \gamma_3 \ }{} \cdots
%\frac{\ \gamma_{2d-1} \ }{} x_{2d}=x_0 .
%$$

%%%%%
\paragraph{Module et courbes}
La proposition suivante fournit une deuxi\`eme d\'efinition possible
du module. Contrairement \`a la premi\`ere d\'efinition, celle-ci ne
n\'ecessite pas l'existence de domaines de translation. Elle permettra donc
de g\'en\'eraliser la d\'efinition du module lorsque l'indice du point fixe est \'egal \`a $1$ (voir section~\ref{ss.jordan-geodesiques} de l'appendice).
\begin{prop}~
\label{p.module-courbes}
\begin{enumerate}
\item Le module de $h$ est \'egal au minimum des $h$-longueurs 
des courbes ferm\'ees $\gamma$ de degr\'e $1$. 
\item  Il existe une courbe de Jordan qui r\'ealise ce minimum.
\end{enumerate}
\end{prop}
%\begin{defi}
%Pour tout \'el\'ement $h$ de ${\cal H}_{0}$, on d\'efinit le module de $h$ comme le 
%minimum des $h$-longueurs 
%des courbes ferm\'ees $\gamma$ de degr\'e $1$.
%\end{defi}
\begin{defi}
Une courbe $\gamma$ r\'ealisant ce minimum est appel\'ee \emph{courbe
ferm\'ee g\'eod\'esique}.
\end{defi}
%\begin{demo}[de la proposition~\ref{p.module-courbes}]
%La proposition vient des deux remarques suivantes~:
%\begin{itemize}
%\item toute courbe libre est inclus dans un domaine de translation~;
%\item si $x$ et $y$ appartienne \`a un m\^eme domaine de translation,
%alors il existe un arc libre allant de $x$ \`a $y$, sauf si $x$ et $y$
%sont dans la m\^eme orbite de $h$. 
%\end{itemize}
%Les d\'etails sont les m\^emes que pour la preuve du corollaire 3.5 du texte~\cite{I}.
%\end{demo}

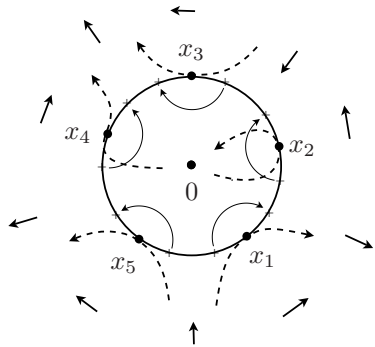
\begin{figure}[htbp]
%\begin{center}
\def\JPicScale{0.8}
\hspace{2cm}
\ifx\JPicScale\undefined\def\JPicScale{1}\fi
\psset{unit=\JPicScale mm}
\psset{linewidth=0.3,dotsep=1,hatchwidth=0.3,hatchsep=1.5,shadowsize=1}
\psset{dotsize=0.7 2.5,dotscale=1 1,fillcolor=black}
\psset{arrowsize=1 2,arrowlength=1,arrowinset=0.25,tbarsize=0.7 5,bracketlength=0.15,rbracketlength=0.15}
\begin{pspicture}(0,0)(100,55.53)
\psbezier(76.47,27.3)(76.47,27.3)(76.47,27.3)(76.47,27.3)
\psbezier(69.78,29.81)(69.78,29.81)(69.78,29.81)(69.78,29.81)
\psdots[](69.78,29.81)
(69.78,29.81)
\psline{->}(70.21,0)(70.11,3.75)
\psline{->}(70.4,55.4)(66.28,55.53)
\psline{->}(95.35,18.14)(100,16)
\psline{->}(85.01,4.61)(89.27,7.7)
\psline{->}(54.61,50.12)(51.92,53.95)
\psline{->}(87.4,48.77)(84.67,45.05)
\psline{->}(96.11,34.16)(95.2,39.72)
\psline{->}(44.75,36.76)(46.56,41.5)
\psline{->}(44.16,21.1)(39.37,19.67)
\psline{->}(54.13,5.44)(50.58,8.1)
\rput(80.24,11.85){}
\rput(99.82,34.64){}
\rput{0}(69.82,29.64){\psellipse[](0,0)(15,15)}
\psdots[](69.82,44.64)
(55.89,35)(61.07,17.5)(79.01,17.95)(84.37,32.94)
\pscustom[linestyle=dashed,dash=1 1]{\psbezier{-}(66.07,7.32)(65.35,13.21)(63.39,15.53)(60.89,17.32)
\psbezier{->}(60.89,17.32)(58.39,19.11)(54.28,19.28)(49.46,17.14)
}
\pscustom[linestyle=dashed,dash=1 1]{\psbezier{-}(73.93,6.43)(75,14.28)(76.6,16.07)(79.1,17.86)
\psbezier{->}(79.1,17.86)(81.6,19.64)(83.21,19.46)(90.18,18.03)
}
\pscustom[linestyle=dashed,dash=1 1]{\psbezier{-}(73.57,28.21)(82.5,25.53)(84.82,29.46)(84.28,32.68)
\psbezier{->}(84.28,32.68)(83.75,35.89)(80.35,37.14)(73.57,32.32)
}
\pscustom[linestyle=dashed,dash=1 1]{\psbezier{-}(80.71,49.46)(76.43,44.64)(73.03,44.82)(69.64,44.82)
\psbezier{->}(69.64,44.82)(66.25,44.82)(63.39,46.43)(60.71,50.36)
}
\pscustom[linestyle=dashed,dash=1 1]{\psbezier{-}(65,29.28)(56.6,29.46)(53.93,30.89)(55.71,35.36)
\psbezier{->}(55.71,35.36)(57.5,39.82)(55.53,41.25)(53.21,44.82)
}
\rput{0}(56.02,35.02){\parametricplot[linewidth=0.1,arrows=->]{-89.53}{46.49}{ t cos 5.65 mul t sin 5.65 mul }}
\rput{90}(78.88,17.86){\parametricplot[linewidth=0.1,arrows=<-]{-30.34}{105.65}{ t cos 5.64 mul t sin 5.64 mul }}
\rput{0}(61.04,17.37){\parametricplot[linewidth=0.1,arrows=->]{-14.62}{121.32}{ t cos 5.65 mul t sin 5.65 mul }}
\rput{0}(84.42,32.87){\parametricplot[linewidth=0.1,arrows=<-]{125.45}{261.3}{ t cos 5.65 mul t sin 5.65 mul }}
\rput{0}(69.9,44.7){\parametricplot[linewidth=0.1,arrows=<-]{-157.21}{-21.26}{ t cos 5.64 mul t sin 5.64 mul }}
\psdots[linewidth=0.1,dotsize=0.7 4.5,dotstyle=+](66.11,15.09)
(57.23,21.56)
\psdots[linewidth=0.1,dotsize=0.7 4.5,dotstyle=+](73.66,15.18)
(82.72,21.87)
\psdots[linewidth=0.1,dotsize=0.7 4.5,dotstyle=+](84.55,27.23)
(82.18,38.12)
\psdots[linewidth=0.1,dotsize=0.7 4.5,dotstyle=+](75.44,43.53)
(64.37,43.57)
\psdots[linewidth=0.1,dotsize=0.7 4.5,dotstyle=+](54.91,29.46)
(59.1,40.18)
\rput(69.82,25.36){$0$}
\rput(69.1,10.71){}
\rput(81.6,14.11){$x_1$}
\rput(88.21,32.68){$x_2$}
\rput(69.82,48.57){$x_3$}
\rput(50.89,34.82){$x_4$}
\rput(58.57,13.21){$x_5$}
\end{pspicture}
\caption{Mot et courbe g\'eod\'esique}
\label{f.mot-courbes}
%\end{center}
\end{figure}

%%%%%%%
\paragraph{Indice et courbes}
L'indice par quarts-de-tour peut \^etre d\'efini en consid\'erant une
courbe ferm\'ee g\'eod\'esique $\gamma$, et  en regardant la dynamique
des sommets d'une d\'ecomposition de $\gamma$.
Commen\c{c}ons par quelques d\'efinitions combinatoires.
On d\'efinit \emph{l'alpha} et \emph{l'om\'ega} d'une fl\`eche verticale~: 
$$
\alpha(\da)=\omega(\ua)=0 \ \ \mbox{ et } \ \  
\alpha(\ua)=\omega(\da)=\infty 
$$
(symboliquement, on voit $0$ ``en haut'' et l'infini ``en bas'', la
fl\`eche $\ua$ repr\'esente une orbite allant de l'infini \`a $0$).
Prenons maintenant un mot (non cyclique) de trois lettres
$M=(m_{-1} m_0 m_1)$, tel que $m_{-1}$ et $m_1$ soient des fl\`eches
verticales, et $m_0$ soit une fl\`eche horizontale
(il y a huit mots possibles).
On d\'efinit  \emph{l'alpha} et \emph{l'om\'ega} du mot $M$ de la fa\c con 
 suivante~:
\begin{itemize}
\item si $m_0=\ra$, alors $\displaystyle
\left\{
\begin{array}{l}
\alpha(M)=\alpha(m_{-1}) \\
\omega(M)=\omega(m_1) \ \ ;
\end{array}
\right.
$
\item  si $m_0=\la$, alors $\displaystyle
\left\{
\begin{array}{l}
\alpha(M)=\alpha(m_{1}) \\
\omega(M)=\omega(m_{-1}).
\end{array}
\right.
$
\end{itemize}
Par exemple, $\alpha(\da \ra \ua)=\omega(\da \ra \ua)=0$
(symboliquement, le mot $\da \ra \ua$ repr\'esente une orbite qui vient
de $0$ et redescend vers $0$ apr\`es avoir fait un demi-tour sur sa gauche).
%D'autre part, on notera $\alpha(x)$ et $\omega(x)$ les ensembles $\alpha$-limite et $\omega$-limite du point $x$ dans la sph\`ere $\bbR^2 \cup \{\infty\}$ (autrement dit, les valeurs d'adh\'erence de la suite des it\'er\'es n\'egatifs (resp. positifs) du point $x$). Rappelons que pour tout point $x$ du plan, $\alpha(x)$ est r\'eduit \`a l'un des points $0$ ou $\infty$ (voir par exemple~\cite{0}, corollaire 3.13 p37).
Le th\'eor\`eme~\ref{t.indice-courbes} caract\'erise le mot $\mathbf{M}(h)$
\`a partir d'une courbe de Jordan g\'eod\'esique.  Plus pr\'ecis\'ement,
les fl\`eches horizontales sont d\'ecrites  en fonction des transitions entre les morceaux d'une d\'ecomposition minimale ;  et les fl\`eches verticales sont d\'ecrites en fonction du type d'orbite des sommets. Ceci est illustr\'e par la figure~\ref{f.mot-courbes}.

\pagebreak[4]
\begin{theo}
\label{t.indice-courbes}
Soit $h$ un hom\'eomorphisme du plan,
pr\'eservant l'orientation, fixant le point $0$,  tel que l'indice de
Poincar\'e-Lefschetz du point fixe est diff\'erent de $1$.
 Soit $\textbf{M}(h)$ l'indice par quarts-de-tour de $h$.
Soit $\gamma$ une courbe de Jordan g\'eod\'esique  munie d'une
d\'ecomposition minimale 
$$
\gamma = \gamma_1 \star \cdots \star \gamma_d
%\gamma=x_0 \frac{\ \gamma_1 \ }{} x_2 \frac{\ \gamma_3 \ }{} \cdots
%\frac{\ \gamma_{2d-1} \ }{} x_{2d}=x_0 .
$$
et $x_1, \dots x_d$  les sommets de cette d\'ecomposition (de fa\c{c}on
 \`a ce que $\gamma_k=[x_k x_{k+1}]_\gamma$).
 Il existe alors une indexation 	 telle que, pour tout entier  $k$~:
\begin{enumerate}
\item la lettre $m_{2k}$ est caract\'eris\'ee par~:
	\begin{itemize}
	\item $
	m_{2k} = \ra \ \ \Longleftrightarrow \ \ h(\gamma_{k-1}) \cap \gamma_k
	\neq \emptyset,
	$
	\item 
	$
	m_{2k} = \la \ \ \Longleftrightarrow \ \ \gamma_{k-1} \cap h(\gamma_k)
	\neq \emptyset \ ;
	$
	\end{itemize}
\item
	\begin{itemize}
	\item $
	\alpha(x_k)= \alpha(m_{2k-1} \ m_{2k} \ m_{2k+1}) , 
	$
	\item $
	\omega(x_k)= \omega(m_{2k-1} \ m_{2k} \ m_{2k+1}).
	$
	\end{itemize}
\end{enumerate}
\end{theo}

%%%%%%%%%%%
\subsection{Lien avec l'indice de Poincar\'e-Lefschetz}
\label{ss.lien-poincare}
On dira qu'un mot $M$ sur l'alphabet $\{\rightarrow, \leftarrow,
\uparrow, \downarrow\}$, cyclique ou non, 
est \emph{autoris\'e} si les lettres de $M$ sont alternativement des
fl\`eches horizontales et des fl\`eches verticales.
On d\'efinit \emph{l'indice partiel symbolique}  des  mots (non cycliques) autoris\'es \`a
deux lettres par  les formules suivantes~: 
$$
\inpa(\da \ra)=\inpa(\ra \ua) = \inpa(\ua \la)= \inpa(\la \da)=
+\frac{1}{4},
$$
$$
\inpa(\ua \ra)= \inpa(\ra \da)=\inpa(\da \la)=\inpa(\la \ua) =
-\frac{1}{4}.
$$
Autrement dit, l'indice partiel symbolique est la variation angulaire
alg\'ebrique, en nombre de tours, quand on passe contin\^ument de la
premi\`ere fl\`eche \`a la seconde (en suivant le chemin le plus court, c'est-\`a-dire en tournant, en valeur absolue, d'un quart de tour et non pas de trois quarts de tour).

On \'etend la d\'efinition \`a un mot non cyclique, autoris\'e,
de $\ell$ lettres par la formule
$$
\inpa(m_1 \dots m_\ell)=\sum_{i=1}^{\ell-1}\inpa(m_im_{i+1}),
$$
et \`a un mot \emph{cyclique}, autoris\'e,
de $\ell$ lettres par la formule
$$
\inpa_{c}((m_{i})_{i \in \bbZ})=\sum_{i=1}^{\ell}\inpa(m_im_{i+1}) = \inpa(m_1 \dots m_\ell m_{1}).
$$
Enfin, \emph{l'indice symbolique}  d'un mot cyclique autoris\'e $\mathbf{M}$ \`a $\ell$ lettres est
$$
\indi(\mathbf{M})=\inpa_{c}(M) +1.
$$
Remarquons que les mots cycliques autoris\'es ont n\'ecessairement un nombre pair de lettres, que cette formule ne d\'epend pas de l'indexation choisie, et que l'indice symbolique d'un mot cyclique autoris\'e est un entier.
Le th\'eor\`eme suivant montre notamment que l'invariant $\mathbf{M}(h)$ est plus fin que l'indice de Poincar\'e-Lefschetz.
%%%%%%%%%%%
\begin{theo}\label{t.lien-indice}
Soit $h$ un hom\'eomorphisme du plan,
pr\'eservant l'orientation, fixant le point $0$,  tel que l'indice de
Poincar\'e-Lefschetz du point fixe est diff\'erent de $1$.
Alors l'indice de Poincar\'e-Lefschetz du point fixe est \'egal \`a l'indice
symbolique du mot cyclique $\mathbf{M}(h)$~:
$$
\indi(h,0)=\indi(\mathbf{M}(h)).
$$
\end{theo}
Ce r\'esultat implique l'in\'egalit\'e suivante :
$$
\diam(h) \geq 2 |\indi(h,0) -1|.
$$ 
Les exemples  dessin\'es sur la figure~\ref{f.indice-PL} ont donc un module minimal parmi les hom\'eomorphismes d'indice de Poincar\'e-Lefschetz donn\'e.
Cette in\'egalit\'e sera un des lemmes importants dans la preuve des r\'esultats pr\'ec\'edents ; il est cependant instructif de voir comment elle d\'ecoule du th\'eor\`eme~\ref{t.lien-indice}.
Notons $\ell$ le nombre de lettres du mot cyclique $\mathbf{M}(h)$~; ce nombre $\ell$ est le double du module de $h$.
L'indice partiel symbolique  du mot cyclique $\mathbf{M}(h)$ est une
somme de $\ell$ termes qui valent chacun $\pm1/4$, par cons\'equent sa valeur absolue est major\'ee par $\frac{\ell}{4}$~:
$$
\frac{\diam(h)}{2} = \frac{\ell}{4} \geq \mid \inpa_{c}(\mathbf{M}(h)) \mid
$$
L'in\'egalit\'e provient maintenant des d\'efinitions de l'indice symbolique et de la formule du th\'eor\`eme~\ref{t.lien-indice}.

On a aussi un r\'esultat de r\'ealisabilit\'e.
\begin{prop}
\label{p.realisabilite}
Soit $\mathbf{M}$ un mot cyclique de longueur paire, constitu\'e de lettres qui sont alternativement des fl\`eches horizontales ($\ra$ ou $\la$) et des fl\`eches verticales ($\ua$ ou $\da$).
Supposons de plus que l'indice symbolique de ce mot soit diff\'erent de $1$.
Alors il existe un hom\'eomorphisme du plan $h$, pr\'eservant l'orientation, fixant uniquement $0$, d'indice diff\'erent de $1$, tel que $\mathbf{M}(h)=\mathbf{M}$.
\end{prop}
Les hom\'eomorphismes de la proposition fournissent donc une infinit\'e d'exemples qui illustrent ce texte. La construction est donn\'ee \`a la section~\ref{ss.realisabilite}, elle est ind\'ependante du reste de l'article.

%
%Ce r\'esultat implique une minoration du module de $h$ en fonction de l'indice de Poincar\'e-Lefschetz.
%\begin{coro}
%\label{c.minoration-diametre}
%Sous les  hypoth\`eses du th\'eor\`eme~\ref{t.lien-indice},  le module de
%$h$ est minor\'e par la formule~:
%$$
%\diam(h) \geq 2 |\indi(h,0) -1|.
%$$ 
%\end{coro}
%Les exemples  dessin\'es sur la figure~\ref{f.indice-PL} ont donc un module minimal parmi les hom\'eomorphismes d'indice de Poincar\'e-Lefschetz donn\'e.
%En r\'ealit\'e, ce ``corollaire'' sera un des lemmes importants dans la preuve des r\'esultats pr\'ec\'edents. Il est cependant instructif de voir pourquoi le th\'eor\`eme~\ref{t.lien-indice} implique cette minoration.
%Notons $\ell$ le nombre de lettres du mot cyclique $\mathbf{M}(h)$~; ce nombre $\ell$ est le double du module de $h$.
%L'indice partiel symbolique  du mot cyclique $\mathbf{M}(h)$ est une
%somme de $\ell$ termes qui valent chacun $\pm1/4$, par cons\'equent sa valeur absolue est major\'ee par $\frac{\ell}{4}$~:
%$$
%\frac{\diam(h)}{2} = \frac{\ell}{4} \geq \mid \inpa_{c}(\mathbf{M}(h)) \mid
%$$
%La minoration provient maintenant des d\'efinitions de l'indice symbolique et de la formule du th\'eor\`eme~\ref{t.lien-indice}.

%

%%%%%%%%%%%%%%%%%%%%%%%%%%%%%%%%%%
\subsection{Cons\'equences dynamiques}
De nombreuses propri\'et\'es dynamiques se lisent sur l'indice par quarts-de-tour. Elles ont \'et\'e rel\'egu\'ees dans la partie~\ref{p.consequences-dynamiques}.

%%%%%%%%%%%%%%%%%%%%%%%%%%%%%%%%%%%%%%%%%
%\subsection{Exemple complet}

%%%%%%%%%%%%%%%%%%%%%%%
\subsection{Structure du texte}
%%%%%%%%%%%%%%%%%%%%%%%
%%%%%%%%%%%%%%%%%%%%%%%
%\subsubsection*{Pr\'erequis}
Le texte s'appuie de mani\`ere fondamentale sur les r\'esultats de~\cite{I}, qui contient la construction d'un invariant du m\^eme type dans le cadre plus simple des hom\'eomorphismes du plan sans point fixe. En particulier, une partie du texte consiste \`a v\'erifier que les r\'esultats obtenus dans~\cite{I} se transposent au cadre des hom\'eomorphismes avec un point fixe. La principale nouveaut\'e concerne la formule qui relie l'indice par quarts-de-tour \`a l'indice de Poincar\'e-Lefschetz. Pour l'obtenir, nous aurons besoin des techniques de construction de droites de Brouwer de P. Le Calvez et A. Sauzet (\cite{sauzet,lecasau}), utilisant les \emph{d\'ecompositions en briques}, et de la notion d'indice partiel d\'evelopp\'ee dans le texte~\cite{0} (sections 3.3 et 3.4 de ~\cite{0}).
Nous utiliserons \'egalement l'existence d'un relev\'e canonique, et le th\'eor\`eme~D du texte~\cite{0} (existence de croissants et de p\'etales). 
%Une preuve plus \'el\'ementaire, \'evitant ces deux derniers recours \`a~\cite{0}, est donn\'ee dans l'appendice~\ref{as.indice-1} ;
%\footnote{??}
% cette preuve alternative permet \'egalement de voir que \emph{tous les r\'esultats de ce texte restent valable en indice $1$, \`a condition de supposer que le module de $h$ est minor\'e par $4$.}

Dans la partie~\ref{p.squelette}, on d\'emontre les r\'esultats annonc\'es \`a partir d'un certain nombre de lemmes. Ces lemmes sont prouv\'es dans la partie~\ref{p.preuves-lemmes}. La partie~\ref{p.consequences-dynamiques} contient les \'enonc\'es et les preuves des liens entre l'indice par quarts-de-tour et la dynamique. 
%Pour all\'eger le texte principal, on a mis en appendice~\ref{as.graphe} la construction d'un invariant plus riche (un graphe simplicial orient\'e, qui est en g\'en\'eral d\'enombrable) ; en 
L'appendice~\ref{as.indice-1} porte sur la g\'en\'eralisation des r\'esultats en indice $1$.

\part{Squelette de preuve}
\label{p.squelette}
 Dans cette partie, on prouve presque tous les r\'esultat annonc\'es \`a la partie pr\'ec\'edente, en admettant provisoirement un certain nombre de  lemmes interm\'ediaires.
%   ;  les exceptions concernent  
%\begin{itemize}
%\item la proposition~\ref{p.module-courbes}, d\'emontr\'ee \`a la section~\ref{s.module-courbes}, qui caract\'erise le module comme la $h$-longueur d'une courbe g\'eod\'esique ;
%\item les r\'esultats dynamiques de la section~\ref{ss.consequences}, dont la preuve fait l'objet de la section~\ref{s.consequences}.
%\end{itemize}
 Dans la partie suivante, on d\'emontrera ces lemmes interm\'ediaires, y compris la proposition~\ref{p.module-courbes}  qui caract\'erise le module comme la $h$-longueur d'une courbe g\'eod\'esique. 
 Pour toute cette partie, on se donne un hom\'eomorphisme $h$ du plan  v\'erifiant les hypoth\`eses de la section pr\'ec\'edente : $h$ pr\'eserve l'orientation, le point $0$ est son unique point fixe, et son indice de Poincar\'e-Lefschetz, not\'e $\indi(h,0)$, est diff\'erent de $1$.
 
% HYPOTHESES NON RAPPELEES
%dans les enonces : les rappeler dans les versions donn\'ees avant les demos. 
%Changer $\wt h $ en $H$.

%%%%%%%%%%%%%%%%%%%%%%%%
%%%%%%%%%%%%%%%%%%%%%%%%
\section{\'Enonc\'es  des lemmes}
\label{s.lemmes}
Dans cette section, on \'enonce les r\'esultats interm\'ediaires, tout en d\'efinissant les objets n\'ecessaires.  Voici les diff\'erentes \'etapes :
\begin{itemize}
\item[\textbf{\ref{ss.releve-canonique}}] d'apr\`es le texte~\cite{0}, l'hom\'eomorphisme $h$ admet  un relev\'e canonique $\wt h$ ;  
\item[\textbf{\ref{ss.geodesiques}}] une courbe de Jordan g\'eod\'esique pour $h$ se rel\`eve alors en une droite g\'eod\'esique pour $\wt h$, qui est essentiellement unique ;  
\item[\textbf{\ref{ss.dynamique-equivariants}}] d'apr\`es le texte~\cite{I}, \`a une telle droite g\'eod\'esique sont associ\'es canoniquement une suite  de composantes de Reeb pour $\wt h$, et un mot infini $\wt{\mathbf{M}}$ qui d\'ecrit grossi\`erement la dynamique de ces composantes ;
\item[\textbf{\ref{ss.compatibilite-composantes}}]  ces composantes de Reeb pour $\wt h$ se projettent sur des composantes de Reeb pour $h$. 
\end{itemize}

%%%%%%%%%%%%%%%%%%%%%%%
\subsection{Relev\'e canonique}
\label{ss.releve-canonique}
Soit $\pi$  le rev\^etement universel\footnote{Le signe `$-$' devant le `$r$' est une commodit\'e,
 introduite pour des raisons essentiellement visuelles : le point $0$ se
retrouve ``\`a l'infini en haut'' du rev\^etement. Par ailleurs, la lettre  ``$\pi$'' \`a l'int\'erieur de la formule n'a bien s\^ur rien \`a voir avec la lettre  ``$\pi$'' d\'esignant l'application de rev\^etement.}

$$
\begin{array}{rcll}
\bbR^2 & \longrightarrow & \bbC \setminus \{0\} &\simeq \bbR^2 \setminus \{0\} \\
(\theta,r) & \longmapsto & \exp(-r + 2i\pi \theta).&
\end{array}
$$
%\footnote{Un \emph{relev\'e} d'un point $x \neq 0$ est un point $\wt x$ tel que
%$\pi(\wt x)=x$. Un \emph{relev\'e} d'une courbe $\gamma : [0,1] \ra \bbR^2
%\setminus\{0\}$ est une courbe $\wt \gamma : [0,1] \ra \bbR^2$ telle que $\pi
%\circ \wt \gamma=\gamma$.}
 Un \emph{relev\'e} de l'hom\'eomorphisme $h$
 est un hom\'eomorphisme  $\wt h$ du plan $\bbR^2$
 tel que  $h \circ \pi=\pi \circ \wt h$.
On note $\tau$ l'automorphisme de rev\^etement $(\theta,r) \mapsto
(\theta+1,r)$.
% \'Etant donn\'e un relev\'e $\wt h $ de $h$, les autres relev\'e
%de $h$ sont les hom\'eomorphismes $\tau^k \circ \wt h$ o\`u $k$
%parcourt l'ensemble des entiers. 
Comme $h$ pr\'eserve l'orientation, il est isotope \`a l'identit\'e parmi les hom\'eomorphismes fixant $0$ (th\'eor\`eme de Kneser), et par cons\'equent tout relev\'e de $h$ commute avec la translation $\tau$.
D'autre part, $\wt h$ n'a aucun point fixe.
Ces remarques permettent de d\'efinir \emph{l'indice de $\wt h$ relativement \`a $\tau$} : c'est  l'indice de $\wt h$ le long de n'importe quelle courbe reliant un point \`a son image par $\tau$ (voir~\cite{0}, p78). C'est un entier, on le note $\indi_{\tau}(\wt h)$.

 Commen\c cons par rappeler l'existence d'un relev\'e canonique.
  L'\'enonc\'e suivant est \'equivalent \`a la proposition 4.17 du texte~\cite{0} (p 80). 
 \begin{prop}[relev\'e canonique]
 \label{p.releve-canonique}
Il existe un unique relev\'e $\wt h$ de $h$ par $\pi$ tel que les indices v\'erifient la relation
$$
\indi_{\tau}(\wt h) = \indi(h,0) -1.
$$
On a de plus :
\begin{enumerate}
%\item\label{i.brouwer} correspondance des droites de Brouwer ;
\item\label{i.libres} si $\wt \gamma$ est un arc libre\footnote{Rappelons que ``libre'' signifie ``disjoint de son image''.}
 pour $\wt h$, alors  $\gamma:=\pi(\wt \gamma)$  est un arc libre pour $h$  ;
\item\label{i.decompositions} si $\cF$ est une d\'ecomposition en briques pour $h$, alors $\wt \cF = \pi^{-1}(\cF) $ est une d\'ecomposition en briques pour $\wt h$.
\end{enumerate}
\end{prop}
 
 Nous rappellerons en temps utile ce qu'est une d\'ecomposition en briques.
Les r\'esultats du texte~\cite{0} nous permettrons \'egalement de minorer le module de $h$.
\begin{prop}[minoration du module]
\label{p.minoration-module}
On a
$$
\diam(h) \geq 2\mid  \indi(h,0)-1       \mid.
$$
En particulier, si l'indice du point fixe est diff\'erent de $0$, $1$ ou $2$, alors le module est sup\'erieur ou \'egal \`a $4$.
\end{prop}
 
%Puisque $0$ est l'unique point fixe de $h$,  le relev\'e canonique $\wt h$ est un hom\'eomorphisme du plan qui pr\'eserve l'orientation et qui n'a aucun point fixe.
% C'est ce qu'on appelle un \emph{hom\'eomorphisme de Brouwer}.
%De plus $\wt h $ a la particularit\'e de commuter avec la translation $\tau$. 

%%%%%%%%%%%%%%%%%%%%%%%
\subsection{G\'eod\'esiques locales et globales}
\label{ss.geodesiques}
On consid\`ere un hom\'eomorphisme de Brouwer $H$, c'est-\`a-dire un hom\'eomorphisme du plan, sans point fixe, pr\'eservant l'orientation  (sans hypoth\`ese de commutation pour le moment). Dans ce cadre, on a d\'efini dans le texte~\cite{I} une distance $d_{H}$, \`a valeurs enti\`eres (voir~\cite{I}, section~3.1 ; cette distance compte le nombre minimal de domaines de translation, ou (de mani\`ere \'equivalente) le nombre minimal d'arcs libres, n\'ecessaires pour connecter deux points donn\'es).
%, et une notion de droite g\'eod\'esique (\cite{I}, section~7.1).
\begin{defi}
\label{d.geodesique-locale}
Soit $I$ un intervalle de $\bbZ$, et $d$ un entier positif. Une suite
$(x_i)_{i  \in I}$ de  points du plan est une \emph{suite
g\'eod\'esique locale d'ordre $d$} si pour tous indices $i,j$ dans
$I$,
$$
|j-i| \leq d \Rightarrow d_H(x_i,x_j)=|j-i|.
$$ 
\end{defi}
%	\begin{rema}\label{r.geo-locale}
%	Consid\'erons $(x_i)_{i  \in I}$  une suite g\'eod\'esique
%	locale d'ordre $1$ (autrement dit, deux points cons\'ecutifs sont \`a
%	$H$-distance $0$ ou $1$). Soit $d \geq 3$. 
%	Pour que cette suite soit une suite g\'eod\'esique
%	locale d'ordre $d$, il suffit que l'implication de la d\'efinition soit
%	v\'erifi\'ee quand $|i-j|=d$.
%	 Cette remarque provient simplement de l'in\'egalit\'e triangulaire.
%	\end{rema}

Rappelons qu'une suite  de points du plan est une \emph{suite g\'eod\'esique}
(``globale'') si c'est une suite g\'eod\'esique locale d'ordre $d$ pour
tout $d$. L'existence de composantes de Reeb, obtenue dans le texte~\cite{I}, nous permettra de montrer que les g\'eod\'esiques locales sont des g\'eod\'esiques globales. 
\begin{prop}[g\'eod\'esiques locales]
\label{p.geodesique-locale}
Soit $H$ un hom\'eomorphisme de Brouwer, et $(x_i)_{i \in I}$ une
suite g\'eod\'esique locale d'ordre $d$, avec $d \geq 3$. Alors $(x_i)_{i \in I}$ est
une suite  g\'eod\'esique.
\end{prop}

Une suite  $(x_{k})_{k \in \bbZ}$  est \emph{$\tau$-\'equivariante}, de p\'eriode $d \geq 1$, si on a  $x_{k+d}=\tau(x_{k})$ pour tout entier $k$.
Rappelons que deux suites g\'eod\'esiques  
%pour l'hom\'eomorphisme $H$ 
$(x_k)_{k \in \bbZ}$ et  $(x_k')_{k \in \bbZ}$ \emph{d\'efinissent les m\^emes
bouts} si il existe un entier $p$ tel que pour tout entier $k$,
$d_H(x_k,x'_k) \leq p$ (voir \cite{I}, section 7.1). Cette relation est une relation d'\'equivalence.
La proposition pr\'ec\'edente entra\^inera l'existence d'une unique classe d'\'equivalence de g\'eod\'esiques $\tau$-\'equivariantes (rappelons que la notion de \emph{droite g\'eod\'esique} est d\'efinie dans~\cite{I}, et que la proposition~\ref{p.module-courbes} nous fournit une courbe de Jordan g\'eod\'esique).
\begin{coro}\label{c.bouts-geodesique}
Soit $\wt h$ le relev\'e canonique fourni par la proposition~\ref{p.releve-canonique}, et $d$ le module de $h$. Soit $\gamma$ une courbe de Jordan g\'eod\'esique pour $h$. Soit $\Gamma$ la droite topologique de $\bbR^2$ qui rel\`eve $\gamma$.
\begin{enumerate}
\item La droite $\Gamma$ est une droite g\'eod\'esique pour $\wt h$. Plus pr\'ecis\'ement, toute d\'ecomposition minimale de $\gamma$ se rel\`eve en une d\'ecomposition minimale de $\Gamma$, et la suite des sommets de cette d\'ecomposition est une suite g\'eod\'esique pour $\wt h$ qui est $\tau$-\'equivariante, de p\'eriode $d$.
\item Si deux suites g\'eod\'esiques  pour $\wt h$  sont $\tau$-\'equivariantes de p\'eriode $d$, alors elles  d\'efinissent les m\^emes bouts de $\wt h$.
\end{enumerate}
\end{coro}

\begin{defi}
On notera $[\tau]$ l'unique classe d'\'equivalence de suites g\'eod\'esiques d\'efinie par le corollaire.
\end{defi}
Voici l'id\'ee de la preuve du premier point du corollaire. Si l'indice du point fixe de $h$ est diff\'erent de $0$, $1$ ou $2$, alors le module $d$ est minor\'e par $4$ (proposition~\ref{p.minoration-module}). D'apr\`es la proposition~\ref{p.module-courbes}, le module est \'egal \`a la $h$-longueur de la courbe de Jordan g\'eod\'esique $\gamma$.
En relevant  les $d$ sommets d'une d\'ecomposition minimale de $\gamma$, on obtient une suite $\tau$-\'equivariante. Cette suite est une g\'eod\'esique locale d'ordre $d$. Puisque $d \geq 4$, la proposition pr\'ec\'edente dit que cette suite est en fait une suite g\'eod\'esique ``globale''.

Le cas o\`u l'indice du point fixe vaut $0$ ou $2$ est plus d\'elicat. N\'eanmoins, on peut se ramener au cas pr\'ec\'edent en montrant que $h$ admet un (unique) relev\'e $h_{2}$, par le rev\^etement  \`a deux feuillets ramifi\'e au-dessus de $0$, qui est d'indice $-1$ ou $3$.

%%%%%%%%%%%%%%%%%%%%%%%
\subsection{Dynamique des hom\'eomorphismes de Brouwer \'equivariants}
\label{ss.dynamique-equivariants}
On va maintenant utiliser les composantes de Reeb associ\'ees \`a une suite g\'eod\'esique.
 Plus pr\'ecis\'ement, nous allons utiliser les r\'esultats des sections 7 et 8 de~\cite{I}~: \`a chaque suite g\'eod\'esique infinie sont associ\'es (canoniquement) une suite de composantes de Reeb et un mot infini, et ce mot infini contient des informations sur la dynamique dans les bords des composantes de Reeb et sur les orbites des points de la suite g\'eod\'esique. De plus, ces objets ne d\'ependent que de la classe d'\'equivalence de la suite g\'eod\'esique modulo la relation ``d\'efinir les m\^emes bouts''.
  Lorsque l'hom\'eomorphisme est \'equivariant, on obtient ainsi un mot p\'eriodique~; dans le th\'eor\`eme~\ref{t.lien-indice}-bis, on reliera l'indice partiel du mot cyclique induit \`a l'indice de $\wt h$ par rapport \`a la translation $\tau$. 

D\'ecrivons maintenant tout ceci en d\'etail. On note toujours $\wt h$ le relev\'e canonique, $d$ le module. On se donne aussi une suite $(\wt x_{k})_{k \in \bbZ}$ g\'eod\'esique pour $\wt h$ qui est $\tau$-\'equivariante, de p\'eriode $d$. Cette suite est la suite des sommets d'une droite g\'eod\'esique $\Gamma$, qui est invariante par $\tau$ (corollaire~\ref{c.bouts-geodesique}).

%\paragraph{Suite et droite g\'eod\'esiques}
% En relevant une courbe de Jordan $\gamma$ g\'eod\'esique pour $h$, ainsi que ses sommets\footnote{On n'a pas \'enonc\'e ceci ???}, on obtient~:
%\begin{itemize}
%\item une droite topologique $\Gamma :  \bbR \ra \bbR^2$ v\'erifiant $\tau \circ \Gamma (t) =\Gamma(t +1)$ ;
%\item une suite $(\wt x_{k})_{k \in \bbZ}$ de points de $\Gamma$ qui est une suite g\'eod\'esique $\tau$-\'equivariante, telle que chaque arc $\Gamma_{k}= [\wt x_{k} \wt x_{k+1}]_{\Gamma}$ est libre.
%\end{itemize}
%Dans l'article~\cite{I}, une telle droite topologique $\Gamma$ \'etait appel\'ee \emph{droite g\'eod\'esique}.

\paragraph{Composantes de Reeb}
 La d\'efinition des composantes de Reeb dans le cadre des hom\'eomorphismes de Brouwer se trouve \`a la section~4.2 de~\cite{I}.
 On note  $(\wt F_{k},\wt G_{k})_{k \in \bbZ}$ la suite de composantes de Reeb associ\'ee
\`a la suite $(\wt x_{k})_{k \in \bbZ}$ (\cite{I}, section 7).
On a notamment 
$$
\wt F_{k} \cap \Gamma \subset [\wt x_{k-1} \wt x_{k}]_{\Gamma} \mbox{ et }  
\wt G_{k} \cap \Gamma \subset [\wt x_{k} \wt x_{k+1}]_{\Gamma}.
$$
Comme la suite $(\wt x_{k})_{k \in \bbZ}$ est \'equivariante, par unicit\'e des composantes de Reeb (th\'eor\`eme~2 de ~\cite{I}), la suite des composantes est elle-m\^eme \'equivariante~: on a $(\wt F_{k+d},\wt G_{k+d}) = (\tau(\wt F_{k}), \tau(\wt G_{k}))$ pour tout $k$. Rappelons aussi que chacun des bords $\wt F_{k}$ et $\wt G_{k}$ est invariant par $\wt h$ (\cite{I}, lemme 7.10).

%\paragraph{D\'ecomposition}
%On consid\`ere maintenant une d\'ecomposition en briques $\mathcal{F}$ pour $h$, g\'eod\'esique, compatible avec la d\'ecomposition de la courbe de Jordan g\'eod\'esique $\gamma$ (lemme~\ref{l.decomposition-geodesique}). L'ensemble $\mathcal{\wt F}= \pi^{-1}(\mathcal{F})$ est une d\'ecomposition en briques pour $H$  (lemme~\ref{l.releve-decomposition})~; elle est invariante par $\tau$~;
%et  $\mathcal{\wt F} \cap \Gamma$  est l'ensemble des points de la suite $(\wt x_{k})_{k \in \bbZ}$.

\paragraph{Mot cyclique}
On note $\wt {\textbf{M}}([\tau])$ le mot infini, sur l'alphabet $\{ \ra, \la, \da, \ua\}$, associ\'e
\`a $[\tau]$ (\cite{I}, th\'eor\`eme~3 et  d\'efinition~8.2).
Comme $\wt h$ commute avec la translation $\tau$, ce mot infini est p\'eriodique de p\'eriode $2d$. On note $\textbf{M}([\tau])$ le mot cyclique \`a $2d$ lettres induit par ce mot p\'eriodique.
\paragraph{Formule d'indice}%\label{ss.indice-brouwer}
 %%%%%%%%%%%%
On pourra alors relier l'indice du relev\'e canonique \`a l'indice symbolique.
\begin{theo*}[\ref{t.lien-indice}-bis]
L'indice de $\wt h$ par rapport \`a $\tau$ est \'egal \`a l'indice partiel
symbolique du mot cyclique $\textbf{M}([\tau])$~:
$$
\indi_{\tau}(\wt h) = \inpa_{c}(\textbf{M}([\tau])) = \inpa(m_{1} \dots m_{2d} m_{1}).
$$
\end{theo*}
%
%Voici l'id\'ee de la preuve du th\'eor\`eme. (...)
%-- dte de Brouwer (def) ; Rq = pouss\'ee sens des fl\`eches.
%\begin{prop}
%\label{p.cRdB}
%Soit $(F,G)$ une composante de Reeb pour un hom\'eomorphisme de Brouwer $H$. On suppose que $F$ et $G$ sont invariants par $H$. Alors il existe une droite de Brouwer qui s\'epare $F$ et $G$.
%\end{prop}

%

%-- dans 0, on a def un IP entre 2 dtes de Brouwer. Par def = indice tau.
%Et additivit\'e.
%-- donc il reste \`a voir que IP entre 2 dtes successives = IP symbolique.
%a) cs indiff\'erent, 0 \'evident.
%b) cas attractif/rep, on a pu choisir delta transverse \`a Gamma, Delta inter gamma = sommet
%(cf figure) ; du coup, arc libre, pouss\'e sens des fleches. Situation de 0 : OK.

%

 %%%%%%%%%%%%%%%%%%%%%%%
\subsection{Compatibilit\'e des composantes de Reeb}
\label{ss.compatibilite-composantes}
Dans le texte~\cite{I}, on construisait une compactification de chaque bord de composante de Reeb par l'ajout de deux points \`a l'infini, not\'es $N$ et $S$ (on pense \`a $N$ comme \'etant l'infini ``en haut'', et $S$ est l'infini ``en bas''). 
Le th\'eor\`eme~3 de~\cite{I} (section 8.1) dit que la lettre $m_{2k+1}$ vaut $\da$ si et seulement si le bord $\wt F_{k}$ est de type dynamique Nord-Sud, c'est-\`a-dire si il existe un ouvert dense de points de $\wt F_{k}$ dont l'orbite va de $N$ \`a $S$.

On peut \'etendre contin\^ument  l'application de rev\^etement $\pi$ en posant $\hat \pi(N)=0$ et $\hat \pi(S)= \infty$, o\`u $\infty$ d\'esigne le point \`a l'infini du compactifi\'e d'Alexandroff de $\bbR^2$.
La proposition suivante \'etablit une correspondance entre les composantes de Reeb de $h$ et celles de son relev\'e canonique $\wt h$. On note $O(\wt F, \wt G)$ l'unique composante connexe de $\bbR^2 \setminus (\wt F \cup \wt G)$  dont l'adh\'erence rencontre \`a la fois $\wt F$ et $\wt G$ (voir l'appendice~\ref{as.ordre-cyclique} ; rappelons que $O(F,G)$ est d\'efini de fa\c con analogue \`a la section~\ref{ss.composantes-de-reeb}).

 \begin{prop}~
\label{p.compatibilite-composantes}

\begin{enumerate}
\item Si un couple $(\wt F, \wt G)$ est une composante de Reeb minimale pour le relev\'e canonique $\wt h$,
associ\'ee \`a $[\tau]$, alors l'application $\pi$ induit un hom\'eomorphisme de $\wt F \cup O(\wt F, \wt G) \cup \wt G \cup \{N,S\}$ sur son image $F \cup O(F,G) \cup G \cup \{0,\infty\}$, et  
le couple $(F,G)=(\pi(\wt F),\pi(\wt G))$ est une
composante de Reeb minimale pour $h$. 
\item R\'eciproquement,  si $(F,G)$ est une composante de Reeb minimale pour $h$,
alors il existe un couple $(\wt F, \wt G)$ qui  est une composante de
Reeb minimale pour $\wt h$,
associ\'ee \`a $[\tau]$, et tel que $\pi(\wt F)=F$ et $\pi(\wt G)=G$.
\end{enumerate}
\end{prop}
Dans la situation d\'ecrite par la proposition, 
 on dira que $(F,G)$ est la \emph{projet\'ee} de la composante $(\wt F,\t
 G)$, et que  $(\wt F,\wt  G)$ est une composante \emph{relev\'ee} de
 $(F,G)$.

%%%%%%%%%%%%%%%%%%%%%%%
%%%%%%%%%%%%%%%%%%%%%%%
\section{Preuve des th\'eor\`emes}
\label{s.preuves}
\`A partir des lemmes de la section pr\'ec\'edente,
nous d\'emontrons tous les r\'esultats \'enonc\'es dans la section~\ref{s.definitions}
(mise \`a part la proposition~\ref{p.module-courbes}, dont la preuve fait l'objet de la section suivante). Nous reprenons les notations de la section pr\'ec\'edente.

\subsection{Existence et unicit\'e des composantes de Reeb pour $h$ (th\'eor\`eme~\ref{t.compo-Reeb})}
  Notons $F_k:=\pi(\wt F_k)$, $G_k:=\pi(\wt G_k)$ les projet\'es des bords des composantes de Reeb pour $\wt h$ associ\'ees \`a $[\tau]$.
D'apr\`es la proposition~\ref{p.compatibilite-composantes}, les
composantes de Reeb minimales pour $h$ sont exactement les couples
$(F_1,G_1), \dots , (F_d,G_d)$. La propri\'et\'e d'ordre cyclique vient de
l'ordre des composantes $(\wt F_i,\wt G_i)$ (\cite{I}, lemme 7.10).
Le  th\'eor\`eme~\ref{t.compo-Reeb} en d\'ecoule.
%\end{demo}

\subsection{Dynamique dans les bords des composantes (proposition~\ref{p.dynamique-bord} et \ref{p.bords-adjacents})}
%\begin{demo}[de la proposition~\ref{p.dynamique-bord}]
Soit $(F,G)$ une composante de Reeb minimale pour $h$, et $(\wt F, \t
G)$ une composante relev\'ee.
D'apr\`es la proposition~\ref{p.compatibilite-composantes},
% $\wt F$ v\'erifie les hypoth\`eses du lemme
%topologique~\ref{l.compatibilite-topologique}. D'apr\`es le lemme,
l'application $\hat \pi : \wt F \cup \{N,S\} \ra \hat F = F \cup \{ 0,
\infty\}$ est un hom\'eomorphisme qui envoie $N$ sur $0$ et  $S$ sur $\infty$. 
La proposition~\ref{p.dynamique-bord}, concernant la dynamique dans le bord $F$, est alors une simple
traduction du r\'esultat sur la dynamique dans $\wt F$ (proposition~7.17
de \cite{I}).
%\end{demo}

%\begin{demo}[de la proposition~\ref{p.bords-adjacents}]
%L\`a encore, la proposition est une traduction du fait que deux bords
%adjacents $\wt G_k$ et $\wt F_{k+1}$ sont du m\^eme type dynamique pour
%$\wt h$~: ce r\'esultat est contenu dans le th\'eor\`eme~3 de~\cite{I}.
%\end{demo}
D'apr\`es le th\'eor\`eme~3 de~\cite{I}, deux bords adjacents $\wt G_k$ et $\wt F_{k+1}$ sont du m\^eme type dynamique pour $\wt h$. On en d\'eduit imm\'ediatement la proposition~\ref{p.bords-adjacents}~: deux bords adjacents $G_k$ et $F_{k+1}$ sont du m\^eme type dynamique pour $h$.
%De plus, il est clair que le mot $\textbf{M}([\tau])$ (d\'efini en~\ref{ss.dynamique-equivariants}, d'apr\`es la dynamique dans les composantes de Reeb de $\wt h$) et le mot $\textbf{M}(h)$ (d\'efini en~\ref{ss.definition-indice}, d'apr\`es la dynamique dans les composantes de Reeb de $h$) sont identiques.

%Comme not\'e plus haut, l
L'indice par quarts-de-tour, d\'efini \`a la section~\ref{ss.definition-indice} \`a partir de la dynamique dans les composantes de Reeb minimales de $h$,  est alors clairement  un invariant de conjugaison, ce qui est l'\'enonc\'e du th\'eor\`eme 2.

\subsection{Caract\'erisation de l'indice en termes de courbes (th\'eor\`eme~\ref{t.indice-courbes})}
L\`a encore, il s'agit principalement de transposer \`a $h$ les liens entre composantes de Reeb et droite g\'eod\'esique \'etablis pour $\wt h$ dans le texte~\cite{I}, et notamment de montrer que le mot cyclique $\mathbf{M}([\tau])$, d\'efini en~\ref{ss.dynamique-equivariants} \`a partir du texte~\cite{I}, co\"incide  avec le mot cyclique $\textbf{M}(h)$ d\'efini en~\ref{ss.definition-indice}.
On consid\`ere  une courbe de Jordan  $\gamma$,  g\'eod\'esique pour $h$,  munie d'une d\'ecomposition minimale 
$$
\gamma = \gamma_1 \star \cdots \star \gamma_d
%\gamma=x_0 \frac{\ \gamma_1 \ }{} x_2 \frac{\ \gamma_3 \ }{} \cdots
%\frac{\ \gamma_{2d-1} \ }{} x_{2d}=x_0 .
$$
et $x_1, \dots x_d$  les sommets de cette d\'ecomposition (de fa\c{c}on
 \`a ce que $\gamma_k=[x_k x_{k+1}]_\gamma$).
Soit $\Gamma$ la droite topologique qui rel\`eve $\gamma$ ; la d\'ecomposition de $\gamma$ se rel\`eve en une d\'ecomposition
$$
\Gamma = \cdots \star \wt \gamma_1 \star \cdots \star \wt  \gamma_d \star \cdots
$$
qui est \'equivariante pour $\tau$ : $\tau(\wt \gamma_{k}) = \wt \gamma_{k + d}$. De m\^eme, on peut relever la suite cyclique des sommets de $\gamma$ en une suite $(\wt x_{k})_{k \in \bbZ}$ de sommets  de $\Gamma$, de mani\`ere compatible.
On note encore  $(\wt F_{k}, \wt G_{k})$ la composante de Reeb minimale associ\'ee au sommet $\wt x_{k}$, et $(F_{k}, G_{k})$ sa projection par $\pi$.
\begin{affi}
Les six propri\'et\'es suivantes sont \'equivalentes :
\begin{enumerate}
\item $\wt h(\wt \gamma_{k-1}) \cap \wt \gamma_{k} \neq \emptyset$ ;
\item dans le mot infini $\wt {\textbf{M}}([\tau])$, la fl\`eche correspondant au sommet $\wt x_{k}$ vaut $\ra$ ; 
\item la dynamique de $\wt h$ va de $\wt F_{k}$ vers $\wt G_{k}$ ; 
\item la dynamique de $h$ va de $F_{k}$ vers $G_{k}$ ; 
\item\label{5} dans le mot cyclique $\mathbf{M}(h)$ , la fl\`eche correspondant au sommet $x_{k}$ vaut $\ra$ ; 
\item\label{6} $ h( \gamma_{k-1}) \cap \gamma_{k} \neq \emptyset$.
\end{enumerate}
\end{affi}
\begin{demo}[de l'affirmation]~

\begin{itemize}
\item L'\'equivalence entre les propri\'et\'es 1 et 2 d\'ecoule de la d\'efinition des fl\`eches horizontales dans le texte~\cite{I} (d\'efinition~3.11).
\item Le th\'eor\`eme~1-bis de~\cite{I} donne l'\'equivalence entre les propri\'et\'es 2 et 3.
\item  L'\'equivalence entre les propri\'et\'es 3 et 4 vient de la compatibilit\'e des composantes de Reeb de $\wt h$ avec celles de $h$ (proposition~\ref{p.compatibilite-composantes}).
\item  La d\'efinition~\ref{d.mot-cyclique} des fl\`eches horizontales  fournit l'\'equivalence entre les propri\'et\'es 4 et 5. 
\end{itemize}
 Il reste \`a montrer l'\'equivalence entre les propri\'et\'es 1 \`a 5 d'une part, et la propri\'et\'e 6 d'autre part.
Il est clair que 1 implique 6. Supposons maintenant que les propri\'et\'es 1 \`a 5 ne sont pas v\'erifi\'ees ; notamment, dans le mot infini $\wt {\textbf{M}}([\tau])$, la fl\`eche correspondant au sommet $\wt x_{k}$ vaut $\la$. Par d\'efinition des fl\`eches horizontales de $\wt h$,  ceci entra\^ine  $\wt \gamma_{k-1} \cap \wt h(\wt \gamma_{k}) \neq \emptyset$, et par cons\'equent $\gamma_{k-1} \cap h(\gamma_{k}) \neq \emptyset$. 
Si la propri\'et\'e 6 \'etait v\'erifi\'ee, en \'epaississant un peu les arcs $\gamma_{k-1}$ et $\gamma_{k}$, on obtiendrait une cha\^ine p\'eriodique de deux disques, et l'indice du point fixe de $h$ serait \'egal \`a $1$ d'apr\`es le lemme de Franks
(voir par exemple~\cite{0}, lemme~3.18 p39).\footnote{On peut montrer, alternativement, que le module de $h$ est alors \'egal \`a $2$. En effet, puisque d'apr\`es la propri\'et\'e~\ref{6} $h(\gamma_{k-1})$ rencontre $\gamma_{k} \subset \gamma$, $\wt h (\wt \gamma_{k-1})$ doit rencontrer $\Gamma$. Comme $\Gamma$ est g\'eod\'esique, la rencontre ne peut se faire que sur $\wt \gamma_{k-2}$ ou $\wt \gamma_{k}$.
Comme (par hypoth\`ese) la propri\'et\'e 1 n'est pas v\'erifi\'ee,  $\wt h (\wt \gamma_{k-1})$ rencontre en fait 
 $\wt \gamma_{k-2}$, et aucun autre arc de la d\'ecomposition de $\Gamma$. Puisque 
 $h(\gamma_{k-1})$ rencontre $\gamma_{k}$, on en d\'eduit que l'arc $\wt \gamma_{k-2}$ est un relev\'e de l'arc $\gamma_{k}$, qui \'etait d\'ej\`a relev\'e par $\wt \gamma_{k}$ ; par suite le module vaut $2$. Cette remarque, 
  qui \'evite le recours au lemme de Franks pour $h$,
  sera exploit\'ee dans l'appendice~\ref{as.indice-1} pour
   les points fixes d'indice $1$.}  
\end{demo}

Nous pouvons maintenant terminer la preuve du th\'eor\`eme~\ref{t.indice-courbes}. L'\'equivalence entre les points~\ref{5} et~\ref{6} de l'affirmation nous donne le premier point du th\'eor\`eme. 
Pour obtenir le deuxi\`eme point, on applique le th\'eor\`eme~3 de~\cite{I} (section 8.1), qui \'etablit une correspondance entre le type dynamique des bords de la composante de Reeb $(\wt F_{k}, \wt G_{k})$ d'une part, et la dynamique du sommet $\wt x_{k}$. L'application de rev\^etement $\pi$ se restreint en un hom\'eomorphisme entre  $\wt F \cup O(\wt F, \wt G) \cup \wt G \cup \{N,S\}$ et son image.
Comme le point $\wt x _{k}$ est dans l'ouvert $O(\wt F, \wt  G)$, on a pour l'hom\'eomorphisme $h$ une correspondance parall\`ele entre  le type dynamique des bords de la composante de Reeb $(F_{k}, G_{k})$ et la dynamique du sommet $x_{k}$.
Ceci prouve le deuxi\`eme point.

\subsubsection*{Lien avec l'indice de Poincar\'e-Lefschetz (th\'eor\`eme~\ref{t.lien-indice})}
Pour montrer le th\'eor\`eme~\ref{t.lien-indice}, qui relie l'indice par quarts-de-tour \`a l'indice de Poincar\'e-Lefschetz, on applique successivement :
\begin{itemize}
\item la formule qui relie l'indice de Poincar\'e-Lefschetz de $h$  \`a l'indice du relev\'e canonique $\wt h$ par rapport \`a $\tau$  (proposition~\ref{p.releve-canonique}) :
$$
 \indi(h,0)  = \indi_{\tau}(\wt h) +1 ;
$$
\item le th\'eor\`eme~\ref{t.lien-indice}-bis qui relie l'indice de $\wt h$ \`a l'indice partiel symbolique du mot $\textbf{M}([\tau])$ (qui est \'egal \`a  $\textbf{M}(h)$ d'apr\`es le point pr\'ec\'edent) :
$$
\indi_{\tau}(\wt h) = \inpa_{c}(\textbf{M}(h)) ;
$$
\item le lien qui existe, par d\'efinition, entre l'indice partiel symbolique et l'indice symbolique :
$$
\indi(\mathbf{M})=\inpa_{c}(\textbf{M}) +1.
$$
\end{itemize}
On en d\'eduit  la formule attendue,
$$
\indi(h,0)=\indi(\mathbf{M}(h)).
$$

%%%%%%%%%%%%%%%%%%%%%%%
\subsection{R\'ealisabilit\'e}
\label{ss.realisabilite}
On explique ici comment construire une famille d\'enombrable d'exemples illustrant la th\'eorie, qui r\'ealisent tous les indices par quarts-de-tour possibles.

\begin{demo}[de la proposition~\ref{p.realisabilite}]
On commence par se donner des dynamiques mod\`eles sur des secteurs du plan. Plus pr\'ecis\'ement, on consid\`ere l'ensemble des mots autoris\'es de longueur trois, dont les premi\`ere et derni\`ere lettres sont des fl\`eches verticales (la fl\`eche centrale \'etant donc horizontale). On va associer \`a chacun de ces huit mots un hom\'eomorphisme d\'efini sur un secteur  du plan issu de $0$, et fixant uniquement $0$.
La dynamique associ\'ee au mot $(\ua \ra \da)$ est dessin\'ee sur la figure~\ref{f.modeles-secteurs} (\`a gauche). On associe au mot $(\da \la \ua)$ la dynamique inverse de la pr\'ec\'edente ; ces deux mod\`eles sont appel\'es 
 \emph{secteurs hyperboliques}. De m\^eme, on associe aux deux mots 
$(\ua \la \da)$ et $(\da \ra \ua)$ les deux \emph{secteurs elliptiques} correspondant au dessin de droite de la figure~\ref{f.modeles-secteurs} et \`a son inverse.
\footnote{Si on veut \^etre plus pr\'ecis, le mod\`ele hyperbolique est donn\'e par la restriction $S_{h}$ de l'application $(x,y) \mapsto (2x, y/2)$ au quart-de-plan ferm\'e  $\{ x \geq0, y\leq0 \}$ ; 
le mod\`ele de secteur elliptique s'obtient en \'echangeant les r\^oles de $0$ et de $\infty$ : on prend la  restriction $S_{e}$  de la m\^eme application $(x,y) \mapsto (2x, y/2)$ au quart-de-plan ferm\'e  $\{ x \geq0, y\leq0 \} \setminus \{0\} \cup \{\infty \}$.}
\begin{figure}[htbp]
\begin{center}
\def\JPicScale{0.8}
\ifx\JPicScale\undefined\def\JPicScale{1}\fi
\psset{unit=\JPicScale mm}
\psset{linewidth=0.3,dotsep=1,hatchwidth=0.3,hatchsep=1.5,shadowsize=1}
\psset{dotsize=0.7 2.5,dotscale=1 1,fillcolor=black}
\psset{arrowsize=1 2,arrowlength=1,arrowinset=0.25,tbarsize=0.7 5,bracketlength=0.15,rbracketlength=0.15}
\begin{pspicture}(0,0)(90,30)
\psline{->}(12.5,4.64)(17.5,4.64)
\psline{->}(5,4.64)(5,9.64)
\psline{->}(25,4.64)(25,-0.36)
\psline[linewidth=0.1,linestyle=dashed,dash=1 1](0,4.64)(30,4.64)
\pscustom[]{\psbezier{-}(7.63,16.88)(11.84,22.53)(13.11,26.03)(15.24,25.99)
\psbezier{->}(15.24,25.99)(17.37,25.96)(18.95,23.06)(22.5,17.14)
}
\psbezier{->}(10.53,16.35)(13.73,21.69)(16.71,21.61)(19.43,16.02)
\psline{<-}(25,17.47)(15.6,30)
\psline(15.6,30)(5.46,17.14)
\psline{->}(9.35,22.09)(10.52,23.58)
\psdots[](15.58,29.91)
(15.58,29.91)
\psline{<-}(83.6,16.73)(74.62,29.87)
\psline(74.62,29.87)(64.94,16.39)
\psline{->}(68.65,21.57)(69.77,23.13)
\psdots[](74.61,29.77)
(74.61,29.77)
\psbezier{->}(74.47,29.77)(78.16,24.25)(83.88,14.48)(73.75,14.61)
\psbezier{->}(74.34,29.77)(78.29,22.53)(80.1,18.06)(73.78,18.06)
\psbezier(74.74,29.77)(70.79,23.03)(68.29,18.39)(73.86,18.06)
\psbezier(74.74,29.9)(70.69,23.26)(63.98,15.23)(73.98,14.61)
\psline{->}(77.5,4.64)(72.5,4.64)
\psline{->}(65,4.64)(65,9.64)
\psline{->}(85,4.64)(85,-0.36)
\psline[linewidth=0.1,linestyle=dashed,dash=1 1](60,4.64)(90,4.64)
\end{pspicture}
\caption{Mod\`eles hyperbolique et elliptique}
\label{f.modeles-secteurs}
\end{center}
\end{figure}
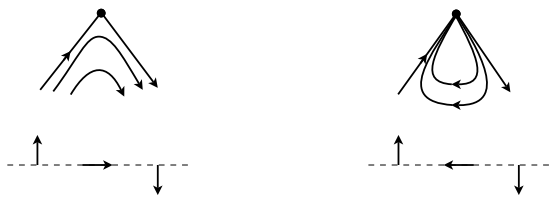
 Il nous reste \`a construire les quatre mod\`eles de secteurs \emph{indiff\'erents}, ce qui  est plus d\'elicat. On peut reprendre un exemple construit par M. Brown, E. Slaminka et W. Transue (\cite{BST}, \cite{brown2}).
On commence par d\'efinir une dynamique $g$ sur la bande $[-1, 1] \times \bbR$, en composant deux hom\'eomorphismes $\phi$ et $\psi$ (figure~\ref{f.modeles-secteurs-indifferents}). Pour cela, on consid\`ere l'ensemble $F$ constitu\'e des deux droites bordant la bande, $\{-1\} \times \bbR$ et $\{1\} \times \bbR$, ainsi que des points d'abscisse $0$ et d'ordonn\'ee enti\`ere. L'hom\'eomorphisme
 $\phi$ translate le point $(\theta,r)$ verticalement d'une hauteur $2\mid \theta \mid-1$ (on a donc $\phi(F) =F$). L'hom\'eomorphisme 
$\psi$ translate le point $(\theta,r)$ horizontalement, vers la droite, d'une distance \'egale \`a la moiti\'e de la distance euclidienne de $(\theta,r)$ \`a l'ensemble $F$ (il fixe donc chaque point de l'ensemble $F$). 
\begin{figure}[htbp]
\begin{center}
\def\JPicScale{1}
\ifx\JPicScale\undefined\def\JPicScale{1}\fi
\psset{unit=\JPicScale mm}
\psset{linewidth=0.3,dotsep=1,hatchwidth=0.3,hatchsep=1.5,shadowsize=1}
\psset{dotsize=0.7 2.5,dotscale=1 1,fillcolor=black}
\psset{arrowsize=1 2,arrowlength=1,arrowinset=0.25,tbarsize=0.7 5,bracketlength=0.15,rbracketlength=0.15}
\begin{pspicture}(0,0)(130,40)
\psline{->}(112.5,2.76)(117.5,2.76)
\psline{->}(105,2.5)(105,7.5)
\psline{->}(125,2.5)(125,7.76)
\psline[linewidth=0.1,linestyle=dashed,dash=1 1](100,2.76)(130,2.76)
\psdots[dotsize=0.7 0.5](114.54,11.25)
(114.54,16.89)(114.61,22.04)(114.67,27.69)(114.74,31.63)(114.67,33.65)(114.67,35.15)
\psbezier[linewidth=0.15,dotsize=0.7 0.5]{->}(115.26,21.74)(115.68,20.4)(115.72,19.01)(115.39,17.11)
\psbezier[linewidth=0.15,dotsize=0.7 0.5]{->}(115.39,16.49)(115.81,15.15)(115.85,13.76)(115.53,11.85)
\pscustom[linewidth=0.15,linestyle=dashed,dash=1 1,dotsize=0.7 0.5]{\psbezier{-}(105.86,12.77)(107.61,17.49)(108.88,20.73)(110.1,23.56)
\psbezier(110.1,23.56)(111.32,26.38)(112.28,27.16)(113.32,26.13)
\psbezier(113.32,26.13)(114.36,25.11)(115.32,23.73)(116.54,21.54)
\psbezier(116.54,21.54)(117.77,19.35)(118.41,19.06)(118.68,20.58)
\psbezier(118.68,20.58)(118.96,22.1)(118.81,23.87)(118.19,26.48)
\psbezier{->}(118.19,26.48)(117.57,29.1)(116.96,31.31)(116.18,33.85)
}
\psdots[](114.67,39.8)
(114.67,39.8)
\psline(114.75,39.69)(125.33,12.77)
\psline{->}(121.27,23.12)(120.05,26.25)
\psline(5,40)(5,10)
\psline(25,40)(25,10)
\psdots[linewidth=0.15,dotsize=0.7 1.5](15,40)
(15,30)(15,20)(15,10)
\psline[linewidth=0.15](7.5,40)(7.5,10)
\psline[linewidth=0.15](10,40)(10,10)
\psline[linewidth=0.15](12.5,40)(12.5,10)
\psline[linewidth=0.15](15,40)(15,10)
\psline[linewidth=0.15](17.5,40)(17.5,10)
\psline[linewidth=0.15](20,40)(20,10)
\psline[linewidth=0.15](22.5,40)(22.5,10)
\psline(114.69,40)(105,13.17)
\psline{->}(108.72,23.49)(109.84,26.6)
\psline(35,40)(35,10)
\psline(55,40)(55,10)
\psdots[linewidth=0.15,dotsize=0.7 1.5](45,40)
(45,30)(45,20)(45,10)
\psline(65,40)(65,10)
\psline(85,40)(85,10)
\psdots[linewidth=0.15,dotsize=0.7 1.5](75,40)
(75,30)(75,20)(75,10)
\psline[linewidth=0.15,dotsize=0.7 0.5]{->}(5,26.25)(5,28.75)
\psline[linewidth=0.15,dotsize=0.7 0.5]{->}(7.5,23.75)(7.5,26.25)
\psline[linewidth=0.15,dotsize=0.7 0.5]{->}(22.5,22.5)(22.5,25)
\psline[linewidth=0.15,dotsize=0.7 0.5]{->}(25,25)(25,27.5)
\psline[linewidth=0.15,dotsize=0.7 0.5]{->}(12.5,26.25)(12.5,23.75)
\psline[linewidth=0.15,dotsize=0.7 0.5]{->}(17.5,26.25)(17.5,23.75)
\psbezier[linewidth=0.15,dotsize=0.7 0.5]{->}(15.72,29.01)(16.32,26.62)(16.32,24.19)(15.72,20.92)
\psbezier[linewidth=0.15,dotsize=0.7 0.5]{->}(15.72,19.08)(16.32,16.69)(16.32,14.26)(15.72,10.99)
\psbezier[linewidth=0.15,dotsize=0.7 0.5]{->}(15.72,38.95)(16.32,36.56)(16.32,34.13)(15.72,30.86)
\psline[linewidth=0.15](35,10)(55,10)
\psline[linewidth=0.15](35,20)(55,20)
\psline[linewidth=0.15](35,17.5)(55,17.5)
\psline[linewidth=0.15](35,15)(55,15)
\psline[linewidth=0.15](35,12.5)(55,12.5)
\psline[linewidth=0.15,dotsize=0.7 0.5]{->}(38.75,20)(41.25,20)
\psline[linewidth=0.15,dotsize=0.7 0.5]{->}(47.5,20)(50,20)
\psline[linewidth=0.15,dotsize=0.7 0.5]{->}(42.5,17.5)(45,17.5)
\psline[linewidth=0.15,dotsize=0.7 0.5]{->}(43.75,15)(46.25,15)
\psline[linewidth=0.15,dotsize=0.7 0.5]{->}(42.5,12.5)(45,12.5)
\psline[linewidth=0.15,dotsize=0.7 0.5]{->}(42.5,22.5)(45,22.5)
\psline[linewidth=0.15,dotsize=0.7 0.5]{->}(42.5,32.5)(45,32.5)
\psline[linewidth=0.15](35,30)(55,30)
\psline[linewidth=0.15](35,27.5)(55,27.5)
\psline[linewidth=0.15](35,25)(55,25)
\psline[linewidth=0.15](35,22.5)(55,22.5)
\psline[linewidth=0.15,dotsize=0.7 0.5]{->}(38.75,30)(41.25,30)
\psline[linewidth=0.15,dotsize=0.7 0.5]{->}(47.5,30)(50,30)
\psline[linewidth=0.15,dotsize=0.7 0.5]{->}(42.5,27.5)(45,27.5)
\psline[linewidth=0.15,dotsize=0.7 0.5]{->}(43.75,25)(46.25,25)
\psline[linewidth=0.15](35,40)(55,40)
\psline[linewidth=0.15](35,37.5)(55,37.5)
\psline[linewidth=0.15](35,35)(55,35)
\psline[linewidth=0.15](35,32.5)(55,32.5)
\psline[linewidth=0.15,dotsize=0.7 0.5]{->}(38.75,40)(41.25,40)
\psline[linewidth=0.15,dotsize=0.7 0.5]{->}(47.5,40)(50,40)
\psline[linewidth=0.15,dotsize=0.7 0.5]{->}(42.5,37.5)(45,37.5)
\psline[linewidth=0.15,dotsize=0.7 0.5]{->}(43.75,35)(46.25,35)
\rput(15,6.25){$\phi$}
\rput(45,6.25){$\psi$}
\rput(16.25,35){}
\psbezier[linewidth=0.15,dotsize=0.7 0.5]{->}(75.72,38.88)(76.32,36.49)(76.32,34.06)(75.72,30.79)
\psbezier[linewidth=0.15,dotsize=0.7 0.5]{->}(75.72,19.14)(76.32,16.75)(76.32,14.32)(75.72,11.05)
\psbezier[linewidth=0.15,dotsize=0.7 0.5]{->}(75.72,28.88)(76.32,26.49)(76.32,24.06)(75.72,20.79)
\psline[linewidth=0.15,dotsize=0.7 0.5]{->}(65,25)(65,27.5)
\psline[linewidth=0.15,dotsize=0.7 0.5]{->}(85,25)(85,27.5)
\pscustom[linewidth=0.15,linestyle=dashed,dash=1 1,dotsize=0.7 0.5]{\psbezier{-}(66.25,10)(66.25,18.75)(66.81,24.38)(68.12,28.75)
\psbezier(68.12,28.75)(69.44,33.12)(70.75,33.5)(72.5,30)
\psbezier(72.5,30)(74.25,26.5)(75.75,23.5)(77.5,20)
\psbezier(77.5,20)(79.25,16.5)(80.56,16.88)(81.88,21.25)
\psbezier{->}(81.88,21.25)(83.19,25.62)(83.75,31.25)(83.75,40)
}
\rput(75,6.25){$\psi \circ \phi$}
\end{pspicture}
\caption{Mod\`ele indiff\'erent}
\label{f.modeles-secteurs-indifferents}
\end{center}
\end{figure}
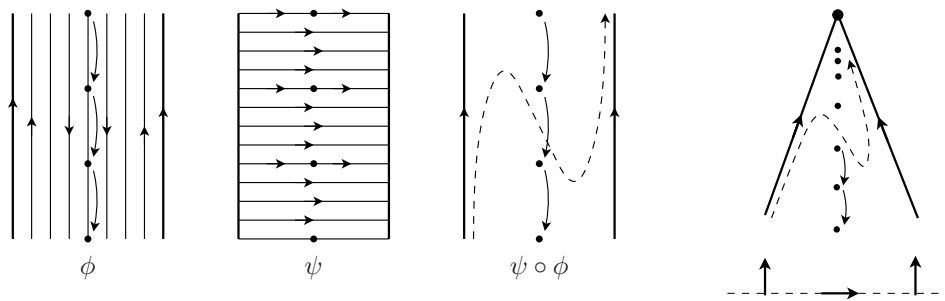
On peut ensuite identifier un secteur du plan \`a la bande $[-1, 1] \times \bbR$ par les coordonn\'ees polaires 
(par exemple $(\theta,r) \mapsto \exp(-r + i\pi \theta/4)$, ce qui
revient \`a ajouter \`a la bande le point \`a l'infini en haut), et transporter ainsi la dynamique $g$ sur un secteur du plan. Ceci nous fournit le secteur mod\`ele associ\'e au mot $(\ua \ra \ua)$. L'inverse donne le mot $(\da \la \da)$, et les deux derniers mod\`eles s'obtiennent en conjugant ces deux-ci par une sym\'etrie orthogonale.

Chaque bord de secteur-mod\`ele est une demi-droite issue de $0$, sur laquelle la dynamique est conjugu\'ee \`a l'homoth\'etie de rapport $2$ ou \`a son inverse, ce qui nous permet de recoller contin\^ument les diff\'erents types de secteurs, pour peu que la dynamique sur les bords \`a recoller se fasse dans le m\^eme sens.

\'Etant donn\'e un mot cyclique $\mathbf{M}$ autoris\'e de longueur $2d$, on \'ecrit la liste des $d$ sous-mots de longueur trois dont la premi\`ere lettre est une fl\`eche verticale, et on recolle dans l'ordre les secteurs mod\`eles associ\'es \`a ces sous-mots, de fa\c con \`a obtenir un hom\'eomorphisme $h$. On peut alors faire deux remarques. D'une part, une courbe qui relie les deux bords d'un secteur mod\`ele en le traversant n'est jamais libre. D'autre part, si $S$ est le secteur mod\`ele associ\'e au mot $(m_{1} m_{2} m_{3})$, si $x$ est un point du bord gauche de $S$ (correspondant \`a la fl\`eche $m_{1}$) et $y$ un point du bord droit (correspondant \`a la fl\`eche $m_{3}$), alors le couple $(x,y)$ est singulier pour $h$ si $m_{2}= \ra$, et il est singulier pour $h^{-1}$ si $m_{2}=\la$. Avec ces deux remarques,
il n'est pas difficile de v\'erifier que le mot  $\mathbf{M}(h)$ est \'egal \`a $\mathbf{M}$.
\end{demo}

 %%%%%%%%%%%%%%%%%%%%%%%%
%%%%%%%%%%%%%%%%%%%%%%%%
 %%%%%%%%%%%%%%%%%%%%%%%%
%%%%%%%%%%%%%%%%%%%%%%%%
 %%%%%%%%%%%%%%%%%%%%%%%%
%%%%%%%%%%%%%%%%%%%%%%%%
 %%%%%%%%%%%%%%%%%%%%%%%%
%%%%%%%%%%%%%%%%%%%%%%%%
 %%%%%%%%%%%%%%%%%%%%%%%%
%%%%%%%%%%%%%%%%%%%%%%%%
 
% \clearpage
 \part{Preuves des lemmes}
 \label{p.preuves-lemmes}
 %%%%%%%%%%%%%%%%%%%%%%%%
%%%%%%%%%%%%%%%%%%%%%%%%
Cette partie contient notamment les preuves des \'enonc\'es de la partie pr\'ec\'edente.
On commence par  caract\'eriser le module de $h$ comme le minimum des $h$-longueurs des courbes de Jordan qui entourent le point fixe (proposition~\ref{p.module-courbes},  section~\ref{s.module-courbes}). Ceci permet de  minorer le module en fonction de l'indice du point fixe (proposition~\ref{p.minoration-module}, section~\ref{s.minoration-module}). On en d\'eduit l'existence d'une suite g\'eod\'esique \'equivariante pour le relev\'e canonique $\wt h$, qui est essentiellement unique (corollaire~\ref{c.bouts-geodesique}, section~\ref{s.geodesiques}). 
Par ailleurs, nous montrons \`a la section~\ref{s.cRdB} que chaque composante de Reeb contient une droite de Brouwer qui s\'epare les deux bords ; ceci permet de d\'emontrer la formule d'indice pour le relev\'e $\wt h$, \`a l'aide de la notion d'indice partiel (th\'eor\`eme~\ref{t.lien-indice}-bis, section~\ref{s.lien-indice}). Nous terminons cette partie avec la proposition~\ref{p.compatibilite-composantes}, qui met en correspondance les composantes de Reeb de $\wt h$ et celles de $h$ (section~\ref{s.composantes-de-reeb}).

\section{Module et courbes}
\label{s.module-courbes}
%Dans cette section, on donne une preuve de la proposition~\ref{p.module-courbes}
%caract\'erisant le module comme la $h$-longueur des courbes de Jordan g\'eod\'esiques.
On se donne un hom\'eomorphisme $h$ du plan,
pr\'eservant l'orientation, fixant le point $0$, sans autre point fixe, et  tel que l'indice de
Poincar\'e-Lefschetz du point fixe est diff\'erent de $1$.
%de la premi\`ere partie 

\begin{demo}[de la proposition~\ref{p.module-courbes}]
\paragraph{Le module comme longueur minimale des courbes}
Il s'agit d'abord de montrer que le module, d\'efini en termes de domaines de translation (d\'efinition~\ref{d.module}), est \'egal au minimum des $h$-longueurs des courbes de degr\'e $1$.
Cette \'egalit\'e provient des deux remarques suivantes~:
\begin{itemize}
\item toute courbe libre est incluse dans un domaine de translation~;
\item si $x$ et $y$ appartiennent \`a un m\^eme domaine de translation,
alors il existe un arc libre allant de $x$ \`a $y$ (sauf si $x$ et $y$
sont dans la m\^eme orbite de $h$). 
\end{itemize}
Les d\'etails sont les m\^emes que pour la preuve du corollaire 3.5 du texte~\cite{I}.
%\end{demo}

%%%%%%%%%%%%%%%%%%%%%%%
%\subsection{Courbes de Jordan g\'eod\'esiques}
%\label{ss.jordan-geodesiques}
%Dans cette section, on consid\`ere un hom\'eomorphisme du plan $h$, pr\'eservant l'orientation, fixant uniquement le point $0$ (sans hypoth\`ese d'indice).
%Quand l'indice de Poincar\'e-Lefschetz du point fixe est diff\'erent de $1$, les notions de d\'ecomposition et de $h$-longueur d'une courbe ferm\'ee ont \'et\'e d\'efinies \`a la section~\ref{s.definition-courbes}. On adopte les m\^emes d\'efinitions lorsque l'indice est \'egal \`a $1$. La caract\'erisation du module donn\'ee (en indice diff\'erent de $1$) par la premi\`ere partie de la proposition~\ref{p.module-courbes} permet alors de g\'en\'eraliser la notion de module. 
%\begin{defi}
%Le module de $h$ est le 
%minimum des $h$-longueurs 
%des courbes ferm\'ees $\gamma$ de degr\'e $1$.
%\end{defi}
%On reprend aussi la d\'efinition des courbes g\'eod\'esiques.
%On peut alors prouver un \'enonc\'e qui g\'en\'eralise la deuxi\`eme partie de la proposition~\ref{p.module-courbes}.
%\begin{prop}
%\label{p.jordan-geodesique}
%Soit $h$ un hom\'eomorphisme du plan, pr\'eservant l'orientation, fixant 
%uniquement le point $0$.
%Alors il existe une courbe de Jordan g\'eod\'esique pour $h$.
%\end{prop}

%\begin{demo}[de la proposition~\ref{p.jordan-geodesique}]
\paragraph{Courbe de Jordan g\'eod\'esique}
Passons maintenant \`a la preuve de la seconde partie de la proposition~\ref{p.module-courbes}.
Il est classique de construire une courbe de Jordan \`a partir d'une
courbe contenant des points doubles en ``effa\c{c}ant les boucles''~; il
faut juste faire attention \`a obtenir une courbe de degr\'e $1$, et
surtout \`a conserver le caract\`ere g\'eod\'esique.
On part d'une courbe g\'eod\'esique $\gamma$, que l'on va modifier pour
obtenir une courbe de Jordan g\'eod\'esique.
On voit $\gamma$ comme une application du cercle $\bbS^1 = \bbR/\bbZ$
dans $\bbR^2 \setminus \{0\}$.
%\pourmoi{INUTILE ? Tout d'abord, toute courbe ferm\'ee de degr\'e $>1$ contient une sous- courbe
%ferm\'e\footnote{def ?} de degr\'e $1$ (...). Quitte \`a remplacer $\gamma$ par une de ses
%sous-courbe ferm\'ee , on peut dont supposer que $\gamma$ ne contient aucune
%sous-courbe ferm\'ee de degr\'e diff\'erent de $0$ ou $1$.}
 Remarquons   que toute
courbe assez proche d'une courbe libre est encore libre~; par
cons\'equent, toute courbe assez voisine de $\gamma$ a une $h$-longueur
inf\'erieure ou \'egale \`a celle de $\gamma$. D'autre part, comme $\gamma$
est une courbe g\'eod\'esique, elle est d\'ej\`a de $h$-longueur minimale~:
finalement, toute courbe  assez voisine de $\gamma$ a une $h$-longueur
 \'egale \`a celle de $\gamma$.
Quitte \`a perturber $\gamma$ arbitrairement peu, on peut supposer que
$\gamma$  n'a qu'un nombre fini de points
multiples (tout en restant une courbe g\'eod\'esique d'apr\`es ce qui
pr\'ec\`ede).  
Disons qu'un intervalle ferm\'e non trivial $[t t']$ du cercle est une
\emph{boucle de $\gamma$} (de degr\'e $0$) si $\gamma(t)=\gamma(t')$, et si la
courbe  $[t t']_\gamma$  est une courbe ferm\'ee de degr\'e $0$ dans
$\bbR^2 \setminus\{0\}$.
 Une boucle $[t t']$  est \emph{maximale} si 
 il n'existe pas d'intervalle contenant strictement $[t t']$ qui soit encore une boucle de $\gamma$. Puisque $\gamma$ n'a qu'un nombre fini de   points multiples, elle n'a
qu'un nombre fini de boucles, et toute boucle est incluse dans une boucle maximale (non unique).
 
 Consid\'erons une famille $\cF$ d'intervalles ferm\'es du cercle, deux \`a deux disjoints, qui sont des boucles maximales de $\gamma$. On peut alors d\'efinir
\begin{itemize}
\item un espace $\overline{\bbS^1}$, hom\'eomorphe au cercle, obtenu en identifiant chaque
intervalle  $[t_i t'_i]$ de $\cF$ \`a un point (on note $p : \bbS^1 \rightarrow
\overline{\bbS^1}$ l'application quotient)~;
\item une application $\gamma_{\cF} : \overline{\bbS^1} \ra \bbR^2 \setminus \{0 \}$ en posant
$$
\left\{
\begin{array}{rcll}
{\gamma_{\cF}} \circ p ([t_{i}t'_{i}]) & = &\gamma(t_{i}) = \gamma(t'_{i}) & \mbox{ pour chaque intervalle } [t_{i} t'_{i}] \in \cF \ ; \\
{\gamma_{\cF}} \circ p (t) & = &\gamma(t) & \mbox{ en dehors de  } \cup\{ I \mid I \in \cF \}.
\end{array}
\right.
$$
\end{itemize}
Cette application est encore une courbe de degr\'e $1$ dans $\bbR^2 \setminus \{ 0 \}$.
Puisque $\gamma$ n'a qu'un nombre fini de boucles, on peut choisir une telle  famille $\cF$  qui est maximale, autrement dit, telle qu'il n'existe aucun intervalle $J$, qui est une boucle maximale, et qui est disjoint de tous les intervalles de $\cF$. On voit facilement que dans ce cas, $\gamma_{\cF}$ est une courbe qui n'admet aucune boucle. On peut alors utiliser le lemme suivant.
\footnote{Voici un sch\'ema de preuve (les d\'etails sont laiss\'es au lecteur).
On consid\`ere les \emph{boucles de degr\'e $p$} de $\alpha$ (avec une d\'efinition \'evidente) :
il s'agit de montrer qu'il n'y en a pas.
Par hypoth\`ese $\alpha$ n'a pas de boucle de degr\'e $0$. Elle n'a pas non plus de boucle $\alpha'$ de degr\'e $1$, sans quoi on aurait $\alpha=\alpha' * \alpha''$ avec $\alpha''$ boucle de degr\'e $0$.
Consid\'erons alors une boucle $\beta$ de degr\'e $p>1$, qui est minimale : c'est une courbe de Jordan. On conclut avec le lemme classique :
\begin{lemm}
Il n'existe pas de courbe de Jordan $\beta$
% : \bbS^1 \ra \bbR^2 \setminus \{0\}$ 
de degr\'e $p>1$.
\end{lemm}}
%Pour cela,  on consid\`ere le relev\'e $B : \bbR \ra \bbR^2$ d'une telle courbe $\beta$. Les hypoth\`eses sur $\beta$ disent que  l'ensemble $B([0,1])$ est disjoint de $B([1,2]), qui est son image par la translation $\tau$. D'apr\`es le lemme de Franks appliqu\'e \`a l'hom\'eomorphisme de Brouwer $\tau$, on a $\tau^p(B([0,1])) \cap B([0,1]) = \emptyset$. Ceci contredit $B(1) = \tau^p(B(0))$.}
\begin{lemm}
Si $\alpha$ est une courbe ferm\'ee de degr\'e $1$ dans $\bbR^2 \setminus \{0\}$,
et si elle n'a pas de boucle, alors c'est une courbe de Jordan.
\end{lemm}
La d\'efinition de $\gamma_{\cF}$ entra\^ine que pour tout intervalle $[ss']$ du
cercle $\bbS^1$,   ${\gamma_{\cF}} \circ p([ss'])$ est inclus dans
$\gamma([ss'])$.
En particulier, si $\gamma([ss'])$ est libre, alors
${\gamma_{\cF}}\circ p ([ss'])$ est encore libre, ce qui signifie que la
$h$-longueur de ${\gamma_{\cF}}$ est inf\'erieure ou \'egale \`a celle de
$\gamma$. Comme $\gamma$ \'etait d\'ej\`a g\'eod\'esique, ${\gamma_{\cF}}$
est aussi une courbe g\'eod\'esique. On a ainsi obtenu une courbe de Jordan g\'eod\'esique.
\end{demo}

%%%%%%%%%%%%%%%%%%%%%%%%
%%%%%%%%%%%%%%%%%%%%%%%%
%\section{Relev\'e canonique}

%A ECRIRE.

%%%%%%%%%%%%%%%%%%%%%%%%
%%%%%%%%%%%%%%%%%%%%%%%%
\section{Minoration du module}
\label{s.minoration-module}
Dans cette section, on d\'emontre la proposition~\ref{p.minoration-module}, qui permet de minorer le module de $h$ en fonction de l'indice du point fixe : plus l'indice est diff\'erent de $1$, plus le module doit \^etre grand. Pour cela, \emph{nous appliquons une bonne partie des r\'esultats du texte~\cite{0}.}
%\pourmoi{??? Une cons\'equence imm\'ediate est que le module de $h$ est sup\'erieur ou \'egal \`a $4$ d\`es que l'indice du point fixe est diff\'erent de $0$, $1$ ou $2$.
%La preuve des r\'esultats de la section~\ref{s.definitions} n'utilise que cette cons\'equence (c'est ce qui est expliqu\'e dans l'appendice~\ref{}).
%On peut aussi en donner une  preuve \'el\'ementaire, ind\'ependante du texte~\cite{0}.}

%%%%%%%%%%
\subsection{D\'ecompositions en briques}
\subsubsection*{G\'en\'eralit\'es}
Nous commen\c{c}ons par rappeler bri\`evement la d\'efinition des d\'ecompositions en briques (voir les travaux de P. Le Calvez et A. Sauzet (\cite{sauzet,lecasau}) ou \cite{0} pour des motivations et des d\'efinitions plus
d\'etaill\'ees).
 Soit $H$ un hom\'eomorphisme du
plan, pr\'eservant l'orientation, 
et n'ayant qu'un nombre fini de points fixes (en pratique, dans ce
texte, $H$ aura $0$ ou $1$ point fixe).
On pose $U=\bbR^2 \setminus \fixe(H)$, o\`u $\fixe(H)$ d\'esigne
l'ensemble des points fixes de $H$.
Soit $\cF$ un graphe triadique plong\'e dans $U$. L'adh\'erence (dans $U$)
d'une composante connexe  du compl\'ementaire de $\cF$ est appel\'e
\emph{brique}. D'un point de vue topologique, on suppose que
chaque  brique est compacte. D'un point de vue dynamique, on demande
que chaque brique soit libre par $H$, mais que la r\'eunion de deux briques
adjacentes ne soit pas libre. Si toutes ces hypoth\`eses sont
v\'erifi\'ees, l'ensemble $\cF$ est appel\'e \emph{d\'ecomposition
en briques} pour $H$.

\subsubsection*{D\'ecomposition g\'eod\'esique}
On suppose maintenant que  $h$ est un hom\'eomorphisme du plan, pr\'eservant l'orientation, dont  $0$ est l'unique point fixe, avec un indice diff\'erent de $1$.
%\pourmoi{??? Rappelons que le module est d\'efini sans
%condition d'indice (\`a la section~\ref{s.jordan-courbes}).}
 Soit $\cF$ une d\'ecomposition en briques pour $h$. Le
\emph{module} de $\cF$ est le nombre minimal de briques dont la
r\'eunion s\'epare $0$ et $\infty$ dans la sph\`ere  $\bbR^2 \cup
\{\infty\}$
(autrement dit, dont la r\'eunion contient une courbe ferm\'ee de
degr\'e $1$ dans $\bbR^2 \setminus \{0\}$).
\begin{defi}
Une d\'ecomposition en brique $\cF$ pour $h$ est
\emph{g\'eod\'esique} si son module est \'egal au module
de $h$. Elle est \emph{compatible} avec une courbe g\'eod\'esique ferm\'ee $\gamma$ si 
l'intersection $\gamma \cap \cF$ co\"incide avec l'ensemble des sommets d'une d\'ecomposition de $\gamma$.
\end{defi}

\begin{lemm}\label{l.decomposition-geodesique}
%On suppose que le module de $h$ est sup\'erieur ou \'egal \`a deux.
Pour toute courbe de Jordan g\'eod\'esique $\gamma$, il existe une d\'ecom\-po\-sition en briques, g\'eod\'esique, compatible avec $\gamma$.
\end{lemm}

\begin{demo}[du lemme~\ref{l.decomposition-geodesique}]
Notons $d$ le module de $h$. Commen\c cons par remarquer que si on a une d\'ecomposition en briques de module $d'$, comme chaque brique est libre,  il existe une courbe de Jordan de $h$-longueur \'egale \`a $d'$. Ceci prouve que $d' \geq d$. Il reste donc \`a trouver une d\'ecomposition dont le module satisfasse l'in\'egalit\'e inverse.
%D'apr\`es la proposition~\ref{p.jordan-geodesique}, il existe une
%courbe de Jordan $\gamma$ qui est une courbe ferm\'ee
%g\'eod\'esique. 
Soit maintenant $\gamma$ une courbe de Jordan g\'eod\'esique.
En consid\'erant une d\'ecomposition de $\gamma$, et en
\'epaississant $\gamma$ en un anneau, on trouve un graphe triadique $\cF$
tel que (voir figure~\ref{f.decomposition-geodesique})
\begin{itemize}
\item le compl\'ementaire de $\cF$ a exactement $d+2$ composantes
connexes~;
\item la premi\`ere contient $0$, la deuxi\`eme est non born\'ee~;
\item les adh\'erences des $d$ composantes restantes, que l'on note
$B_1, \dots , B_d$ sont libres, et
leur r\'eunion contient $\gamma$ (voir figure~\ref{f.decomposition-geodesique}).
\end{itemize}
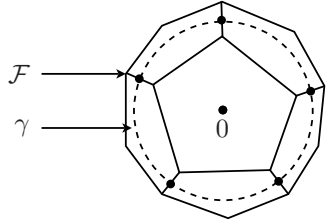
\begin{figure}[htbp]
\begin{center}
\def\JPicScale{0.8}
\ifx\JPicScale\undefined\def\JPicScale{1}\fi
\psset{unit=\JPicScale mm}
\psset{linewidth=0.3,dotsep=1,hatchwidth=0.3,hatchsep=1.5,shadowsize=1}
\psset{dotsize=0.7 2.5,dotscale=1 1,fillcolor=black}
\psset{arrowsize=1 2,arrowlength=1,arrowinset=0.25,tbarsize=0.7 5,bracketlength=0.15,rbracketlength=0.15}
\begin{pspicture}(0,0)(55,36.03)
\psbezier(38.14,17.99)(38.14,17.99)(38.14,17.99)(38.14,17.99)
\psdots[](38.14,17.99)
(38.14,17.99)
\rput{0}(38.18,17.82){\psellipse[linestyle=dashed,dash=1 1](0,0)(15,15)}
\psdots[](37.93,32.93)
(24.25,23.18)(29.43,5.68)(47.37,6.13)(52.73,21.12)
\rput(38,14.84){$0$}
\pspolygon(38,30.24)(26.55,22.21)
(26.55,22.21)(30.89,7.61)
(30.89,7.61)(46,8)
(46,8)(50.37,20.37)
(50.37,20.37)(38,30.24)
\pspolygon(38,36)(28,32)
(28,32)(22,24)
(22,24)(22,12)
(22,12)(28,4)
(28,4)(39,0)
(39,0)(49,4)
(49,4)(54,12)
(54,12)(55,22)
(55,22)(50,31)
(50,31)(38,36)
\psline(37.87,36.03)(38.07,30.17)
\psline(26.49,22.21)(22.14,24.12)
\psline(50.37,20.37)(54.97,21.95)
\psline(45.96,8)(48.99,4.05)
\psline(28,4.05)(30.89,7.61)
\psline{->}(8,15)(23,15)
\psline{->}(8,24)(22,24)
\rput[r](6,24){$\cal F$}
\rput[r](6,15){$\gamma$}
\rput[r](52,0){}
\end{pspicture}
\caption{\'Epaississement d'une courbe de Jordan g\'eod\'esique}
\label{f.decomposition-geodesique}
\end{center}
\end{figure}
D'apr\`es le corollaire~{3.63} de~\cite{0}, il existe une d\'ecomposition
en briques $\cF'$ qui ``\'etend'' $\cF$, au sens o\`u chaque ensemble $B_1,
\dots B_d$ est inclus dans une brique de $\cF'$.
%Puisque la courbe $\gamma$ \'etait g\'eod\'esique, la r\'eunion de deux brique $B_{i}$ adjacentes et $B_{i+1}$ ne peut pas \^etre libre ; par cons\'equent 
Il est clair que le module $d'$ de $\cF'$ v\'erifie $d' \leq d$, ce que l'on voulait.
\end{demo}

%\footnote{def compatible ?? Si l'intersection $\gamma \cap \cF$ co\"incide avec l'ensemble des sommets d'une d\'ecomposition de $\gamma$, comme dans la preuve ci-dessus, 
%on dira que $\cF$ et $\gamma$ sont \emph{compatibles}. On voit facilement que pour toute d\'ecomposition g\'eod\'esique $\cF$, il existe une courbe de Jordan  g\'eod\'esique $\gamma$ compatible avec $\cF$.
%}

%%%%%%%%%%%%%%%%%%%%%%%
\subsection{Preuve de la minoration}
On consid\`ere une d\'ecomposition en briques $\cal F$ qui est g\'eod\'esique pour $h$.
Supposons d'abord que l'indice du point fixe soit plus petit que $1$, on \'ecrit
$$
\indi(h,0) = 1-p < 1.
$$
On identifie $\bbR^2 \cup \{\infty\}$ \`a la sph\`ere $\bbS^2$, en envoyant les points $0$ et $\infty$ respectivement sur deux  points $N$ et $S$ de la sph\`ere.
On peut alors appliquer le r\'esultat principal du texte~\cite{0} (th\'eor\`eme F, p80) : il existe $2p$ croissants attractifs et r\'epulsifs simpliciaux, d'int\'erieurs deux \`a deux disjoints. La seule chose qui nous importe ici est que toute famille de briques de $\cF$ dont la r\'eunion s\'epare $0$ et $\infty$ doit rencontrer chaque croissant ; par cons\'equent, une telle famille doit contenir au moins $2p$ briques. Ceci prouve que le module de la d\'ecomposition $\cal F$ est sup\'erieur ou \'egal \`a $2p$. Comme la d\'ecomposition est g\'eod\'esique, on en d\'eduit la m\^eme minoration pour  le module de $h$.
Dans le cas o\`u l'indice est $1+p>1$, on inverse les r\^oles de $0$ et $\infty$ (ou, ce qui revient au m\^eme, on conjugue $h$ par l'application $z \mapsto -\frac{1}{z}$). D'apr\`es le th\'eor\`eme de Lefschetz, l'indice du point $\infty$ vaut alors $1-p<1$, ce qui nous ram\`ene au cas pr\'ec\'edent.

%%%%%%%%%%%%%%%%%%%%%%%%
%%%%%%%%%%%%%%%%%%%%%%%%
\section{G\'eod\'esiques}
\label{s.geodesiques}
Dans cette section, on d\'emontre qu'une g\'eod\'esique locale est une g\'eod\'esique globale (proposition~\ref{p.geodesique-locale}). On en d\'eduit l'existence d'une suite g\'eod\'esique pour le relev\'e canonique $\wt h$, qui est $\tau$-\'equivariante (corollaire~\ref{c.bouts-geodesique}).
\emph{Cette section utilise de mani\`ere essentielle les r\'esultats de l'article~\cite{I}}.

\begin{figure}[htbp]
\begin{center}
\def\JPicScale{0.8}
\ifx\JPicScale\undefined\def\JPicScale{1}\fi
\psset{unit=\JPicScale mm}
\psset{linewidth=0.3,dotsep=1,hatchwidth=0.3,hatchsep=1.5,shadowsize=1}
\psset{dotsize=0.7 2.5,dotscale=1 1,fillcolor=black}
\psset{arrowsize=1 2,arrowlength=1,arrowinset=0.25,tbarsize=0.7 5,bracketlength=0.15,rbracketlength=0.15}
\begin{pspicture}(0,0)(130,60)
\psline{->}(0,30)(130,30)
\pscustom[linewidth=0.5]{\psbezier{-}(50,0)(53.5,14)(58,20)(65,20)
\psbezier{->}(65,20)(72,20)(76.5,14)(80,0)
}
\pscustom[]{\psbezier{-}(15,0)(18.5,10.5)(21.5,15)(25,15)
\psbezier{->}(25,15)(28.5,15)(31.5,10.5)(35,0)
}
\pscustom[]{\psbezier{-}(60,0)(62.1,5.6)(63.75,8)(65.5,8)
\psbezier{->}(65.5,8)(67.25,8)(68.6,5.6)(70,0)
}
\pscustom[linewidth=0.5]{\psbezier{-}(10,0)(13.5,14)(18,20)(25,20)
\psbezier{->}(25,20)(32,20)(36.5,14)(40,0)
}
\pscustom[]{\psbezier{-}(55,0)(58.5,10.5)(61.5,15)(65,15)
\psbezier{->}(65,15)(68.5,15)(71.5,10.5)(75,0)
}
\pscustom[]{\psbezier{-}(20,0)(22.1,5.6)(23.75,8)(25.5,8)
\psbezier{->}(25.5,8)(27.25,8)(28.6,5.6)(30,0)
}
\pscustom[linewidth=0.5]{\psbezier{-}(90,0)(93.5,14)(98,20)(105,20)
\psbezier{->}(105,20)(112,20)(116.5,14)(120,0)
}
\pscustom[]{\psbezier{-}(95,0)(98.5,10.5)(101.5,15)(105,15)
\psbezier{->}(105,15)(108.5,15)(111.5,10.5)(115,0)
}
\pscustom[]{\psbezier{-}(100,0)(102.1,5.6)(103.75,8)(105.5,8)
\psbezier{->}(105.5,8)(107.25,8)(108.6,5.6)(110,0)
}
\pscustom[linewidth=0.5]{\psbezier{-}(90,60)(93.5,46)(98,40)(105,40)
\psbezier{->}(105,40)(112,40)(116.5,46)(120,60)
}
\pscustom[]{\psbezier{-}(95,60)(98.5,49.5)(101.5,45)(105,45)
\psbezier{->}(105,45)(108.5,45)(111.5,49.5)(115,60)
}
\pscustom[]{\psbezier{-}(100,60)(102.1,54.4)(103.75,52)(105.5,52)
\psbezier{->}(105.5,52)(107.25,52)(108.6,54.4)(110,60)
}
\pscustom[linewidth=0.5]{\psbezier{-}(50,60)(53.5,46)(58,40)(65,40)
\psbezier{->}(65,40)(72,40)(76.5,46)(80,60)
}
\pscustom[]{\psbezier{-}(55,60)(58.5,49.5)(61.5,45)(65,45)
\psbezier{->}(65,45)(68.5,45)(71.5,49.5)(75,60)
}
\pscustom[]{\psbezier{-}(60,60)(62.1,54.4)(63.75,52)(65.5,52)
\psbezier{->}(65.5,52)(67.25,52)(68.6,54.4)(70,60)
}
\pscustom[linewidth=0.5]{\psbezier{-}(10,60)(13.5,46)(18,40)(25,40)
\psbezier{->}(25,40)(32,40)(36.5,46)(40,60)
}
\pscustom[]{\psbezier{-}(15,60)(18.5,49.5)(21.5,45)(25,45)
\psbezier{->}(25,45)(28.5,45)(31.5,49.5)(35,60)
}
\pscustom[]{\psbezier{-}(20,60)(22.1,54.4)(23.75,52)(25.5,52)
\psbezier{->}(25.5,52)(27.25,52)(28.6,54.4)(30,60)
}
\pscustom[]{\psbezier{-}(0,44)(2.1,48.2)(3.6,50)(5,50)
\psbezier(5,50)(6.4,50)(8.95,47.6)(13.5,42)
\psbezier(13.5,42)(18.05,36.4)(21.5,34)(25,34)
\psbezier(25,34)(28.5,34)(31.95,36.4)(36.5,42)
\psbezier(36.5,42)(41.05,47.6)(43.6,50)(45,50)
\psbezier(45,50)(46.4,50)(48.95,47.6)(53.5,42)
\psbezier(53.5,42)(58.05,36.4)(61.5,34)(65,34)
\psbezier(65,34)(68.5,34)(71.95,36.4)(76.5,42)
\psbezier(76.5,42)(81.05,47.6)(83.6,50)(85,50)
\psbezier(85,50)(86.4,50)(88.95,47.6)(93.5,42)
\psbezier(93.5,42)(98.05,36.4)(101.5,34)(105,34)
\psbezier(105,34)(108.5,34)(111.5,35.65)(115,39.5)
\psbezier{->}(115,39.5)(118.5,43.35)(120.9,45.9)(123,48)
}
\pscustom[]{\psbezier{-}(0,51)(2.1,57.3)(3.75,60)(5.5,60)
\psbezier(5.5,60)(7.25,60)(8.75,57.6)(10.5,52)
\psbezier(10.5,52)(12.25,46.4)(14.05,43.1)(16.5,41)
\psbezier(16.5,41)(18.95,38.9)(21.5,38)(25,38)
\psbezier(25,38)(28.5,38)(31.05,38.9)(33.5,41)
\psbezier(33.5,41)(35.95,43.1)(37.9,46.4)(40,52)
\psbezier(40,52)(42.1,57.6)(43.6,60)(45,60)
\psbezier(45,60)(46.4,60)(47.9,57.6)(50,52)
\psbezier(50,52)(52.1,46.4)(54.05,43.1)(56.5,41)
\psbezier(56.5,41)(58.95,38.9)(61.5,38)(65,38)
\psbezier(65,38)(68.5,38)(71.05,38.9)(73.5,41)
\psbezier(73.5,41)(75.95,43.1)(77.9,46.4)(80,52)
\psbezier(80,52)(82.1,57.6)(83.6,60)(85,60)
\psbezier(85,60)(86.4,60)(87.9,57.6)(90,52)
\psbezier(90,52)(92.1,46.4)(94.05,43.1)(96.5,41)
\psbezier(96.5,41)(98.95,38.9)(101.5,38)(105,38)
\psbezier(105,38)(108.5,38)(111.05,38.9)(113.5,41)
\psbezier(113.5,41)(115.95,43.1)(117.6,45.5)(119,49)
\psbezier{->}(119,49)(120.4,52.5)(121.6,54.9)(123,57)
}
\pscustom[]{\psbezier{-}(0,16)(2.1,11.8)(3.6,10)(5,10)
\psbezier(5,10)(6.4,10)(8.95,12.4)(13.5,18)
\psbezier(13.5,18)(18.05,23.6)(21.5,26)(25,26)
\psbezier(25,26)(28.5,26)(31.95,23.6)(36.5,18)
\psbezier(36.5,18)(41.05,12.4)(43.6,10)(45,10)
\psbezier(45,10)(46.4,10)(48.95,12.4)(53.5,18)
\psbezier(53.5,18)(58.05,23.6)(61.5,26)(65,26)
\psbezier(65,26)(68.5,26)(71.95,23.6)(76.5,18)
\psbezier(76.5,18)(81.05,12.4)(83.6,10)(85,10)
\psbezier(85,10)(86.4,10)(88.95,12.4)(93.5,18)
\psbezier(93.5,18)(98.05,23.6)(101.5,26)(105,26)
\psbezier(105,26)(108.5,26)(111.5,24.35)(115,20.5)
\psbezier{->}(115,20.5)(118.5,16.65)(120.9,14.1)(123,12)
}
\pscustom[]{\psbezier{-}(0,9)(2.1,2.7)(3.75,0)(5.5,0)
\psbezier(5.5,0)(7.25,0)(8.75,2.4)(10.5,8)
\psbezier(10.5,8)(12.25,13.6)(14.05,16.9)(16.5,19)
\psbezier(16.5,19)(18.95,21.1)(21.5,22)(25,22)
\psbezier(25,22)(28.5,22)(31.05,21.1)(33.5,19)
\psbezier(33.5,19)(35.95,16.9)(37.9,13.6)(40,8)
\psbezier(40,8)(42.1,2.4)(43.6,0)(45,0)
\psbezier(45,0)(46.4,0)(47.9,2.4)(50,8)
\psbezier(50,8)(52.1,13.6)(54.05,16.9)(56.5,19)
\psbezier(56.5,19)(58.95,21.1)(61.5,22)(65,22)
\psbezier(65,22)(68.5,22)(71.05,21.1)(73.5,19)
\psbezier(73.5,19)(75.95,16.9)(77.9,13.6)(80,8)
\psbezier(80,8)(82.1,2.4)(83.6,0)(85,0)
\psbezier(85,0)(86.4,0)(87.9,2.4)(90,8)
\psbezier(90,8)(92.1,13.6)(94.05,16.9)(96.5,19)
\psbezier(96.5,19)(98.95,21.1)(101.5,22)(105,22)
\psbezier(105,22)(108.5,22)(111.05,21.1)(113.5,19)
\psbezier(113.5,19)(115.95,16.9)(117.6,14.5)(119,11)
\psbezier{->}(119,11)(120.4,7.5)(121.6,5.1)(123,3)
}
\pscustom[]{\psbezier{-}(0,39)(1.4,39.7)(3.2,40)(6,40)
\psbezier(6,40)(8.8,40)(11.5,38.8)(15,36)
\psbezier(15,36)(18.5,33.2)(21.5,32)(25,32)
\psbezier(25,32)(28.5,32)(31.5,33.2)(35,36)
\psbezier(35,36)(38.5,38.8)(41.5,40)(45,40)
\psbezier(45,40)(48.5,40)(51.5,38.8)(55,36)
\psbezier(55,36)(58.5,33.2)(61.5,32)(65,32)
\psbezier(65,32)(68.5,32)(71.5,33.2)(75,36)
\psbezier(75,36)(78.5,38.8)(81.5,40)(85,40)
\psbezier(85,40)(88.5,40)(91.5,38.8)(95,36)
\psbezier(95,36)(98.5,33.2)(101.5,32)(105,32)
\psbezier(105,32)(108.5,32)(111.5,33.2)(115,36)
\psbezier{->}(115,36)(118.5,38.8)(121.5,40)(125,40)
}
\pscustom[]{\psbezier{-}(0,21)(1.4,20.3)(3.2,20)(6,20)
\psbezier(6,20)(8.8,20)(11.5,21.2)(15,24)
\psbezier(15,24)(18.5,26.8)(21.5,28)(25,28)
\psbezier(25,28)(28.5,28)(31.5,26.8)(35,24)
\psbezier(35,24)(38.5,21.2)(41.5,20)(45,20)
\psbezier(45,20)(48.5,20)(51.5,21.2)(55,24)
\psbezier(55,24)(58.5,26.8)(61.5,28)(65,28)
\psbezier(65,28)(68.5,28)(71.5,26.8)(75,24)
\psbezier(75,24)(78.5,21.2)(81.5,20)(85,20)
\psbezier(85,20)(88.5,20)(91.5,21.2)(95,24)
\psbezier(95,24)(98.5,26.8)(101.5,28)(105,28)
\psbezier(105,28)(108.5,28)(111.5,26.8)(115,24)
\psbezier{->}(115,24)(118.5,21.2)(121.5,20)(125,20)
}
\pscustom[linewidth=0.2,linestyle=dashed,dash=1 1]{\psline(25,10)(25,50)
\psline(25,50)(65,10)
\psline(65,10)(65,50)
\psline(65,50)(105,10)
\psline(105,10)(105,50)
\psbezier(105,50)(105,50)(105,50)(105,50)
}
\psdots[linewidth=0.2,linestyle=dashed,dash=1 1,dotsize=0.7 3.5](25,10)
(25,50)(65,10)(65,50)(105,10)(105,50)(105,50)
\psline[linewidth=0.2,linestyle=dashed,dash=1 1](25,10)(10,25)
\psline[linewidth=0.2,linestyle=dashed,dash=1 1](105,50)(120,35)
\end{pspicture}
\caption{Une g\'eod\'esique locale d'ordre $2$ qui n'est pas une g\'eod\'esique globale}
\label{f.geo-locale}
\end{center}
\end{figure}
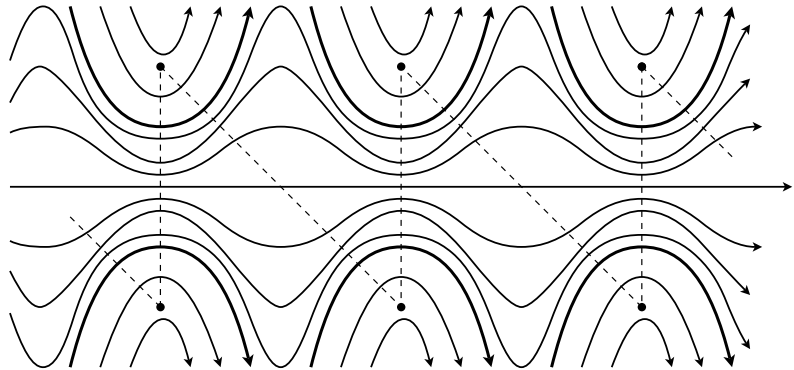

Remarquons que la proposition~\ref{p.geodesique-locale} est fausse pour $d=2$ : la figure~\ref{f.geo-locale} repr\'esente un hom\'eo\-mor\-phisme de Brouwer ayant une g\'eod\'esique locale d'ordre $2$ qui n'est pas une g\'eod\'esique globale.

%%%%%%%%%%%%%%%%%%%%%%%
\subsection{G\'eod\'esiques locales et globales}
On se donne un hom\'eomorphisme de Brouwer $H$.
On commence par une remarque \'el\'ementaire.
\begin{rema}\label{r.geo-locale}
Consid\'erons $(x_i)_{i  \in I}$  une suite g\'eod\'esique
locale d'ordre $1$ (autrement dit, deux points cons\'ecutifs sont \`a
$H$-distance $0$ ou $1$). Soit $d \geq 3$. 
Pour que cette suite soit une suite g\'eod\'esique
locale d'ordre $d$, il suffit qu'on ait, pour tout entier $i$, $d_H(x_i,x_{i+d})=d$ (autrement dit, 
 il suffit que l'implication de la d\'efinition~\ref{d.geodesique-locale} soit
v\'erifi\'ee quand $|i-j|=d$).
 Cette remarque provient simplement de l'in\'egalit\'e triangulaire.
\end{rema}

La preuve de la proposition~\ref{p.geodesique-locale} repose sur l'existence des composantes de
Reeb (th\'eor\`eme 2 de \cite{I}), et sur le lemme suivant.
\begin{lemm}[``d\'eviation des g\'eod\'esiques'']
\label{l.deviation-geodesiques}
Soit $(x,y)$ un couple de points singulier pour $H$ ou pour $H^{-1}$.
On consid\`ere deux autres points $x'$ et $y'$.
 
Si $\displaystyle \left\{\begin{array}{l} d_H(x',x) \leq 1 \\
d_H(y,y') \leq 1 \end{array}\right.$ alors $d_H(x',y') \leq 2$.
\end{lemm}
\begin{demo}[du lemme]
L'hypoth\`ese $d_H(x',x) \leq 1$ signifie qu'il existe un domaine de
translation $O_x$ contenant $x$ et $x'$. De m\^eme, si $d_H(y,y') \leq 1$,
il existe un domaine de translation $O_y$ contenant $y$ et
$y'$. Puisque le couple $(x,y)$ est singulier pour $H$ ou $H^{-1}$, il
existe un point $z$ de $O_x$ qui a un it\'er\'e $H^n(z)$ dans $O_y$.
Le r\'esultat d\'ecoule alors de l'in\'egalit\'e triangulaire et de
l'invariance de $d_H$ (lemme 3.7 de \cite{I})~:
$$
\begin{array}{rcl}
d_H(x',y') &\leq & d_H(x',z) + d_H(z,y') \\
	& \leq 	& d_H(x',z) + d_H(H^n(z),y') \\
	& \leq & 1+1 =2.
\end{array}
$$
\end{demo}

Avant d'entamer la d\'emonstration de la proposition, rappelons
bri\`evement le contenu du th\'eor\`eme~2 de~\cite{I}. On consid\`ere une
suite g\'eod\'esique $(x_0, \dots x_d)$, $d \geq 3$.
 Alors il existe $d-1$ composantes de Reeb (canoniques)  $(F_1,G_1,)
\dots , (F_{d-1},G_{d-1 })$ pour le couple $(x_0,x_d)$. En
particulier, chaque produit $F_i \times G_i$ est inclus dans
l'ensemble singulier de  $H$ ou de $H^{-1}$. De plus, sauf
\'eventuellement pour $i=1$, le bord $F_i$
s\'epare les points $x_j,j < i$ des points $x_j, j \geq i$, et 
$d_H(x_{i-1},F_i)=d_H(x_{i},F_i)=1$.\footnote{On utilise la notation
usuelle $d_H(x,E)$ pour le 
minimum de la distance entre $x$ et un point de l'ensemble $E$.}
Ceci est encore vrai pour $i=1$ sauf si la composante $(F_1,
G_1)$ est d\'eg\'en\'er\'ee, c'est-\`a-dire si $x_0$ appartient \`a $F_1$. On a des propri\'et\'es sym\'etriques pour les bords $G_i$.

\begin{demo}[de la proposition~\ref{p.geodesique-locale}]
%Soit $(x_i)_{i  \in I}$ une suite de points.
%D'apr\`es la remarque~\ref{r.geo-locale},
Il suffit de voir que pour tout $d \geq 3$, une suite g\'eod\'esique locale d'ordre $d$ est
encore une suite g\'eod\'esique locale d'ordre $d+1$. Ceci se voit sur les
sous-suites de longueur $d+2$~;
on consid\`ere donc une suite $(x_{0}, \dots , x_{d+1})$ g\'eod\'esique locale d'ordre $d$. On veut montrer que $d_H(x_0,x_{d+1}) \geq d+1$
(l'autre in\'egalit\'e est claire).

Puisque cette suite est g\'eod\'esique locale d'ordre $d$, on a
$d_H(x_0,x_d)=d$. 
On note  $(F_1,G_1,)
\dots , (F_{d-1},G_{d-1 })$ les $d-1$ composantes de Reeb pour le
couple $(x_0,x_{d})$ fournies par le th\'eor\`eme~2 de~\cite{I}.

Montrons, dans un premier temps,  que $d_H(x_{d+1}, G_{d-1}) >
1$.  
Pour cela, rappelons que $d_H(x_{d-2},F_{d-1}) =1$ (voir
ci-dessus). Si on avait \'egalement  $d_H(x_{d+1}, G_{d-1}) \leq 1$, 
  le lemme de ``d\'eviation des g\'eod\'esiques'' donnerait
$d_H(x_{d-2},x_{d+1}) \leq 2$, ce qui contredirait le fait que la suite 
$(x_i)$ est une g\'eod\'esique locale d'ordre $d$ avec $d \geq 3$.

Ceci montre en particulier que le point $x_d$ n'appartient pas au bord
$G_{d-1}$, et par cons\'equent $G_{d-1}$ s\'epare $x_d$ des points $x_i$,
$i < d$.
Notons $O_{x_{d}}(G_{d-1})$ la composante connexe de $\bbR^2 \setminus
G_{d-1}$ qui contient $x_d$. Puisque $d_H(x_d,x_{d+1})=1$, il existe un arc libre
reliant ces deux points~; et d'apr\`es ce qui pr\'ec\`ede, cet arc ne peut
pas rencontrer l'ensemble $G_{d-1}$. Par cons\'equent le point $x_{d+1}$
est dans $O_{x_{d}}(G_{d-1})$.

Consid\'erons enfin un arc g\'eod\'esique $\gamma$ allant de $x_0$ \`a
$x_{d+1}$~; nous souhaitons montrer que cet arc est de $H$-longueur au
moins $d+1$.
 Cet arc doit  rencontrer successivement tous les  bords 
des composantes de Reeb~:  $F_1,G_1, \dots , F_{d-1},G_{d-1}$.
 Consid\'erons une d\'ecomposition de $\gamma$. Pour chaque entier
$i=1, \dots, d-1$,  cette d\'ecomposition a  un sommet entre $F_i$ et
$G_i$ (voir par exemple~\cite{I}, lemme~4.8).
 S'il n'y avait aucun autre sommet, alors le dernier sous-arc de
la d\'ecomposition rencontrerait le bord $G_{d-1}$, ce qui contredirait 
l'in\'egalit\'e $d_H(x_{d+1}, G_{d-1}) > 1$.  Par cons\'equent la
d\'ecomposition a au moins $d$ sommets, et $\gamma$ a une $H$-longueur
au moins \'egale \`a $d+1$. Ceci prouve, comme attendu, que
$d_H(x_0,x_{d+1}) \geq d+1$. 
\end{demo}

%%%%%%%%%%%%%%%%%%%%%%%%
\subsection{Classe d'\'equivalence canonique de suites g\'eod\'esiques}
\label{ss.geodesique-canonique}
On se place maintenant sous les hypoth\`eses habituelles : $h$ est un hom\'eomorphisme du plan, pr\'eservant l'orientation,  fixant uniquement $0$, dont l'indice est diff\'erent de $1$. 
Il s'agit principalement de construire une suite g\'eod\'esique pour le relev\'e canonique $\wt h$, $\tau$-\'equivariante. Nous commen\c cons par le cas simple, o\`u l'indice est assez loin de $1$ ; puis nous examinerons les cas d'indice $0$ ou $2$ ; enfin, nous montrerons l'unicit\'e de la classe d'\'equivalence des suites g\'eod\'esiques \'equivariantes.

\begin{demo}[du premier point du corollaire~\ref{c.bouts-geodesique}, lorsque l'indice est diff\'erent de $0$ ou $2$]
Soit $\gamma$ une courbe de Jordan de degr\'e $1$ g\'eod\'esique pour $h$. Notons $d$ le module de $h$, et $(s_{i})_{i \in \bbZ/d\bbZ}$ les sommets d'une d\'ecomposition de $\gamma$, index\'es dans l'ordre cyclique correspondant \`a l'orientation positive de $\gamma$.
Soit $\Gamma : \bbR \ra \bbR^2$ le relev\'e de $\gamma$ au rev\^etement universel du compl\'ementaire du point fixe. Il existe  une suite $(\wt s_{i})_{i \in \bbZ}$ qui ``rel\`eve'' la suite finie $(s_{i})_{i \in \bbZ/d\bbZ}$, au sens o\`u
\begin{itemize}
\item $\pi(\wt s_{i})= s_{i \ \mathrm{mod} \ d}$ pour tout $i$ ;
\item les points $s_{i}, {i \in \bbZ}$ sont index\'es dans l'ordre correspondant \`a l'orientation de $\Gamma$ (on a notamment $\tau(\wt s_{i}) = \wt s_{i+d}$ pour tout $i$).
\end{itemize}
Il est clair que le relev\'e d'une courbe qui est libre pour $h$ est libre pour $\wt h$. On en d\'eduit que, pour tout entier $i$,
$$
d_{\wt h}(\wt s_{i}, \wt s_{i+1}) =1.
$$ 
En particulier, d'apr\`es l'in\'egalit\'e triangulaire, on a $d_{\wt h}(\wt s_{i+d}, \wt s_{i}) \leq d$.
D'autre part, puisque $\wt h$ est le relev\'e canonique de $h$, tout arc libre pour $\wt h$ se projette en un arc libre pour $h$ (proposition~\ref{p.releve-canonique}). On en d\'eduit facilement que, pour tout entier $i$, 
$$
d_{\wt h}(\wt s_{i+d}, \wt s_{i}) \geq d.
$$ 
Finalement, on a $d_{\wt h}(\wt s_{i+d}, \wt s_{i}) = d$ pour tout entier $i$.
D'apr\`es la remarque~\ref{r.geo-locale}, la suite $(\wt s_{i})_{i \in \bbZ}$ est une suite g\'eod\'esique locale d'ordre $d$ pour $\wt h$.

Lorsque l'indice du point fixe est diff\'erent de $0$, $1$, $2$, le module $d$ est minor\'e par $4$ 
(proposition~\ref{p.minoration-module}). On peut alors appliquer la proposition~\ref{p.geodesique-locale}, et la suite 
$(\wt s_{i})_{i \in \bbZ}$ est aussi une suite g\'eod\'esique ``globale''.
\end{demo}

\begin{demo}[du corollaire~\ref{c.bouts-geodesique}, lorsque l'indice vaut $0$ ou $2$]
On consid\`ere l'anneau $\bbR^2/\tau^2$, qui est un rev\^etement interm\'ediaire, \`a deux feuillets, de l'anneau $\bbR^2/\tau$ :
%\footnote{ DIAGRAMME xypic  ??}
$$
%\left(
\begin{array}{cccc}
 \bbR^2   \\
 \da \\
\bbR^2/\tau^2 & \cup \{0\}  &  \simeq &  \bbR^2  \\
 
 \da &  & & \da   \\
\bbR^2/\tau &  \cup \{0\}  &    \simeq &  \bbR^2 
\end{array}
%\right)
$$
 (dans ce diagramme, la derni\`ere fl\`eche est  un rev\^etement \`a deux feuillets ramifi\'e au-dessus du point $0$). 
L'hom\'eomorphisme  $\wt h$ induit un  hom\'eomorphisme $h_{2}$ sur l'anneau 
$\bbR^2/\tau^2$, qui s'etend en un hom\'eomorphisme du plan $\bbR^2$ qui fixe uniquement le point $0$, et que l'on note encore $h_{2}$. Cet hom\'eomorphisme est, bien s\^ur, l'un des deux relev\'es de l'hom\'eomorphisme $h$ par le rev\^etement \`a deux feuillets.
On suppose que l'indice du point fixe de $h$ est diff\'erent de $1$. Montrons alors que, pour  le nouvel hom\'eomorphisme $h_{2}$, le point fixe est d'indice diff\'erent de $0,1,2$ ; plus pr\'ecis\'ement, on va obtenir les relations suivantes.
\begin{affi}
\label{a.egalite-indices}
$$
\indi(h_{2},0) = \indi_{\tau^2}(\wt h) +1 = 2 \times \indi_{\tau}(\wt h) +1 = 2 \times (\indi(h,0) -1) +1.
$$
\end{affi}
La deuxi\`eme \'egalit\'e suit imm\'ediatement des d\'efinitions de l'indice de $\wt h $ par rapport aux translations $\tau$ et $\tau^2$. La derni\`ere \'egalit\'e caract\'erise le relev\'e canonique $\wt h$ de $h$ (proposition~\ref{p.releve-canonique}). Il nous reste \`a montrer la premi\`ere. 
%Rappelons qu'une \emph{droite de Brouwer d'extr\'emit\'es  $0$ et l'infini} pour $h$ est une droite topologique qui va de $0$ \`a l'infini\footnote{Autrement dit, son adh\'erence contient $0$, et elle est non born\'ee.}, et qui est libre pour $h$.
Commen\c cons par faire deux remarques.
\begin{enumerate}
\item Un arc libre pour $\wt h$ se projette par $\pi$ en un arc libre pour $h$ (propri\'et\'e du relev\'e canonique) ; par cons\'equent sa projection au rev\^etement interm\'ediaire est encore un arc libre pour $h_{2}$.
\item Il existe une droite de Brouwer
\footnote{\label{f.droite-de-brouwer} Rappelons la d\'efinition. On appelle \emph{droite topologique} une application continue injective $\Delta : \bbR \ra \bbR^2 \setminus\{ 0\}$ qui se prolonge contin\^ument en une application de la droite achev\'ee $\{-\infty \}\cup \bbR \cup \{+\infty \}$ dans la sph\`ere $\bbR^2 \cup\{\infty\}$ en envoyant $-\infty$ et $+\infty$ sur les points $0$ ou $\infty$. 
Si la droite topologique $\Delta$ va de $0$ \`a $\infty$, et si elle est disjointe de son image $h(\Delta)$, on dit que c'est une \emph{droite de Brouwer d'extr\'emit\'es $0$ et $\infty$}.}
 pour $h$ d'extr\'emit\'es  $0$ et $\infty$ : ceci suit par exemple du th\'eor\`eme principal de~\cite{0} (th\'eor\`eme~D, p31). Elle se rel\`eve clairement  en une droite de Brouwer pour $h_{2}$.
\end{enumerate}
La premi\`ere \'egalit\'e ci-dessus d\'ecoule alors  du lemme suivant (dont la preuve est rel\'egu\'ee plus bas).
Le lemme est appliqu\'e en rempla\c cant $h$ par $h_{2}$, $\tau$ par $\tau^2$, et donc $\pi$ par l'application quotient $\bbR \ra \bbR/\tau^2$. Les deux remarques pr\'ec\'edentes consistaient en fait \`a en v\'erifier les  hypoth\`eses dans ce cadre.
\begin{lemm}\label{l.lien-indices}
Soit $h$ un hom\'eomorphisme du plan, pr\'eservant l'orientation, fixant
uniquement le point $0$. Supposons que $h$ admette un relev\'e $\wt h$
tel que  l'image par $\pi$ de tout arc libre pour $\wt h$ est un arc libre
pour $h$. 
Supposons de plus que $h$ admette une droite de
Brouwer $\Delta$ d'extr\'emit\'es $0$ et l'infini.
Alors on a la relation
$$
\indi(h,0) = \indi_{\tau}(\wt h) +1.
$$
En particulier, si $\indi(h,0)$ est diff\'erent de $1$, alors $\wt h$ est le relev\'e canonique de $h$.
\end{lemm}
Nous pouvons maintenant achever la preuve du corollaire~\ref{c.bouts-geodesique}.
On reprend les notations de la preuve du corollaire dans le ``bon'' cas d'indice diff\'erent de $0,1,2$ : en relevant les sommets d'une courbe de Jordan $\gamma$ g\'eod\'esique pour $h$,  on a construit une
suite  $(\wt s_{i})_{i \in \bbZ}$, et il s'agit de voir qu'elle est une suite g\'eod\'esique globale pour $\wt h$.
Pour fixer les id\'ees, supposons que le point fixe de $h$ soit d'indice $2$.
Remarquons tout d'abord qu'il peut arriver que le module de $h$ soit sup\'erieur ou \'egal \`a $3$ ; dans ce cas, la suite est g\'eod\'esique locale d'ordre $3$, et la proposition~\ref{p.geodesique-locale} s'applique comme avant, qui nous dit que c'est aussi une suite g\'eod\'esique globale. On peut donc se restreindre au cas o\`u le module de $h$ est \'egal \`a $2$.
D'apr\`es l'affirmation~\ref{a.egalite-indices},
  le point fixe de $h_{2}$ est d'indice $3$. Le module de $h_{2}$ est donc minor\'e par $4$.
  D'autre part, comme $h$ a pour module $2$, et qu'un arc libre pour $h$ se rel\`eve clairement en un arc libre pour $h_{2}$, le module de $h_{2}$ est aussi major\'e par $4$. On a donc
  $$
  \diam(h_{2}) = 4 = 2 \times \diam(h).
  $$
Par cons\'equent, la courbe $\gamma$ se rel\`eve au rev\^etement \`a deux feuillets en une courbe de Jordan $\gamma_{2}$ g\'eod\'esique pour $h_{2}$. On peut maintenant appliquer le ``bon'' cas \`a la courbe $\gamma_{2}$ : la suite de ses sommets se rel\`eve en une suite g\'eod\'esique globale pour $\wt h$. Cette suite n'est autre que la suite 
$(\wt s_{i})_{i \in \bbZ}$ relev\'ee de la suite des sommets de $\gamma$.
Ceci termine la preuve.
%\footnote{Noter $\wt h $ relev\'e canonique de $h_{2}$ ??}
\end{demo}

\begin{demo}[du lemme]
La preuve est une variante d'un morceau de la construction du relev\'e canonique  de ~\cite{0} (proposition~4.17). 
Soit $\wt \Delta$ une droite topologique qui rel\`eve $\Delta$~; c'est une droite de Brouwer pour $\wt h$.
On commence par montrer qu'on a l'une des deux situations
$$
\tau^{-1} (\wt \Delta) < \wt h(\wt \Delta) < \wt \Delta \mbox{ ou } 
\wt \Delta < \wt h(\wt \Delta) < \tau(\wt \Delta).
$$
En effet, dans le cas contraire, puisque $h(\Delta)$ est disjointe de $\Delta$, on aurait
$  \wt h(\wt \Delta) < \tau^{-1} (\wt \Delta)  < \wt \Delta $ ou  
$\wt \Delta < \tau(\wt \Delta) < \wt h(\wt \Delta) $. On en d\'eduit que la bande d\'elimit\'ee par $\wt \Delta$ et $\tau(\wt \Delta)$ est libre par $\wt h$, et on trouverait ainsi un arc libre pour $\wt h$ qui se projette en un arc qui n'est pas libre pour $h$. Ceci contredit nos hypoth\`eses.

Pour fixer les id\'ees, on se place dans le cas o\`u $\wt \Delta < \wt h(\wt \Delta) < \tau(\wt \Delta)$. Tous les it\'er\'es positifs de $\wt \Delta$ sont alors inclus dans la bande d\'elimit\'ee par $\wt \Delta$ et $\tau(\wt \Delta)$~: dans le cas contraire, le domaine de translation engendr\'e par $\wt \Delta$ rencontrerait $\tau(\wt \Delta)$, et on trouverait un arc libre pour $\wt h$ reliant un point $\wt x$ sur $\wt \Delta$ au point $\tau(\wt x)$ (voir l'affirmation~3.4 de~\cite{I})~; la projection par $\pi$ de cet arc ne serait pas un arc, contrairement aux hypoh\`eses.

Ce qui   pr\'ec\`ede montre que la situation est la suivante~: 
$$
\wt \Delta < \wt h(\wt \Delta) < \wt h^2(\wt \Delta) <\wt  h^3(\wt \Delta) < \tau(\wt \Delta).
$$
%La th\'eorie de Brouwer nous dit que les it\'er\'es de $\Delta$ par $h$
% sont deux \`a deux disjoints, et l'ordre
%cyclique est donn\'e par
%$$
%\Delta < h(\Delta) < h^2(\Delta) < h^3(\Delta) < \Delta
%$$ 
%(voir par exemple~\cite{0}, remarque 3.24 et affirmation 3.25).
%Soit $\wt \Delta$ un relev\'e de $\Delta$.
% Montrons qu'on a 
%$$
%\wt \Delta < \wt h(\wt \Delta) < \wt h^2(\wt \Delta) < \wt h^3(\wt \Delta) <
%\tau(\wt \Delta).
%$$
%Pour cela, notons $\wt \Delta'$ le relev\'e de $h(\Delta)$
%situ\'e entre $\wt \Delta$ et $\tau(\wt \Delta)$. Par l'absurde, on suppose par exemple que
%$\wt h(\wt \Delta)\neq \wt \Delta'$.  
% On voit facilement que
%la bande situ\'ee entre $\wt \Delta$ et $\wt \Delta'$ est alors libre par $\t
%h$. Ceci permet de trouver un arc libre pour $\wt h$ qui se projette en
%un arc qui n'est pas libre pour $h$, ce qui contredit la premi\`ere hypoth\`ese du lemme.
%Les autres in\'egalit\'es s'obtiennent de mani\`ere analogue.
Quitte \`a conjuguer $h$, on peut  supposer que les ensembles $\Delta, \dots h^3(\Delta)$
sont des demi-droites euclidienne issues de $0$, dispos\'ees \`a angles droits,
et qu'on a
$$
\wt h^i(\wt \Delta)= \left\{\frac{i}{4}\right\} \times \bbR, \ \ \ i=0 \dots 3
$$
(et aussi $\tau(\wt \Delta) = \{1\} \times \bbR$).
%droits. On suppose aussi que le rev\^etement $\pi$ est donn\'e par la
%formule
%$$
%\begin{array}{rcl}
%\pi : \bbR^2 & \longrightarrow & \bbR^2 \setminus\{0\} \\
%(\wt \theta, r ) & \longmapsto & \exp(-r + 2i \pi \wt \theta)
%\end{array}
%$$
%l'automorphisme de rev\^etement \'etant alors $\tau(\wt \theta, r) =
%(\wt \theta +1, r)$. 
En  notant  $\wt \theta : \bbR^2 \ra \bbR$ la premi\`ere coordonn\'ee, 
on en d\'eduit facilement la propri\'et\'e $(\star)$~:
$$
  \mbox{ Pour tout } \wt x \in \bbR^2, \quad \quad |\wt \theta (\wt x)-\wt \theta
(\wt h(\wt x))| < \frac{1}{2}. 
$$
La relation d'indice en d\'ecoule  : c'est exactement l'objet du
lemme~4.41 de \cite{0}.
\end{demo}

%%%%%%%%%%%%%%%%%%%%%%%%
%%%%%%%%%%%%%%%%%%%%%%%%
\section{Droites de Brouwer}
\label{s.cRdB}
Dans cette section, nous montrons comment construire une droite de Brouwer \`a l'int\'erieur de chaque composante de Reeb. D'une part, ceci sera un outil important dans la preuve de la formule d'indice (section~\ref{s.lien-indice}). D'autre part, l'existence d'une certaine famille de  droites de Brouwer sera l'une des propri\'et\'es dynamiques lues sur l'indice par quarts-de-tour (section~\ref{s.croissants-petales}).

%Nous aurons besoin d'une version ``simpliciale'' de la proposition~\ref{p.cRdB}, qui tient compte de l'existence d'une d\'ecomposition en briques.

%%%%%%%%%%%%%%%
%\subsection{Composantes de Reeb et droites de Brouwer}
%\pourmoi{def droite de Brouwer en 2.5 ???}
Rappelons qu'une \emph{droite de Brouwer} pour un hom\'eomorphisme de Brouwer $H$ est une droite topologique $\Delta$, disjointe de son image par $H$, et qui s\'epare son image $H(\Delta)$ de son image r\'eciproque $H^{-1}(\Delta)$. 

\begin{prop}
\label{p.cRdB}
%[\ref{p.cRdB}-bis]
Soit $H$ un hom\'eomorphisme de Brouwer,
 et $\cal F$ une d\'ecomposition en
briques pour $H$. On consid\`ere une composante de Reeb $(F,G)$ pour $H$
(pour un couple de points $(x,y)$ quelconques).
On suppose que $F$ et $G$ sont invariants par
 $H$.
Alors il existe une droite de Brouwer pour $H$,
incluse dans $\cF$, qui s\'epare $F$ et $G$.
\end{prop}

On rappelle qu'une \emph{cha\^\i ne de briques} est une suite
$(B_0, \dots B_k)$ de briques de la d\'ecomposition $\cF$, telle que $H(B_{0})$ rencontre $B_{1}$, $H(B_{1})$ rencontre $B_{2}$, \dots, $H(B_{k-1})$ rencontre $B_{k}$. Le lemme de Franks pour les briques  dit que les briques d'une cha\^\i ne de briques sont toutes distinctes (\cf\ \cite{0}, lemme 3.65 p64).
\begin{lemm}\label{l.CDD}
Soit $(x,y)$ un couple de points singulier pour un hom\'eomorphisme de Brouwer $H$.
Alors il n'existe pas de cha\^\i ne de briques
$(B_0, \dots B_k)$ telle que $B_0$ contienne $y$ et $B_k$ contienne $x$.
\end{lemm}

\begin{demo}[du lemme~\ref{l.CDD}]
On suppose, par l'absurde, qu'il existe une cha\^\i ne de briques comme dans l'\'enonc\'e.
Dans un premier temps, pour simplifier, supposons que $x$ et $y$
appartiennent respectivement \`a l'\emph{int\'erieur} des briques $B_k$ et
$B_0$. Si le couple $(x,y)$ est singulier, par d\'efinition, 
il existe un it\'er\'e positif de $B_k$
qui rencontre $B_0$. On peut alors prolonger la suite $(B_0, \dots
B_k)$ en une cha\^{\i}ne  de briques de $B_0$ \`a $B_0$, ce qui est interdit
par le lemme de Franks rappel\'e ci-dessus.

Traitons le cas g\'en\'eral. Remarquons qu'un point est toujours inclus dans l'int\'erieur de la r\'eunion des briques qui le contiennent.
D'autre part le couple $(H(x), H^{-1}(y))$ est
aussi singulier. Il existe donc une brique $B_{k+1}$ qui contient $H(x)$,
et une brique $B_{-1}$ qui contient $H^{-1}(y)$, telles qu'un it\'er\'e
positif de $B_{k+1}$ rencontre $B_{-1}$. La suite $(B_{-1}, B_0,
\dots, B_k, B_{k+1})$ est encore une cha\^ine de briques, et on peut la 
prolonger  en une cha\^{\i}ne  de briques de $B_{-1}$ \`a $B_{-1}$,
ce qui est interdit.
\end{demo}

Rappelons encore quelques d\'efinitions du texte~\cite{0}. Un \emph{attracteur strict} est un ensemble $E$ tel que $H(E)$ est inclus  dans l'int\'erieur de $E$.
 \'Etant donn\'ee une brique $B$,
l'ensemble $A^+(B)$ est le plus petit ensemble qui contient $B$, qui est une r\'eunion de briques de $\cF$,  et qui est un attracteur strict. C'est un ensemble connexe.

\begin{demo}[de la proposition~\ref{p.cRdB}]
On se place dans le cas o\`u $F \times G$ est singulier pour $H$
(le cas o\`u $F \times G$ est singulier pour $H^{-1}$ est bien s\^ur sym\'etrique).
On consid\`ere l'ensemble
$$
A^+(G) := \bigcup\left\{ A^+(B), B \mbox{ brique rencontrant } G \right\}.
$$
On montre facilement que cet ensemble contient $G$ dans son int\'erieur, qu'il est connexe,
que c'est un attracteur strict. D'apr\`es le lemme~\ref{l.CDD}, il est
disjoint de $F$.
Soit alors $\Delta$ la composante connexe de la fronti\`ere de
$A^+(G)$ qui s\'epare $F$ et $G$. C'est une droite topologique,
elle est incluse dans le 1-squelette de la d\'ecomposition $\cF$.
 Comme elle est incluse dans la fronti\`ere
d'un attracteur strict, elle est libre. On v\'erifie facilement
que $\Delta$ s\'epare $H(\Delta)$ de $H^{-1}(\Delta)$
(en utilisant que $F$ et $G$ sont connexes et
invariants par $H$, et que $\Delta$ s\'epare $F$ et $G$). C'est donc une
droite de Brouwer.
\end{demo}

%%%%%%%%%%%%%%%%%%%%%%%%
%%%%%%%%%%%%%%%%%%%%%%%%
\section{Formule d'indice}
\label{s.lien-indice}
Dans cette section, nous prouvons le th\'eor\`eme~\ref{t.lien-indice}-bis,
qui relie l'indice du relev\'e canonique \`a l'indice symbolique. La preuve utilise la construction de droites de Brouwer dans les composantes de Reeb (section~\ref{s.cRdB}), 
et la notion d'indice partiel entre deux droites de Brouwer (\cite{0}, section 3.3).

On consid\`ere un hom\'eomorphisme $h$ du plan, v\'erifiant les hypoth\`eses habituelles
($h$ pr\'eserve l'orientation, le point $0$ est son unique point fixe, l'indice de Poincar\'e-Lefschetz est diff\'erent de $1$). On reprend les notations du paragraphe~\ref{ss.dynamique-equivariants}
($\wt h$ est le relev\'e canonique, $\Gamma$ une droite g\'eod\'esique relev\'ee d'une courbe de Jordan g\'eod\'esique $\gamma$,  $(\wt x_{k})_{k \in \bbZ}$ la suite des sommets, $(\wt F_{k},\wt G_{k})_{k \in \bbZ}$ la suite de composantes de Reeb, et 
$\wt {\textbf{M}}([\tau])$ le mot p\'eriodique). L'ordre des bords des composantes de Reeb est donn\'e par $\cdots \leq \wt F_{k} < \wt G_{k} \leq \wt F_{k+1} < \wt G_{k+1} \leq \cdots$ (voir l'appendice~\ref{as.ordre-cyclique}).

%$\wt h$ est le relev\'e canonique de $h$, $\Gamma$ est une droite g\'eod\'esique pour $\wt h$ qui est invariante par $\tau$, $(\wt x_{k})_{k \in \bbZ}$ est la suite des sommets d'une d\'ecomposition de $\Gamma$ qui est $\tau$-\'equivariante, de p\'eriode $d$. Enfin, 
% on note  $(\wt F_{k},\wt G_{k})_{k \in \bbZ}$ la suite de composantes de Reeb associ\'ee
%\`a la suite $(\wt x_{k})_{k \in \bbZ}$, et 
%On note $\wt {\textbf{M}}([\tau])$ le mot infini, sur l'alphabet $\{ \ra, \la, \da, \ua\}$, associ\'e
%\`a $[\tau]$ (th\'eor\`eme~I.3 et  d\'efinition~I.8.2).
%Comme $\wt h$ commute avec la translation $\tau$, ce mot infini est p\'eriodique de p\'eriode $2d$. On note $\textbf{M}([\tau])$ le mot cyclique \`a $2d$ lettres induit par ce mot p\'eriodique.

\begin{demo}[du th\'eor\`eme~\ref{t.lien-indice}-bis]
On se donne une d\'ecomposition en briques g\'eod\'esiques $\cal F$
pour $h$ (lemme~\ref{l.decomposition-geodesique}),
qui est compatible avec la courbe de Jordan g\'eod\'esique $\gamma$.
L'ensemble $\wt {\cal F}= \pi^{-1}({\cal F})$ est une d\'ecomposition en briques pour $\wt h$ (proposition~\ref{p.releve-canonique}). L'intersection de cette d\'ecomposition  avec la droite $\Gamma$ est une suite de sommets d'une d\'ecomposition de $\Gamma$, et on pourra  supposer que
$$
\wt {\cal F} \cap \Gamma = \{\wt x_{k} , k \in \bbZ  \}.
$$

On applique maintenant la proposition~\ref{p.cRdB} : pour chaque entier 
$k$, il existe une droite de Brouwer $\wt \Delta_{k}$, incluse dans $\mathcal{\wt F}$, qui s\'epare les deux bords de la composante de Reeb $(\wt F_{k},\wt G_{k})$. On a donc n\'ecessairement $\wt \Delta_{k} \cap \Gamma = \{\wt x_{k}\}$~: notamment, l'arc $\Gamma_{k}= [\wt x_{k} \wt x_{k+1}]_{\Gamma}$ est un arc libre qui va de $\wt \Delta_{k}$ \`a $\wt \Delta_{k+1}$.
D'autre part, puisque la d\'ecomposition 
$\mathcal{\wt F}$ est invariante par $\tau$, et que la famille 
des composantes de Reeb est \'equivariante,
on peut supposer que $\wt \Delta_{k+d}=\tau(\wt \Delta_{k})$.
%\footnote{C'est d'ailleurs le cas si on d\'efinit les droites de Brouwer par la formule donn\'ee dans la preuve de la proposition~\ref{p.cRdB}.}
%D'apr\`es le lemme 7.10 de~\cite{I},
L'ordre des composantes de Reeb et des droites de Brouwer est alors donn\'e par~:
$$
\cdots \leq \wt F_{k} < \wt \Delta_{k} < \wt G_{k} \leq \wt F_{k+1} < \wt \Delta_{k+1} < \wt G_{k+1} \leq \cdots
$$
%Une partie de ces informations est repr\'esent\'ee sur la figure~\ref{f.dynamique-equivariant}.
\begin{figure}[htbp]
\begin{center}
\def\JPicScale{0.8}
\ifx\JPicScale\undefined\def\JPicScale{1}\fi
\psset{unit=\JPicScale mm}
\psset{linewidth=0.3,dotsep=1,hatchwidth=0.3,hatchsep=1.5,shadowsize=1}
\psset{dotsize=0.7 2.5,dotscale=1 1,fillcolor=black}
\psset{arrowsize=1 2,arrowlength=1,arrowinset=0.25,tbarsize=0.7 5,bracketlength=0.15,rbracketlength=0.15}
\begin{pspicture}(0,0)(135,76)
\psline(10,40)(50,40)
\psdots[](10,40)
(50,40)
\psline(50,40)(90,40)
\psdots[](50,40)
(90,40)
\psline(90,40)(130,40)
\psdots[](90,40)
(130,40)
\psline[linestyle=dashed,dash=1 1](10,40)(5,40)
\psline[linestyle=dashed,dash=1 1]{->}(130,40)(135,40)
\psline(50,69)(50,10)
\psline(90,70)(90,10)
\pscustom[linewidth=0.5]{\psbezier(35,70)(30,60)(30,60)(32,60)
\psbezier(32,60)(34,60)(33,61)(35,61)
\psbezier(35,61)(37,61)(29,48)(31,48)
\psbezier(31,48)(33,48)(32,50)(34,50)
\psbezier(34,50)(36,50)(30,38)(32,38)
\psbezier(32,38)(34,38)(33,41)(35,41)
\psbezier(35,41)(37,41)(31,25)(33,25)
\psbezier(33,25)(35,25)(34,31)(36,31)
\psbezier(36,31)(38,31)(30,12)(32,12)
\psbezier(32,12)(34,12)(33,14)(35,14)
\psbezier(35,14)(37,14)(35,10)(35,10)
}
\pscustom[linewidth=0.5]{\psbezier(105,10)(108.5,20)(108.5,20)(107.1,20)
\psbezier(107.1,20)(105.7,20)(106.4,19)(105,19)
\psbezier(105,19)(103.6,19)(109.2,32)(107.8,32)
\psbezier(107.8,32)(106.4,32)(107.1,30)(105.7,30)
\psbezier(105.7,30)(104.3,30)(108.5,42)(107.1,42)
\psbezier(107.1,42)(105.7,42)(106.4,39)(105,39)
\psbezier(105,39)(103.6,39)(107.8,55)(106.4,55)
\psbezier(106.4,55)(105,55)(105.7,49)(104.3,49)
\psbezier(104.3,49)(102.9,49)(108.5,68)(107.1,68)
\psbezier(107.1,68)(105.7,68)(106.4,66)(105,66)
\psbezier(105,66)(103.6,66)(105,70)(105,70)
}
\pscustom[linewidth=0.5]{\psbezier(76,70)(71,60)(71,60)(73,60)
\psbezier(73,60)(75,60)(74,61)(76,61)
\psbezier(76,61)(78,61)(70,48)(72,48)
\psbezier(72,48)(74,48)(73,50)(75,50)
\psbezier(75,50)(77,50)(68,33)(70,33)
\psbezier(70,33)(72,33)(74,41)(76,41)
\psbezier(76,41)(78,41)(72,25)(74,25)
\psbezier(74,25)(76,25)(75,31)(77,31)
\psbezier(77,31)(79,31)(71,12)(73,12)
\psbezier(73,12)(75,12)(74,14)(76,14)
\psbezier(76,14)(78,14)(76,10)(76,10)
}
\pscustom[linewidth=0.5]{\psbezier(65.1,10)(68.6,20)(68.6,20)(67.2,20)
\psbezier(67.2,20)(65.8,20)(66.5,19)(65.1,19)
\psbezier(65.1,19)(63.7,19)(71.4,36)(70,36)
\psbezier(70,36)(68.6,36)(67.2,30)(65.8,30)
\psbezier(65.8,30)(64.4,30)(68.6,42)(67.2,42)
\psbezier(67.2,42)(65.8,42)(66.5,39)(65.1,39)
\psbezier(65.1,39)(63.7,39)(67.9,55)(66.5,55)
\psbezier(66.5,55)(65.1,55)(65.8,49)(64.4,49)
\psbezier(64.4,49)(63,49)(68.6,68)(67.2,68)
\psbezier(67.2,68)(65.8,68)(66.5,66)(65.1,66)
\psbezier(65.1,66)(63.7,66)(65.1,70)(65.1,70)
}
\rput(35,6){$\t{\cal F}_0$}
\rput(76,6){$\t{\cal F}_1$}
\rput(65,6){$\t{\cal G}_0$}
\rput(45,36){$\t x_0$}
\rput(128,36){$\t x_2$}
\rput(86,36){$\t x_1$}
\rput(25,42.5){$\Gamma_{-1}$}
\rput(61.25,42.5){$\Gamma_{0}$}
\rput(100,42.5){$\Gamma_{1}$}
\rput(134,43){$\Gamma$}
\rput(50,75){$\t \Delta_0$}
\rput(90,73){$\t \Delta_1$}
\rput(105,6){$\t{\cal G}_1$}
\pscustom[]{\psbezier(4,40)(1,40)(2,36)(2,33)
\psbezier(2,33)(2,30)(3,27)(2,26)
\psbezier(2,26)(1,25)(1,25)(1,25)
\psbezier(1,25)(1,25)(1,25)(2,24)
\psbezier(2,24)(3,23)(2,20)(2,17)
\psbezier(2,17)(2,14)(2,10)(2,10)
}
\rput[r](-1,25){$O_S(\Gamma)$}
\pscustom[]{\psbezier(50,70)(50,73)(46,72)(43,72)
\psbezier(43,72)(40,72)(37,71)(36,72)
\psbezier(36,72)(35,73)(35,73)(35,73)
\psbezier(35,73)(35,73)(35,73)(34,72)
\psbezier(34,72)(33,71)(30,72)(27,72)
\psbezier(27,72)(24,72)(20,72)(20,72)
}
\psline[linestyle=dashed,dash=1 1](2,10)(2,4)
\psline[linestyle=dashed,dash=1 1](20,72)(14,72)
\rput(35,76){$O_g(\t \Delta_0)$}
\end{pspicture}
\caption{Preuve de la formule d'indice pour $\wt h$}
\label{f.dynamique-equivariant}
\end{center}
\end{figure}
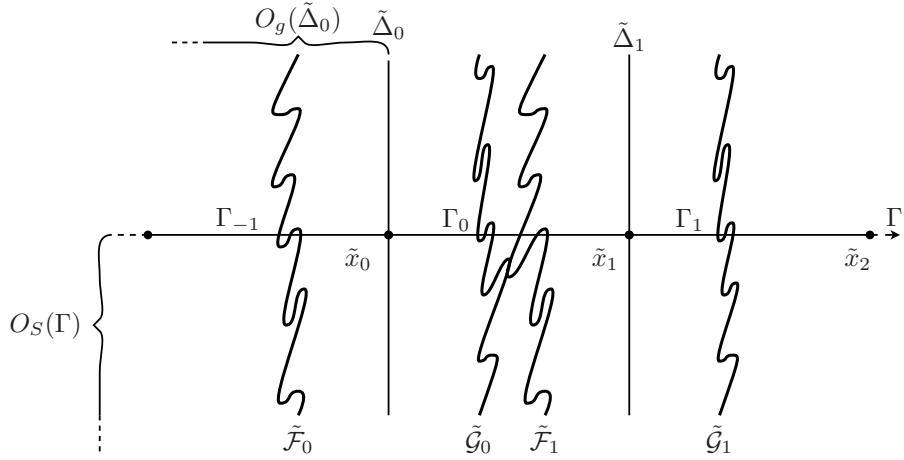

On oriente chaque droite $\wt \Delta_{k}$ du bas vers le haut (voir la figure~\ref{f.dynamique-equivariant}), et on note $O_{g}(\wt \Delta_{k})$ le demi-plan topologique ouvert situ\'e \`a gauche de $\wt \Delta_{k}$. La droite $\Gamma$ est orient\'ee par l'action de $\tau$ (de la gauche vers la droite sur la figure) ; on note $O_{S}(\Gamma)$ le demi-plan topologique ouvert situ\'e \`a droite de $\Gamma$. On consid\`ere l'indexation $(m_{i})_{i \in \bbZ}$ du mot $\wt {\textbf{M}}([\tau])$
compatible avec les notations qui pr\'ec\'edent, autrement dit~:
\begin{itemize}
\item $m_{2k}= \ra$ si et seulement si $\wt h(\Gamma_{k-1}) \cap \Gamma_{k} \neq \emptyset$ (d\'efinition 3.11 de~\cite{I}) ;
\item $m_{2k+1}= \da$ si et seulement si $\wt G_{k}$ et $\wt F_{k+1}$ sont de type dynamique Nord-Sud (d\'efinition 8.1 de~\cite{I}).
\end{itemize}
Voici deux autres liens entre le mot $\wt {\textbf{M}}([\tau])$ et la dynamique de $\wt h$ :
%Pour relier l'indice de $\wt h$ \`a l'indice symbolique, nous aurons besoin des liens suivants entre le mot $\wt {\textbf{M}}([\tau])$ et la dynamique de $\wt h$ :
% \footnote{Verifier indices $k$ ??}
\begin{enumerate}
\item la fl\`eche $m_{2k}$ vaut $\ra$ si et seulement si $\wt h(\wt \Delta_{k})$ est situ\'ee \`a droite de $\wt \Delta_{k}$ ;
\item supposons que  $m_{2k}= \ra$ et $m_{2k+2}= \la$ ; alors $m_{2k+1}=\da$ si et seulement si 
$\wt h(\Gamma_{k}) \subset O_{S}(\Gamma)$. 
\end{enumerate}
Le deuxi\`eme point est l'affirmation 8.5 de~\cite{I}.
Expliquons le premier point, par contrappos\'ee. Si  $\wt h(\wt \Delta_{k})$ est  \`a gauche de $\wt \Delta_{k}$, alors le demi-plan $O_{g}(\wt \Delta_{k})$ est attractif, \emph{i. e.} on a
$$
\wt h(O_{g}(\wt \Delta_{k})) \subset O_{g}(\wt \Delta_{k}).
$$
Dans ce cas, la composante $(\wt F_{k},\wt G_{k})_{k \in \bbZ}$ ne peut pas \^etre singuli\`ere pour $\wt h $. Par d\'efinition des composantes de Reeb, elle est donc singuli\`ere pour $\wt h^{-1}$, et la fl\`eche $m_{2k}$ vaut $\la$ (\cite{I}, th\'eor\`eme 1-bis, fin de la section 4).

Nous pouvons maintenant commencer \`a calculer l'indice de $\wt h $ par rapport \`a la translation $\tau$.
$$
\begin{array}{rcl}
 \indi(\wt h, \tau) & = & \inpa (\wt \Delta_0, \wt \Delta_{d}) \\
		& = & \displaystyle \sum_{k=0}^{d-1} \inpa (\wt \Delta_k, \wt \Delta_{k+1}).
\end{array}
$$
La premi\`ere \'egalit\'e est une cons\'equence imm\'ediate des d\'efinitions de
l'indice de $\wt h$ par rapport \`a $\tau$ et de l'indice partiel.
La deuxi\`eme vient de la relation d'additivit\'e (\cite{0}, lemme 3.38, p49).
D'autre part, par d\'efinition, l'indice partiel du mot cyclique est donn\'e par
$$
\begin{array}{rcl}
\inpa_{c}((m_{k})_{k \in \bbZ}) & = &  \displaystyle \sum_{i=0}^{2d-1}\inpa(m_im_{i+1}) \\
 & = &\displaystyle \sum_{k=0}^{d-1}\inpa(m_{2k}m_{2k+1}m_{2k+2}).
\end{array}
$$
La formule d\'ecoulera alors imm\'ediatement du lemme-cl\'e ci-dessous.
\end{demo}
\begin{lemm}\label{l.inpa}
Pour tout entier $k$,
$$
\inpa (\wt \Delta_k, \wt \Delta_{k+1})= \inpa(m_{2k}m_{2k+1}m_{2k+2}).
$$
\end{lemm}
\begin{demo}[du lemme]
Il s'agit d'utiliser les liens entre le mot et la dynamique.
\'Etudions d'abord le cas o\`u $m_{2k}= m_{2k+2}$. Alors l'indice partiel symbolique $\inpa(m_{2k}m_{2k+1}m_{2k+2})$ est nul (pratiquement par d\'efinition).
D'autre part, d'apr\`es le point $1$ ci-dessus, 
$\wt h$ pousse les deux droites $\wt \Delta_k$ et $\wt \Delta_{k+1}$ dans le m\^eme sens ;
autrement dit, avec le vocabulaire du texte~\cite{0}, le couple $(\wt \Delta_k, \wt \Delta_{k+1})$ est indiff\'erent. De plus, l'arc $\Gamma_{k}$ est libre et joint les deux droites. Par cons\'equent, l'indice partiel entre les deux droites est nul (lemme 3.39 de~\cite{0}, p50).
La conclusion du lemme est donc v\'erifi\'ee dans le cas indiff\'erent.

Il reste \`a traiter les cas attractif et r\'epulsif. 
Pour fixer les id\'ees, supposons que le mot est $(\ra \da \la)$, les
autres cas se traitant par un raisonnement similaire. L'indice partiel
du mot symbolique vaut alors $-1/2$. D'apr\`es le premier point ci-dessus,
le couple $(\wt \Delta_k, \wt \Delta_{k+1})$ est attractif ;
d'apr\`es le deuxi\`eme point, l'image de l'arc libre $\Gamma_{k}$ est situ\'ee dans le demi-plan $O_{S}(\Gamma)$, en dessous-de $\Gamma$. On peut encore appliquer le lemme~3.39 de~\cite{0}, qui affirme que l'indice partiel vaut $-1/2$, comme attendu.
\end{demo}

%%

%%%%%%%%%%%%%%%%%%%%%%%%
%%%%%%%%%%%%%%%%%%%%%%%%
\section{Compatibilit\'e des composantes de Reeb}
 \label{s.composantes-de-reeb}
Dans cette section, on montre que le rev\^etement $\pi$ \'etablit une correspondance entre les composantes de Reeb pour $\wt h$ associ\'ees \`a la g\'eod\'esique $[\tau]$ et les composantes de Reeb minimales de $h$.

%%%%%%%%%%%%%%%%
\subsection{Un lemme de topologie}
On consid\`ere une courbe de Jordan $\gamma: \bbS^1 \ra \bbR^2 \setminus \{0\}$ de degr\'e $1$, et son relev\'e
$\Gamma : \bbR \ra \bbR^2$. Rappellons quelques notations et d\'efinitions de~\cite{I}.
La droite topologique $\Gamma$ \'etant orient\'ee par son param\'etrage, on note $O_N(\Gamma)$ et $O_S(\Gamma)$ les deux  composantes connexes
du compl\'ementaire de $\Gamma$, situ\'ees respectivement \`a gauche et \`a
droite de $\Gamma$. Notons que $\pi(O_N(\Gamma))$ et $\pi(O_S(\Gamma))$ sont respectivement des voisinages point\'es des points $0$ et $\infty$ dans la sph\`ere $\bbR^2 \cup \{\infty\}$.
 On se donne une partie ferm\'ee   $\wt C$ du plan, qui
rencontre $\Gamma$ selon un ensemble compact.
On d\'efinit une compactification
``Nord-Sud'' de $\wt C$ en munissant l'ensemble  
 $\wt C \cup\{N,S\}$ d'une topologie compatible avec celle de $\wt C$, 
telle que des bases respectives de voisinages de $N$ et $S$ sont donn\'ees par
$
\{O_N(\Gamma) \setminus K, K \mbox{ compact du plan} \}
$ et
$
\{O_S(\Gamma) \setminus K, K \mbox{ compact du plan} \}.
$
Comme dans~\cite{I}, on dira que \emph{$\wt C$ s\'epare les
deux bouts de $\Gamma$} si $\wt C$ s\'epare les deux composantes
connexes non born\'ees de $\Gamma \setminus \wt C$.
Dans le cas o\`u $\wt C$ est connexe, ceci revient \`a dire que l'ensemble $\wt C \cup \{N,S\}$ est encore connexe, ou encore que $N$ et $S$ sont dans l'adh\'erence de $\wt C$.
%$N$ est le point \`a l'infini du compactifi\'e d'Alexandroff de $\wt F \cap
%O_N(\Gamma)$,
% et  $S$ est le point \`a l'infini du compactifi\'e d'Alexandroff de $\wt F \cap
%O_S(\Gamma)$. VERIF??
On \'etend  le rev\^etement $\pi$ en posant 
$
\hat \pi(N)= 0 \ \ , \hat \pi (S)=\infty.
$
On d\'efinit les ensembles 
$
C := \pi(\wt C) \mbox{ et } \hat C := C \cup \{0, \infty\} = \hat \pi (\wt C \cup\{N,S\}).
$
\begin{lemm}[compatibilit\'e topologique]
\label{l.compatibilite-topologique}
Soit $\wt C$ un ferm\'e du plan, connexe, qui s\'epare les bouts de $\Gamma$.
On suppose de plus que \begin{enumerate}
\item les ensembles de la famille $\{\tau^q(\wt C), q
\in \bbZ\}$ sont deux \`a deux disjoints,
\item  cette famille est
localement finie (aucun compact du plan ne rencontre une infinit\'e
d'\'el\'ements de cette famille).
\end{enumerate}
On a alors~:
\begin{enumerate}
\item L'ensemble $C$ est ferm\'e (dans $\bbR^2 \setminus \{0\}$) et connexe~;
%\item l'ensemble $\hat C$ est connexe~;
\item l'application $\hat \pi$ induit un hom\'eomorphisme  
de $\wt C \cup\{N,S\}$ sur $\hat C$.
\end{enumerate}
\end{lemm}
On appliquera ce lemme aux bords des composantes de Reeb de $\wt h$, et
aussi aux composantes ``remplies'', voir plus loin.

\begin{demo}[du lemme~\ref{l.compatibilite-topologique}]
La d\'emonstration est laiss\'ee au lecteur.
\end{demo}

%%%%%%%%%%%%%%%%
\subsection{Compatibilit\'e des composantes de Reeb}
Soit $h$ un hom\'eomorphisme du plan, pr\'eservant l'orientation, fixant uniquement $0$, d'indice diff\'erent de $1$. On consid\`ere  le relev\'e $\wt h$ fourni par la proposition~\ref{p.releve-canonique}.
%\begin{prop}~
%\label{p.compatibilite-composantes}
%\begin{enumerate}
%\item Si un couple $(\wt F, \wt G)$ est une composante de Reeb minimale pour le relev\'e canonique $\wt h$,
%associ\'ee \`a $[\tau]$, alors l'application $\pi$ induit un hom\'eomorphisme de $\wt F \cup \wt G \cup \{N,S\}$ sur son image, et  
%le couple $(F,G)=(\pi(\wt F),\pi(\wt G))$ est une
%composante de Reeb minimale pour $h$. 
%\item R\'eciproquement,  si $(F,G)$ est une composante de Reeb minimale pour $h$,
%alors il existe un couple $(\wt F, \wt G)$ qui  est une composante de
%Reeb minimale pour $\wt h$,
%associ\'ee \`a $[\tau]$, et tel que $\pi(\wt F)=F$ et $\pi(\wt G)=G$.
%\end{enumerate}
%\end{prop}
%Dans la situation d\'ecrite par la proposition, 
% on dira que $(F,G)$ est la \emph{projet\'ee} de la composante $(\wt F,\t
% G)$, et que  $(\wt F,\wt  G)$ est une composante \emph{relev\'ee} de
% $(F,G)$.

\begin{demo}[de la proposition~\ref{p.compatibilite-composantes}]
On reprend les notations de la section~\ref{ss.dynamique-equivariants}.
Soit $(\wt F, \wt G)$  une composante de Reeb minimale pour $\wt h$ associ\'ee \`a $[\tau]$.
%Soit $(\wt F_{k}, \wt G_{k})$  une composante de Reeb minimale pour $\wt h$ associ\'ee \`a $[\tau]$.
%Soit $\gamma$ une courbe de Jordan  g\'eod\'esique pour $h$, $\Gamma$ la droite topolgique relev\'ee de $\gamma$. Soit $(\wt x_k)_{k \in \bbZ}$ une suite g\'eod\'esique
%  obtenue en relevant une d\'ecomposition de
%$\gamma$ (corollaire~\ref{c.releve-geodesique}). On consid\`ere une composante de Reeb $(\t
%F,  \wt G)$ associ\'e \`a la courbe $\Gamma$ (\cite{I}, d\'efinition
%7.3)~: autrement dit,  $(\wt F,\wt G)$ est l'unique
%composante de Reeb minimale pour un couple $(\wt x_{i-1},\t
%x_{i+1})$. Pour fixer les id\'ees, on supposera que $i=0$.
%Par unicit\'e des composantes de Reeb, le couple $(\wt F, \wt G)$ co\"i ncide avec l'une des composantes de la suite ${\cal R}([\tau])= (\wt F_{k},\wt G_{k})_{k \in \bbZ}$~; 
Quitte \`a d\'ecaler les indices, on supposera que $(\wt F, \wt G)= (\wt F_{0}, \wt G_{0})$.
Rappellons que  $O(\wt F, \wt G)$ d\'esigne l'unique composante connexe de $\bbR^2
\setminus (\wt F \cup \wt G)$ dont l'adh\'erence rencontre \`a la fois $\wt F$
et $\wt G$. On notera $\wt C$ l'adh\'erence de $O(\wt F, \wt G)$ ; $\wt C$ est en fait \'egal \`a $\wt F \cup O(\wt F, \wt G) \cup \wt G$
\footnote{En effet, la $\wt h$-boule $B_{1}(\wt x_{1})$ est un ouvert connexe disjoint de $\wt F$, et son adh\'erence contient $\wt F$ par construction des composantes de Reeb : voir~\cite{I}, section~4.3.}.
%Une propri\'et\'e de topologie g\'en\'erale indique que 
%$\wt C \subset \wt F \cup O(\wt F, \wt G) \cup  \wt G$.
%Par minimalit\'e de la composante de
%Reeb, on a en fait $\wt C=\wt F \cup O(\wt F, \wt G) \cup  \wt G$.
%~; on appelle
%\emph{composante de Reeb remplie} l'ensemble $\wt C$).
On veut appliquer le lemme~\ref{l.compatibilite-topologique} \`a l'ensemble $\wt C$.
Il est clair que $\wt C$ est ferm\'e, connexe, et s\'epare les bouts de $\Gamma$. Montrons que les \'el\'ements de la famille $\{\tau^q(\wt C), q \in \bbZ\}$
sont deux \`a deux disjoints, et que la famille est localement finie.
 Comme $(\wt F,\wt G)$ est la
composante de Reeb minimale associ\'ee au couple $(\wt x_{-1},\t
x_{1})$, par unicit\'e des composantes de Reeb minimales,
   $(\tau^q(\wt F),\tau^q(\wt G))$ est
la composante de Reeb minimale associ\'ee au couple $(\tau^q(\wt x_{-1}),\tau^q(\t
x_{+1})) = (\wt x_{qd-1},\wt x_{qd+1})$, o\`u $d$ est le module de $h$.
Puisque $d \geq 2$, on a donc
$$
\wt F < O(\wt F, \wt G) < \wt G < \tau^q(\wt F) < \tau^q(O(\wt F, \wt G)) < \tau^q(\wt G).
$$
%(voir \cite{I}, lemme 7.10 et figure 17).
On en d\'eduit que les ensembles $\wt C$ et $\tau^q(\wt C)$ sont disjoints.   Supposons maintenant  que la famille n'est pas localement finie.
Il existe alors un point $\wt z$ dont tout voisinage $V$
  rencontre une infinit\'e d'\'el\'ements de la famille. Choisissons un
voisinage $V$ connexe et libre,  et notons $q_1$ et
$q_2$ deux entiers tels que $\tau^{q_1}(\wt C)$ et 
$\tau^{q_2}(\wt C)$ rencontrent $V$, avec $q_1 < q_2$.
 Comme le module  $d$ est plus grand que $2$,
il existe une composante de reeb $(\wt F_k , \wt G_k)$ situ\'ee entre $\tau^{q_1}(\wt C)$ et $\tau^{q_2}(\wt C)$~: 
$$
\tau^{q_1}(\wt C) < \wt F_k < \wt G_k < \tau^{q_2}(\wt C).
$$
Puisque $V$ est connexe et rencontre les deux ensembles extr\^emes, il
doit rencontrer aussi $\wt F_k$ et $\wt G_k$. Comme $V$ est libre, on a
trouv\'e un point $\wt x$ de $\wt F_k$ et un point $\wt y$ de $\wt G_k$ tels
que $d_{\wt h}(\wt x, \wt y) =1$. Mais puisque $(\wt F_k, \wt G_k)$ est une
composante de Reeb, le couple $(\wt x, \wt y)$ est singulier pour $\wt h$
ou $\wt h^{-1}$, et on devrait avoir $d_{\wt h}(\wt x, \wt y) =2$ (\cite{I},
lemme 4.2), ce qui est absurde. La famille est donc localement finie, et le lemme~\ref{l.compatibilite-topologique} s'applique : on en d\'eduit que l'application $\hat \pi$ est un hom\'eomorphisme de $\wt C \cup \{N,S\}$ sur son image, et que les ensembles $F = \pi(\wt F)$ et $G=\pi(\wt G)$ sont ferm\'es dans $\bbR^2 \setminus\{0\}$, connexes, et $\hat F$ et $\hat G$ sont encore connexes. Ceci prouve que $F$ et $G$ ont la topologie requise pour former une composante de Reeb (deux premiers points
de la d\'efinition~\ref{d.reeb}).

Montrons alors qu'il existe $\varepsilon=\pm 1$ tel que tout couple $(\wt x,\wt y) \in \wt  F \times \wt G$ est singulier \emph{via} l'ouvert $O(\wt F,\wt G)$ pour l'hom\'eomorphisme $h^\varepsilon$. Rappelons que les ensembles $\wt F$, $O(\wt F, \wt G)$ et $\wt G$ sont
invariants par $\wt h$ (\cite{I}, lemme 7.10, point 3). 
Par d\'efinition des composantes de Reeb,
tout couple $(\wt x,\wt y)$ dans $\wt F \times \wt G$ est singulier
pour $\wt h$ ou $\wt h^{-1}$, et il s'agit de voir que  $(\wt x,\wt y)$
est singulier \textit{via} l'ouvert $O(\wt F, \wt G)$.
Soit $O_{\wt G}(\wt F)$ la composante connexe du compl\'ementaire de
$\wt F$ qui contient $\wt G$. On note ${\cal F}$ le compl\'ementaire de 
 $O_{\wt G}(\wt F)$. On d\'efinit sym\'etriquement l'ensemble $\cal G$. 
Les ensembles  ${\cal F}$ et  $\cal G$ sont ferm\'es et contiennent
respectivement les ensembles $\wt F$ et $\wt G$, ils sont invariants par $\t
h$, et on a une  partition du plan
$$
\bbR^2 = {\cal F} \cup O(\wt F, \wt G) \cup {\cal G}.
$$ 
 Soient $V_{\wt x}$  un voisinage de $\wt x$
disjoint de $\cal G$, et  $V_{\wt y}$ un voisinage de $\wt y$
disjoint de $\cal F$. Pour tout  entier $n$,   $\wt h^n(V_{\t
x}) \cap V_{\wt y} \subset O(\wt F, \wt G)$. Ceci prouve que le couple 
 singulier $(\wt x,\wt y)$ est singulier \textit{via} l'ouvert $O(\wt F, \wt G)$. 

Les ensembles $\wt F$ et $\wt G$ sont invariants par $\wt h$, par cons\'equent $F$ et $G$ sont invariants par $h$.
 L'ensemble $\pi (O(\wt F, \wt G))$ co\"incide avec le domaine
 $O(F,G)$ (par d\'efinition, voir l'appendice~\ref{as.ordre-cyclique}). 
Il est alors clair que $F
 \times G$ est inclus dans  l'ensemble des couples singuliers pour $h^{\pm 1}$ \emph{via} l'ouvert
 $O(F,G)$.
  Le couple $(F,G)$ est donc une composante de
Reeb. On montre facilement que la minimalit\'e de la composante de Reeb
$(\wt F,\wt G)$ entra\^ine celle de $(F, G)$.

%%%%%%%
Nous avons montr\'e la premi\`ere implication de la proposition. Montrons
la deuxi\`eme. Soit $(F,G)$ une composante de Reeb pour $h$.
On a l'ordre cyclique
$$
F < O(F,G) < G < O(G,F) < F.
$$
Par d\'efinition de l'ordre cyclique, on peut choisir une
composante connexe $\wt F$ de $\pi^{-1}(F)$ et
 une composante connexe $\wt G$ de $\pi^{-1}(G)$
 de fa\c{c}on \`a ce que 
 $$
\wt F < O(\wt F, \wt G) < \wt G < O(\wt G, \tau(\wt F)) < \tau(\wt F),
$$
et on a $\pi (O(\wt F, \wt G)) = O(F, G)$.
%(comme avant, on a not\'e $O(\wt F, \wt G)$ l'unique composante connexe de $\bbR^2
%\setminus (\wt F \cup \wt G)$ dont l'adh\'erence rencontre $\wt F$ et $\wt G$).
L'ensemble  $\wt C:= \wt F \cup  O(\wt F, \wt G) \cup \wt G$ est alors disjoint de
ses images par les applications $\tau^q$. D'autre part, puisque son image par $\pi$ est ferm\'ee, il est ferm\'e, et la famille de ses it\'er\'es par $\tau$ est localement finie : d'apr\`es le lemme~\ref{l.compatibilite-topologique}, 
l'application $\pi$ induit un hom\'eomorphisme de 
$\wt C$ sur son image.

D'autre part, 
les ensembles $F$ et $G$ sont invariants par
$h$ (par hypoth\`ese)~; comme $h$ pr\'eserve l'orientation, l'ensemble
$O(F,G)$ est aussi invariant par $h$. Par cons\'equent, l'ensemble 
$\wt C$ est invariant ou libre 
pour $\wt h$. Montrons qu'il n'est pas libre~: dans le cas contraire,
on consid\`ere un arc $\gamma$, dans l'ouvert connexe invariant $O(F,G)$, reliant
un point \`a son image par $h$~; le relev\'e de $\gamma$ dans 
 $O(\wt F, \wt G)$ serait un arc libre, bien que $\gamma$ ne soit pas
libre, ce qui contredit la propri\'et\'e du relev\'e canonique (proposition~\ref{p.releve-canonique}). On en d\'eduit que $\wt F$ et $\wt G$ sont invariants par $\wt h$ ;  
l'application $\pi$ conjugue les restrictions respectives de $\wt h$ et $h$. Notamment, 
le produit $\wt F \times \wt G$ est singulier pour $\wt h$ ou pour $\wt h^{-1}$.
On consid\`ere alors la suite g\'eod\'esique  $(\wt x_k)_{k \in \bbZ}$. Pour tout
$k$ assez grand, l'ensemble $\wt F$  s\'epare le point
$\wt x_{-k}$ de l'ensemble $\wt G \cup \{\wt x_k\}$~; de m\^eme,
l'ensemble $\wt G$  s\'epare le point 
$\wt x_{k}$ de l'ensemble $\wt F \cup \{\wt x_{-k}\}$
 Ceci montre que le couple $(\wt F, \wt G)$ est une
composante de Reeb pour $(\wt x_{-k}, \wt x_{k})$. 
On montre facilement que la minimalit\'e de la composante de Reeb
$(F,G)$ entra\^ine celle de $(\wt F, \wt G)$.
D'apr\`es le lemme~7.8 de
\cite{I},  $(\wt F, \wt G)$ est bien une composante de Reeb associ\'ee \`a
la suite g\'eod\'esique $(\wt x_i)_{i \in \bbZ}$, ce qui termine la preuve
de la deuxi\`eme implication de la proposition.
\end{demo}

\part{Cons\'equences dynamiques}
\label{p.consequences-dynamiques}
Dans toute cette partie, on consid\`ere un hom\'eomorphisme $h$ du plan, v\'erifiant les hypoth\`eses habituelles ($h$ pr\'eserve l'orientation, fixe uniquement $0$, avec un indice diff\'erent de $1$). 
Dans la section~\ref{s.ouverts-dynamiques}, on construit des ouverts de points ayant la m\^eme dynamique, indiqu\'ee par les fl\`eches du mot $\mathbf{M}(h)$ (th\'eor\`eme~\ref{t.ouverts-dynamiques}). En corollaire, on pr\'ecise un th\'eor\`eme de Pelikan et Slaminka (\cite{pesla})~: l'indice de Poincar\'e-Lefschetz d'un point fixe isol\'e pour un hom\'eomorphisme de surface qui pr\'eserve une bonne mesure est inf\'erieur ou \'egal \`a $1$, et son indice par quarts-de-tour est trivial.
Dans la section~\ref{s.croissants-petales}, on montre comment lire sur le mot $\mathbf{M}(h)$ l'existence de droites de Brouwer organis\'ees en croissants et p\'etales, ce qui pr\'ecise le th\'eor\`eme D de~\cite{0} (p31).
Dans la section~\ref{s.autres}, on montre une propri\'et\'e dynamique des secteurs indiff\'erents, et on propose une conjecture sur l'existence d'ensembles stables et instables canoniques.

%%%%%%%%%%%%%%%%%%%%%%%%
%%%%%%%%%%%%%%%%%%%%%%%%
\section{Ouverts dynamiquement coh\'erents}
\label{s.ouverts-dynamiques}
Dans cette partie, on associe canoniquement \`a $h$ des ouverts, en nombre \'egal au module de $h$, constitu\'es de points dont le comportement dynamique est indiqu\'e par l'indice par quarts-de-tour. On en d\'eduira  que l'indice par quarts-de-tour d'un hom\'eomorphisme conservatif est trivial.

%%%%%%%%%%%%%%%%%%%%%%%
\subsection{D\'efinition}
On considère le relevé canonique $\wt h $ de $h$.  
%Comme dans la section~\ref{s.lemmes}, on note $\gamma$ une courbe de Jordan géodésique pour $h$, $\Gamma$ la droite géodésique pour $\wt h$ qui relève $\gamma$,
%et $(\wt x_k)_{k \in \bbZ}$ la suite géodésique des sommets d'une décomposition minimale de $\Gamma$ qui est $\tau$-équivariante, de période $d$, où $d$ est  le module de $h$.
Dans ce qui pr\'ec\`ede, on a associ\'e (canoniquement) \`a $h$ une classe d'\'equivalence $[\tau]$ de suites g\'eod\'esiques \'equivariantes. La suite de composantes de Reeb associée à $[\tau]$  est désignée par $((\wt F_k, \wt G_k))_{k \in \bbZ}$.
Pour chaque entier $k$, on d\'efinit un ensemble $\wt \delta _{k}$ de la mani\`ere suivante.
%Choisissons une droite g\'eod\'esique $\Gamma$ dans la classe $[\tau]$ : autrement dit,  $((\wt F_k, \wt G_k))_{k \in \bbZ}$ est la suite de composantes de Reeb associ\'ee \`a $\Gamma$.
 Choisissons une courbe de Jordan g\'eod\'esique $\gamma$, et $\Gamma$ la droite topologique qui rel\`eve $\gamma$. 
Donnons-nous, comme avant, une d\'ecomposition minimale de $\Gamma$  qui rel\`eve une d\'ecomposition de $\gamma$ ; cette d\'ecomposition admet exactement un sommet 
$\wt x_{k}$
dans l'ouvert  $O(\wt F_k, \wt G_k)$ entre $\wt F_k$ et $ \wt G_k$.
On d\'efinit  $\wt \delta _{k}$ comme l'ensemble des points $\wt x_{k}$ obtenus de cette mani\`ere (pour tous les choix possibles d'une courbe g\'eod\'esique $\gamma$ munie d'une d\'ecomposition).
\begin{prop}
L'ensemble $\wt \delta_k$ est un ouvert non vide,  invariant par $\wt h$, inclus dans l'ouvert  $O(\wt F_k, \wt G_k)$ situé entre $\wt F_k$ et $\wt G_k$.
La suite $(\wt \delta_k)_{k \in \bbZ}$ est équivariante, elle est canoniquement associée à $h$.
\end{prop}

\begin{demo}
Le fait d'\^etre libre est une propri\'et\'e stable : tout arc assez proche d'un arc libre est libre. Le fait que $\delta_{k}$ soit ouvert provient essentiellement de cette remarque.
Le reste de la proposition est imm\'ediat. 
%Les d\'etails sont laiss\'es au lecteur.
\end{demo}

\begin{defi}
L'ensemble $\delta_{k} := \pi(\wt \delta_k)$ est appel\'e \emph{ouvert des sommets associ\'e \`a la composante de Reeb $( F_k, G_k)$}.
\end{defi}

\begin{theo}
\label{t.ouverts-dynamiques}
La suite d'ouverts
$$
(\delta_k)_{k \in \bbZ/d\bbZ}
$$
est canoniquement associée à $h$ ; chacun de ces ouverts est constitué de points ayant la m\^eme dynamique,
qui est indiquée par l'indice par quarts-de-tour : plus précisément, pour tout point $x$ dans $\delta_k$, on a
$$
\alpha(x) = \alpha(m_{2k-1} \ m_{2k} \ m_{2k+1})  \mbox{ et } \omega(x) = \omega(m_{2k-1} \ m_{2k} \ m_{2k+1}).
$$
\end{theo}

\begin{demo}
Ici encore, tout est \'evident : notamment, la dynamique des sommets
est d\'ej\`a donn\'ee par le th\'eor\`eme~\ref{t.indice-courbes}. 
\end{demo}

%%%%%%%%%%%%%%%%%%%%%%%
\subsection{Indice des hom\'eomorphismes conservatifs}
%%%%%%%%%%%%
L'\'enonc\'e suivant pr\'ecise un r\'esultat de Pelikan et Slaminka
(\cite{pesla}, voir aussi Le Calvez~\cite{leca}, Simon~\cite{simon}, Niki\v sin~\cite{niki}).
\begin{coro}
\label{c.conservatif}
Soit $h_0$ un homéomorphisme d'une surface orientable $S$, préservant l'orientation, et $x_0$ un point fixe isolé de $h_0$.
Supposons que $h_0$ préserve une mesure $\mu_{0}$ qui charge les ouverts, et qui est finie sur un voisinage de $x_0$. Si de plus l'indice du point fixe est diff\'erent de $1$, 
alors
\begin{enumerate}
\item le module  de $h_0$ en $x_0$ est pair ;
\item l'indice par quarts-de-tour de $h_0$ en $x_{0}$ est périodique de période $4$,
$$
\mathbf{M}(h_0,x_0) = (\ua \ra \da \la \cdots ) \ ;
$$
\item  l'indice de Poincaré-Lefschetz du point $x_0$ est égal à $1-\frac{\diam(h)}{2}$.
\end{enumerate}
En particulier, l'indice est toujours inf\'erieur ou \'egal \`a $1$.
\end{coro}

\begin{demo}
On commence par  se ramener au cadre de l'article en appliquant les propositions~\ref{p.cadre1} et~\ref{p.cadre2}, qui nous fournissent un hom\'eomorphisme $h$ du plan, fixant uniquement le point $0$, et  qui est localement conjugu\'e \`a $h_{0}$. De plus, la mesure $\mu_{0}$ se rel\`eve en une mesure $\mu$, invariante par $h$, qui charge les ouverts et est finie sur un voisinage $V$ du point fixe $0$.\footnote{Soit $p$ l'application fournie par la proposition~\ref{p.cadre1}~; alors pour tout disque topologique $D$ sur lequel $p$ est injective, on a $\mu(D) = \mu_{0}(p(D))$. Cette propri\'et\'e caract\'erise la mesure $\mu$.}
Consid\'erons les ensembles stables 
$$
W^s_{loc}(0) : = \{ x \in \bbR^2 \mid  \forall n \geq0 , h^n(x) \in V  \}
$$
et
$$
W^s(0) := \{ x \in \bbR^2 \mid \omega(x) =\{0\} \    \}.
$$
Puisque tout point poss\`ede un voisinage errant, et que la mesure de $V$ est finie,
la mesure de $W^s_{loc}(0)$ est nulle.
Comme l'ensemble $\omega$-limite de tout point est r\'eduit \`a $0$ ou $\infty$, on a 
$$
W^s(0) = \bigcup_{p \geq 0} h^{-p} ( W^s_{loc}(0) ).
$$
On en d\'eduit que l'ensemble $W^s(0)$ est aussi de mesure nulle.

Montrons alors par l'absurde que l'indice $\textbf{M}(h)$ ne contient pas la s\'equence \mbox{$(\ua \la)$}.
Dans le cas contraire, le th\'eor\`eme~\ref{t.ouverts-dynamiques} nous fournit un ouvert $\delta $ inclus dans l'ensemble $W^s(0)$. Comme la mesure $\mu$ charge les ouverts, 
on a une contradiction.

On montre de la m\^eme mani\`ere que $\textbf{M}(h)$ ne contient aucune des s\'equences  $(\ra \ua)$, $(\da \ra)$, $(\la \da)$. On en d\'eduit que le mot
$\textbf{M}(h)$ est  obtenu en concat\'enant un certain nombre de copies de la s\'equence $(\ua \ra \da \la)$. En particulier, il est
p\'eriodique de p\'eriode $4$. Comme sa longueur est le double du module, le module de $h$ est pair. Le th\'eor\`eme~\ref{t.lien-indice} donne alors
$$
\indi(h,0) = \indi(\textbf{M}(h)) = 1 + \frac{\diam(h)}{2} \times \inpa_{c}((\ua \ra \da \la))
= 1-\frac{\diam(h)}{2}.
$$
\end{demo}

%%%%%%%%%%%%%%%%%%%%%%%%
%\subsection{Une conjecture}

%Terminons cette partie par une conjecture. Cas conservatif ou g\'en\'eral ?

%-- def branche stable : cc invariante de $W_{\infty \ra 0}$.

%-- clair : chaque bord rencontre une branche. Probl\`eme = 1) disjointes ? 2) rien d'autre ?

%-- il existe exactement $d$ branches, stables ou instables, 
%chacune des branche rencontre les bords de deux composantes de Reeb adjacentes, et son caract\`ere stable ou instable est donn\'e par la fl\`eche verticale correspondante dans le mot $\textbf{M}(h)$.

%

%%%%%%%%%%

%%%%%%%%%%%%%%%%%%%%%%%
\section{Croissants et p\'etales}
\label{s.croissants-petales}
Dans cette section, nous expliquons comment on peut lire sur l'indice par quarts-de-tour l'existence d'un certain nombre de droites de Brouwer. Ceci pr\'ecise le th\'eor\`eme principal du texte~\cite{0} (th\'eor\`eme D, p31).

%%%%%%%%%%%%%%%%%%%%%%%
\subsection{Croissants}
\label{ss.croissants}
%\footnote{dtes de B $\`a-\infty$ deja utilisees section 7.2, lemme 7.4 ???}
%On appelle \emph{droite topologique} une application continue injective $\Delta : \bbR \ra \bbR^2 \setminus\{ 0\}$ qui se prolonge contin\^ument en une application de la droite achev\'ee $\{-\infty \}\cup \bbR \cup \{+\infty \}$ dans la sph\`ere $\bbR^2 \cup\{\infty\}$ en envoyant $-\infty$ et $+\infty$ sur les points $0$ ou $\infty$. 
%Si la droite topologique $\Delta$ va de $0$ \`a $\infty$, et si elle est disjointe de son image $h(\Delta)$, on dit que c'est une \emph{droite de Brouwer d'extr\'emit\'es $0$ et $\infty$}. 
%La d\'efinition des droites de Brouwer d'extr\'emit\'es  $0$ et $\infty$ a \'et\'e rappel\'ee dans une note de bas de page de la section~\ref{s.cRdB}.
Rappelons qu'une \emph{droite de Brouwer} pour $h$ est une droite topologique qui va de $0$ \`a $\infty$ et qui est disjointe de son image par $h$ (voir la note~\ref{f.droite-de-brouwer} \`a la section~\ref{ss.geodesique-canonique}).
La th\'eorie de Brouwer permet notamment de voir qu'une telle droite  est disjointe de tous ses it\'er\'es (voir~\cite{0}, affirmation 3.25, p42). 
%Le r\'esultat qui suit est contenu dans la preuve du th\'eor\`eme~\ref{t.lien-indice}, \`a la section~\ref{s.lien-indice}.
%%
\begin{prop}
\label{p.croissants}
On reprend les notations du th\'eor\`eme~\ref{t.compo-Reeb}.
 Pour chaque entier $k$ entre $1$ et $d$, il existe une droite de
Brouwer $\Delta_k$, d'extr\'emit\'es $0$ et $\infty$, contenue dans l'ouvert
$O(F_k,G_k)$. De plus,  la fl\`eche $m_{2k}$ est \'egale \`a $\ra$
 si et seulement si
 on a l'ordre cyclique 
 $$h^{-1}(\Delta_{k}) < \Delta_{k} < h(\Delta_{k}) < h^{-1}(\Delta_{k}).$$
\end{prop}

\begin{demo}
Il suffit de poser $\Delta_{k} = \pi(\wt \Delta_{k})$, o\`u $\wt \Delta_{k}$ est la droite de Brouwer pour $\wt h$ obtenue \`a la section~\ref{s.lien-indice}.
Comme $\wt \Delta_{k}$ est libre, et que $\pi$ projette les arcs libres pour $\wt h$ en des arcs, la courbe $\Delta_{k}$ est sans point double. Comme elle est sans point double et incluse dans une d\'ecomposition en briques $\cF$, elle est proprement plong\'ee : c'est donc une droite topologique. Puisqu'elle est  transverse \`a la courbe de Jordan $\gamma = \pi(\Gamma)$, elle va de $0$ \`a $\infty$.
Comme c'est la projet\'ee d'une droite libre, et que $\pi$ projette les arcs libres pour $\wt h$ en des arcs libres,
elle est encore libre~: c'est donc une droite de Brouwer.
Il reste \`a comprendre le lien avec l'indice par quarts-de-tour : 
or nous avons d\'ej\`a vu que $m_{2k} = \ra $ si et seulement si $\wt h (\wt \Delta_{k})$ est situ\'ee \`a droite de $\wt \Delta_{k}$ (voir la section~\ref{s.lien-indice}, preuve du th\'eor\`eme~\ref{t.lien-indice}). Ceci implique la fin de la proposition.
\end{demo}

%%%%%%%%%%%%%%%%%%%%%%%%%%%%
\subsection{P\'etales}
\label{ss.petales}
\begin{figure}[htbp]
\begin{center}
\input{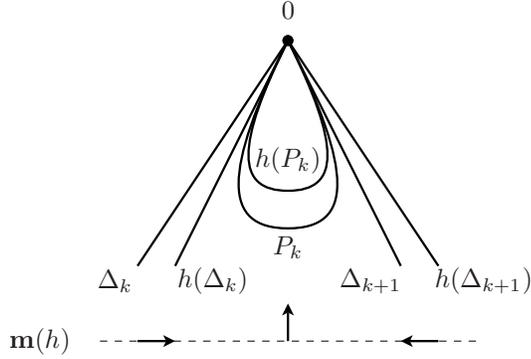}
\caption{P\'etales et croissants}
\label{f.petales}
\end{center}
\end{figure}

Rappelons qu'un \emph{p\'etale attractif  bas\'e en $0$} est un disque topologique ferm\'e $P$ du plan, qui contient le point $0$ dans son bord, et qui est un attracteur strict, au sens o\`u
$h(P) \subset \inte(P) \cup\{0\}$. On d\'efinit de mani\`ere analogue les p\'etales r\'epulsifs bas\'es en $0$, et les p\'etales attractifs ou r\'epulsifs  bas\'es en l'infini. Par d\'efinition, le bord d'un p\'etale bas\'e en $0$ (resp. l'infini) est une \emph{droite de Brouwer d'extr\'emit\'es $0$} (resp. $\infty$). La proposition suivante est illustr\'ee sur la figure~\ref{f.petales}.
%\footnote{Un mot de trois lettres est dit \emph{symboliquement hyperbolique}
%(respectivement \emph{elliptique, indiff\'erent}) si son indice partiel
%symbolique est $-1/2$ (respectivement $1/2$, $0$).
%Une composante de Reeb $(F_k,G_k)$ est  dite \emph{hyperbolique,
%elliptique} ou \emph{indiff\'erente} selon  le type symbolique du
%sous-mot correspondant $(m_{2k-1}\ m_{2k} \ m_{2k+1})$. ???}
\begin{prop}
\label{p.petales}
Soit $k$ un entier modulo $d$. Supposons que le sous-mot 
$(m_{2k}m_{2k+1}m_{2k+2})$ du mot cyclique $\mathbf{M}(h)$ soit \'egal \`a $(\ra \ua \la)$.  Alors il existe
un p\'etale attractif $P_{k}$ bas\'e en $0$. De plus, on peut choisir $P_{k}$ et les droites de Brouwer de la proposition~\ref{p.croissants} pour avoir
$$
P_{k} \subset O(\Delta_{k}, \Delta_{k+1}) \cup \{0\}.
$$ 
%o\`u $O(\Delta_{k}, \Delta_{k+1})$ d\'esigne l'ouvert situ\'e entre les droites de Brouwer $\Delta_{k}$ et $ \Delta_{k+1}$ donn\'ees par la proposition pr\'ec\'edente.
On peut \'egalement choisir le p\'etale $P_{k}$ inclus dans un voisinage donn\'e de $0$.
\end{prop}

Bien s\^ur, on a un r\'esultat analogue pour les autres types de p\'etales, qui sont induits par les s\'equences $(\la \da \ra)$, $(\ra \da \la)$ et $(\la \ua \ra)$. En corollaire, on trouve une nouvelle d\'emonstration de la version topologique du  th\'eor\`eme de la fleur de Leau-Fatou (th\'eor\`eme A de ~\cite{0}) : si l'indice de Poincar\'e-Lefschetz du point fixe est $p+1 > 1$, la formule reliant l'indice \`a l'indice par  quarts-de-tour permet de voir qu'on a trouv\'e au moins $2p$ croissants alternativement  attractifs et r\'epulsifs autour du point fixe (ceci n\'ecessite un raisonnement combinatoire tout \`a fait analogue \`a la preuve du th\'eor\`eme F du texte~\cite{0} (p81).

\begin{demo}
On se donne une d\'ecomposition en briques $\cF$ g\'eod\'esique pour $h$ (lemme~\ref{l.decomposition-geodesique}). On peut choisir les droites $\Delta_{k}, \Delta_{k+1}$ incluses dans $\cF$. Elles d\'elimitent alors ce qu'on a  appel\'e un \emph{croissant attractif simplicial} dans le texte~\cite{0}.
% On applique alors la proposition~4.8 de ce texte :
% ce croissant contient un croissant attractif simplicial minimal, et celui-ci contient un p\'etale simplicial $P_{k}$.
Comme la d\'ecomposition est g\'eod\'esique, il existe une brique $B_{k}$ qui traverse le croissant (\ie\ elle rencontre $\Delta_{k}$ et $\Delta_{k+1}$) ; de plus, la fl\`eche $m_{2k+1} = \ua$ indique que $h(B_{k})$ est au-dessus de $B_{k}$ dans le croissant.
Comme \`a la section~\ref{s.cRdB}, on consid\`ere le plus petit ensemble $A^+(B_{k})$ qui contient $B_{k}$, qui est r\'eunion de briques et qui est un attracteur strict.
Il existe une unique composante connexe $\Delta'_{k}$ du bord de $A^+(B_{k})$ qui rencontre $B_{k}$. On montre facilement que $\Delta'_{k}$ est une droite de Brouwer d'extr\'emit\'es $0$, et qu'elle d\'elimite un p\'etale attractif $P_{k}$ bas\'e en $0$ et inclus dans  le croissant ferm\'e d\'elimit\'e par $\Delta_{k}$ et $\Delta_{k+1}$.
Le lemme 5.1 de~\cite{0} (p104) dit que tout p\'etale attractif simplicial  contient des p\'etales attractifs arbitrairement petits, ce qui explique la derni\`ere phrase de la proposition.
Enfin, si $P_{k}$ est un p\'etale attractif bas\'e en $0$ et inclus dans le croissant ferm\'e d\'elimit\'e par $\Delta_{k}$ et $\Delta_{k+1}$, alors 
l'ensemble $h(P_{k})$ est encore un p\'etale, et est inclus dans $O(\Delta_{k}, \Delta_{k+1}) \cup \{0\}$.
 \end{demo}

%%%%%%%%%%%%%%%%%%%%%%%%
%%%%%%%%%%%%%%%%%%%%%%%%
\section{Autres propri\'et\'es}
\label{s.autres}
%%%%%%%%%%%%%%%%%%%%%%%
\subsection{Orbites \`a rebrousse-poils des secteurs indiff\'erents}

\begin{prop}
Supposons que le mot $\textbf{M}(h)$ contienne deux fl\`eches verticales successives du m\^eme type,
$$
(m_{2k-1} m_{2k} m_{2k+1})  = (\da \  ? \da)
$$
o\`u le point d'interrogation d\'esigne n'importe laquelle des deux fl\`eches $\ra$ ou $\la$.
Alors l'ouvert $O(F_{k}, G_{k})$ contient  un point $z$ dont la dynamique est inverse de celle indiqu\'ee par les fl\`eches :
$$\alpha(z) = \infty \neq \alpha(\da \ ? \da) , \ \ \omega(z) = 0 \neq \omega(\da \ ? \da) .$$
\end{prop}
Bien s\^ur, on a un r\'esultat analogue pour les mots $(\ua \  ? \ua)$.

\begin{demo}
Nous allons utiliser le lemme 9.5 de l'article~\cite{I}, dont voici l'\'enonc\'e adapt\'e \`a notre contexte. 
\begin{lemm}
Supposons que pour tout voisinage $V_{\infty}$ de l'infini dans la sph\`ere $\bbR^2 \cup\{\infty \}$, et pour tout voisinage $V_{0}$ de $0$, il existe un entier positif $n$
et un point $x_{1}$ de $O(F_{k}, G_{k})$ tel que $x_{1} \in V_{\infty}$ et $h^n(x_{1}) \in V_{0}$.
Alors l'adh\'erence de $O(F_{k}, G_{k})$ contient un point de  $W_{\infty \ra 0}$. 
\end{lemm}

V\'erifions tout d'abord  les hypoth\`eses du lemme. On se donne un voisinage $V_{\infty}$ de l'infini et un  voisinage $V_{0}$ de $0$.
Par d\'efinition des fl\`eches verticales, les bords $F_{k}$ et $G_{k}$ sont tous les deux de type dynamique $0 \ra \infty$. En particulier, on peut trouver un point $x$ dans 
$F_{k}$ et un point $y$ dans $G_{k}$, tels que $x$ et $y$ sont tous les deux dans l'ensemble $W_{0 \ra \infty}$. Quitte \`a remplacer $x$ par l'un de ses it\'er\'es positifs, on peut supposer que $x \in V_{\infty}$ ; de m\^eme, quitte \`a remplacer $y$ par l'un de ses it\'er\'es n\'egatifs, on peut supposer que $y \in V_{0}$.
Par d\'efinition des composantes de Reeb, le couple $(x,y)$ est singulier \emph{via} l'ouvert $O(F_{k}, G_{k})$ : en particulier, il existe un point $x_{1}$  v\'erifiant les hypoth\`eses du lemme.

On en conclut que l'adh\'erence de $O(F_{k}, G_{k})$ contient un point $z$ de  $W_{\infty \ra 0}$. D'apr\`es la proposition~\ref{p.dynamique-bord} et la d\'efinition des fl\`eches verticales, les bords $F_{k}, G_{k}$ ne contiennent aucun point allant de $\infty$ \`a $0$ ; par cons\'equent, le point $z$ appartient \`a
$$
\adhe(O(F_{k}, G_{k})) \setminus \left(  F_{k} \cup  G_{k} \right)  = O(F_{k}, G_{k}),
$$
ce que l'on voulait.
%Notons que cette derni\`ere \'egalit\'e est due \`a la minimalit\'e des ensembles $F_{k}$ et $G_{k}$.
\end{demo}

%%%%%%%%%%%%%%%%%%%%%%%%
%%%%%%%%%%%%%%%%%%%%%%%%
\subsection{Une conjecture concernant les branches stables et instables}
%Pour $\alpha, \omega \in \{0,\infty \}$, on note
%$$
%W_{\alpha \ra \omega} : = \{ x \in \bbR^2 \setminus \{0\} \mid \alpha(x) = \alpha \mbox{ et } \omega(x) = \omega   \}.
%$$
\begin{defi}
On appelle \emph{branche h\'et\'erocline} une partie $b$ de $\bbR^2 \setminus \{0\}$ qui est une composante connexe de l'ensemble $W_{0 \ra \infty}$ ou de l'ensemble $W_{\infty \ra 0}$, et qui est invariante par $h$.
\end{defi}

\begin{conj*}
Notons $d$ le module de $h$. Alors il existe exactement $d$ branches h\'et\'eroclines
$$
b_{1}, \dots , b_{d}.
$$
On peut choisir l'indexation pour que la branche $b_{k}$ rencontre les bords $G_{k}$ et $F_{k+1}$ des composantes de Reeb. La dynamique dans $b_{k}$ est alors donn\'ee par le mot $\textbf{M}(h)$ : 
$$
b_{k} \subset W_{0 \ra \infty} \mbox{ si et seulement si } m_{2k+1}=\da,
$$
$$
b_{k} \subset W_{\infty \ra 0} \mbox{ si et seulement si } m_{2k+1}=\ua.
$$
\end{conj*}
Voici ce que l'on sait dans la direction de cette conjecture.
Chaque bord de composante de Reeb rencontre une unique branche h\'et\'erocline (proposition~\ref{p.dynamique-bord}). Notons $b(F_{k})$ la branche qui rencontre $F_{k}$. Supposons que le mot $\textbf{M}(h)$ contienne une sous-suite de fl\`eches alternativement $\ua$ et $\da$ :
$$
(m_{2k_{0}+1}, m_{2k_{1}+1}, \dots , m_{2k_{p}+1} ) = (\ua, \da, \dots , \da)
$$
avec $1 \leq k_{0} < \dots  < k_{p} \leq d$. Alors les branches correspondantes
$$
b(F_{2k_{0}+1}), b(F_{2k_{1}+1}), \dots , b(F_{2k_{p}+1} )
$$
sont deux \`a deux disjointes. En effet, deux branches de types dynamiques distincts sont clairement disjointes ; en utilisant la proposition~\ref{p.dynamique-bord}, on voit que l'ensemble $F_{2k_{0}+1} \cup F_{2k_{2}+1}$ s\'epare les ensembles $b(F_{2k_{1}+1})$ et $b(F_{2k_{3}+1})$ (si $p \geq 4$).
D'autre part, en utilisant la formule d'indice, on montre facilement que le mot $\textbf{M}(h)$ contient une telle sous-suite, avec $p = 2\mid \indi(h,0) -1 \mid$.
On en d\'eduit qu'il existe au moins $2\mid \indi(h,0) -1 \mid$ branches h\'et\'eroclines.
Il resterait \`a montrer : 
\begin{enumerate}
\item  $b(G_{k}) = b(F_{k+1})$ ;
\item $b(F_{k}) \neq b(G_{k})$ ;
\item toute branche h\'et\'erocline rencontre un bord de composante de Reeb.
\end{enumerate}
Le point 3 impliquerait qu'il n'existe qu'un nombre fini de branches h\'et\'eroclines, ce qui n'est  pas connu.
Le plus difficile est probablement le point 2 (dans le cas o\`u les deux branches sont du m\^eme type dynamique). Notons n\'eanmoins que, dans le cas conservatif, ce point est trivial, puisque les types dynamiques des deux bords d'une composante de Reeb sont oppos\'es (corollaire~\ref{c.conservatif}).

%
%\footnote{Grouper par type de sous-mot.
%(+) Invariant trivial en conservatif (et version locale ?) Oui (?) cf $(\ua \la)$ ?
%\begin{enumerate}
%\item Global
%\begin{enumerate}
%\item (+) Dtes de Brouwer dans les cR ;
%\item (+) $(\ua \la)$ : il existe un ouvert positivement invariant de points qui vont \`a $0$.
%\item 
%\item  (+) $\ra \ua \la$ (mots elliptiques) : p\'etales dans tout secteur attractif entre deux telles droites de Brouwer.
%Arbitrairement grand ou petit...
%\item (-) $\da \ra \ua $ (mots elliptiques bis) : ouvert de points allant de $0$ \`a $0$.
%\end{enumerate}

%\item Local

%-- Branches stables : si on restreint un bord $F \cap \wio$ \`a un voisinage $U$, o obtient un ensemble $k$, il existe $n$ tel que $h^n(k) \subset U$ ? cf en pr\'eservant l'aire, branches stables canoniques ????

%-- p\'etales locaux.

%\end{enumerate}
%}

%\pagebreak	 
%%%%%%%%%%%%%%%%%%%%%%%%%%%
%%%%%%%%%%%%%%%%%%%%%%%%%%% 
 \appendix
 \part{Appendices}

%%%%%%%%%%%%%%%%%%%%%%%%
%%%%%%%%%%%%%%%%%%%%%%%%
\section{Points fixes d'indice $1$}
\label{as.indice-1}
Dans cette section, nous consid\'erons un hom\'eomorphisme du plan, pr\'eservant l'orientation, fixant uniquement $0$. Nous montrons que tous les r\'esultats de l'article restent valables quand l'indice de Poincar\'e-Lefschetz du point fixe est \'egal \`a $1$, \`a condition de supposer que le module de $h$ (d\'efini ci-dessous) est assez grand.\footnote{L'hypoth\`ese retenue est ``$\diam(h) \geq 4$''. En fait, avec un peu plus d'effort, on peut voir que la th\'eorie marche encore en supposant seulement ``$\diam(h) \geq 3$''.}
Nous reprendrons implicitement la plupart des d\'efinitions donn\'ees dans la section~\ref{s.definitions} en indice diff\'erent de $1$, en explicitant seulement celles qui posent probl\`eme. En particulier, la section suivante contient une d\'efinition du module de $h$.

Comme dans le cadre du d\'ebut de l'article, les exemples mod\`eles sont obtenus en recollant un certain nombre  de secteurs hyperboliques, elliptiques et indiff\'erents (voir l'introduction, section~\ref{ss.exemple}). Si il y a autant de secteurs hyperboliques que de secteurs elliptiques, alors l'indice du mod\`ele est \'egal \`a $1$. D'autre part, le module est \'egal au nombre total de secteurs. On construit ainsi une infinit\'e de mod\`eles d'indice $1$ et de modules sup\'erieurs ou \'egaux \`a $4$, qui auront des indices par quarts-de-tour deux \`a deux distincts.

\begin{figure}[htbp]
\begin{center}
\psset{xunit=1mm,yunit=1mm,runit=1mm}
\psset{linewidth=0.3,dotsep=1,hatchwidth=0.3,hatchsep=1.5,shadowsize=1}
\psset{dotsize=0.7 2.5,dotscale=1 1,fillcolor=black}
\psset{arrowsize=1 2,arrowlength=1,arrowinset=0.25,tbarsize=0.7 5,bracketlength=0.15,rbracketlength=0.15}
\begin{pspicture}(0,0)(107.5,33.75)
\psline{>-<}(37.5,30)(37.5,5)
\psline{<->}(25,17.5)(50,17.5)
\psline{<->}(78.75,16.25)(103.75,16.25)
\psline{>-<}(91.25,28.75)(91.25,3.75)
\psbezier(25,20)(32.5,20)(35,22.5)(35,30)
\psbezier(25,18.75)(36.25,18.75)(36.25,18.75)(36.25,30)
\psbezier(25,21.25)(28.75,21.25)(33.75,26.25)(33.75,30)
\psbezier(88.75,3.75)(88.75,11.25)(86.25,13.75)(78.75,13.75)
\psbezier(90,3.75)(90,15)(90,15)(78.75,15)
\psbezier(87.5,3.75)(87.5,7.5)(82.5,12.5)(78.75,12.5)
\psbezier(50,15)(42.5,15)(40,12.5)(40,5)
\psbezier(50,16.25)(38.75,16.25)(38.75,16.25)(38.75,5)
\psbezier(50,13.75)(46.25,13.75)(41.25,8.75)(41.25,5)
\psbezier(78.75,18.75)(86.25,18.75)(88.75,21.25)(88.75,28.75)
\psbezier(78.75,17.5)(90,17.5)(90,17.5)(90,28.75)
\psbezier(78.75,20)(82.5,20)(87.5,25)(87.5,28.75)
\psbezier{->}(91.25,16.25)(107.5,16.25)(107.5,22.5)(102.5,27.5)
\psbezier{->}(91.25,16.25)(102.5,16.25)(103.75,21.25)(100,25)
\psbezier{->}(91.25,16.25)(97.5,16.25)(101.25,18.75)(97.5,22.5)
\psbezier(91.25,16.25)(91.25,32.5)(97.5,32.5)(102.5,27.5)
\psbezier(91.25,16.25)(91.25,27.5)(96.25,28.75)(100,25)
\psbezier(91.25,16.25)(91.25,22.5)(93.75,26.25)(97.5,22.5)
\psbezier{->}(91.25,16.25)(107.5,16.25)(107.5,10)(102.5,5)
\psbezier{->}(91.25,16.25)(102.5,16.25)(103.75,11.25)(100,7.5)
\psbezier{->}(91.25,16.25)(97.5,16.25)(101.25,13.75)(97.5,10)
\psbezier(91.25,16.25)(91.25,0)(97.5,0)(102.5,5)
\psbezier(91.25,16.25)(91.25,5)(96.25,3.75)(100,7.5)
\psbezier(91.25,16.25)(91.25,10)(93.75,6.25)(97.5,10)
\psbezier{->}(37.5,17.5)(53.75,17.5)(53.75,23.75)(48.75,28.75)
\psbezier{->}(37.5,17.5)(48.75,17.5)(50,22.5)(46.25,26.25)
\psbezier{->}(37.5,17.5)(43.75,17.5)(47.5,20)(43.75,23.75)
\psbezier(37.5,17.5)(37.5,33.75)(43.75,33.75)(48.75,28.75)
\psbezier(37.5,17.5)(37.5,28.75)(42.5,30)(46.25,26.25)
\psbezier(37.5,17.5)(37.5,23.75)(40,27.5)(43.75,23.75)
\psbezier{->}(37.5,17.5)(21.25,17.5)(21.25,11.25)(26.25,6.25)
\psbezier{->}(37.5,17.5)(26.25,17.5)(25,12.5)(28.75,8.75)
\psbezier{->}(37.5,17.5)(31.25,17.5)(27.5,15)(31.25,11.25)
\psbezier(37.5,17.5)(37.5,1.25)(31.25,1.25)(26.25,6.25)
\psbezier(37.5,17.5)(37.5,6.25)(32.5,5)(28.75,8.75)
\psbezier(37.5,17.5)(37.5,11.25)(35,7.5)(31.25,11.25)
\psline{->}(31.43,25)(30.8,24.29)
\psline{->}(43.39,9.91)(44.37,10.8)
\psline{->}(85.36,24.02)(84.46,22.95)
\psline{->}(85.27,8.66)(84.37,9.55)
\end{pspicture}
\caption{Deux mod\`eles d'indice $1$ et de module $4$ (correspondant aux mots $(\ua \ra \da \ra \ua \ra \da \ra)$ et $(\ua \la \da \ra \ua \ra \da \la)$). Notons que le premier exemple ne v\'erifie pas la conclusion du lemme de Franks (il admet des cha\^ines de disques p\'eriodiques), contrairement au second.}
\label{}
\end{center}
\end{figure}
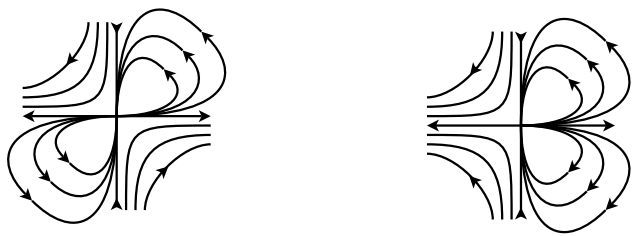

%%%%%%%%%%%%%%%%%%%%%
\subsection{Courbes de Jordan g\'eod\'esiques}
\label{ss.jordan-geodesiques}
%Dans cette section, on consid\`ere un hom\'eomorphisme du plan $h$, pr\'eservant l'orientation, fixant uniquement le point $0$ (sans hypoth\`ese d'indice).
Lorsque l'indice de Poincar\'e-Lefschetz du point fixe est diff\'erent de $1$, les notions de d\'ecomposition et de $h$-longueur d'une courbe ferm\'ee ont \'et\'e d\'efinies \`a la section~\ref{ss.definition-courbes}. On adopte les m\^emes d\'efinitions lorsque l'indice est \'egal \`a $1$. La caract\'erisation du module donn\'ee (en indice diff\'erent de $1$) par la premi\`ere partie de la proposition~\ref{p.module-courbes} permet alors de g\'en\'eraliser la notion de module. 
\begin{defi}
Le module de $h$ est le 
minimum des $h$-longueurs 
des courbes ferm\'ees $\gamma$ de degr\'e $1$.
\end{defi}
On reprend aussi la d\'efinition des courbes g\'eod\'esiques.
En recopiant la preuve de la deuxi\`eme partie de la  proposition~\ref{p.module-courbes}, (section~\ref{s.module-courbes}), on obtient le r\'esultat plus g\'en\'eral suivant.
\begin{prop}
\label{p.jordan-geodesique}
Soit $h$ un hom\'eomorphisme du plan, pr\'eservant l'orientation, fixant 
uniquement le point $0$.
Alors il existe une courbe de Jordan g\'eod\'esique pour $h$.
\end{prop}

%%%%%%%%%%%%%%%%%%%%%%%
\subsection{Construction du relev\'e canonique}
\label{ss.releve-canonique-1}

Pour mimer la construction en indice diff\'erent de $1$, il nous faut maintenant \'etablir l'existence d'un relev\'e canonique de $h$, qui satisfasse les m\^emes propri\'et\'es  qu'avant. Ceci est l'objet de la pr\'esente section.
%On consid\`ere un hom\'eomorphisme du plan, pr\'eservant l'orientation, fixant le point $0$.

%Soit $\pi$  le rev\^etement universel
%$$
%\begin{array}{rcl}
%\bbR^2 & \longrightarrow & \bbC \setminus \{0\} \simeq \bbR^2 \setminus \{0\} \\
%(x,y) & \longmapsto & \exp(-y + 2i\pi x).
%\end{array}
%$$
%\footnote{??? Le signe $-$ est une commodit\'e, pour des raisons
%visuelles : le point $0$ se
%retrouve ``\`a l'infinie en haut'' du rev\^etement. Par ailleurs, la lettre  ``$\pi$'' \`a l'int\'erieur de la formule n'a bien s\^ur rien \`a voir avec la lettre  ``$\pi$'' d\'esignant l'application de rev\^etement.}
%Un \emph{relev\'e} d'un point $x \neq 0$ est un point $\wt x$ tel que
%$\pi(\wt x)=x$. Un \emph{relev\'e} d'une courbe $\gamma : [0,1] \ra \bbR^2
%\setminus\{0\}$ est une courbe $\wt \gamma : [0,1] \ra \bbR^2$ telle que $\pi
%\circ \wt \gamma=\gamma$. Un \emph{relev\'e} de l'hom\'eomorphisme $h$
% est un hom\'eomorphisme  $\wt h$ du plan $\bbR^2$
% tel que  $h \circ \pi=\pi \circ \wt h$.

%On note $\tau$ l'automorphisme de rev\^etement $(x,y) \mapsto
%(x+1,y)$.
% \'Etant donn\'e un relev\'e $\wt h $ de $h$, les autres relev\'e
%de $h$ sont les hom\'eomorphismes $\tau^k \circ \wt h$ o\`u $k$
%parcourt l'ensemble des entiers. D'autre part, comme $h$ pr\'eserve l'orientation, il est isotope \`a l'identit\'e parmi les hom\'eomorphismes fixant $0$ (th\'eor\`eme de Kneser, ~\cite{}), et par cons\'equent, tout relev\'e de $h$ commute avec la translation $\tau$.
 Rappelons que pour choisir un relev\'e de $h$ au rev\^etement universel, il suffit de
choisir un point $x$ et une courbe $\gamma$ reliant $x$ \`a $h(x)$~:
il existe alors un unique relev\'e $\wt h $ de $h$ tel que pour tout
relev\'e $\wt x $ de $x$, le relev\'e  de $\gamma$ issu de $\wt x$ a pour
autre extr\'emit\'e $\wt h (\wt x)$~; on dira alors que $\wt h$ est le
\emph{relev\'e associ\'e \`a $\gamma$}. Les arcs de translation
 sont des cas particuliers de courbes reliant un point \`a son image. 
\begin{defi}
Un arc $\gamma$ du plan est un
\emph{arc de translation} pour $h$ (au point $x$) si $\gamma$ va de
$x$ \`a $h(x)$, et si $\gamma \cap h(\gamma) = \{ h(x) \}$.
\end{defi}

%On note $\pi : \bbR^2 \rightarrow \bbR^2 \setminus \{0\}$ le
%rev\^etement universel, $\tau$ l'automorphisme de rev\^etement,
%qui est (conjugu\'e \`a) une translation ($\tau$ est choisie de mani\`ere
%compatible avec l'orientation de $\bbR^2 \setminus\{0\}$~: plus
%pr\'ecis\'ement, soit $\gamma$ une courbe de degr\'e $1$ dans $\bbR^2
%\setminus\{0\}$, soit $\wt \gamma$  une application de $\bbR$ dans
%$\bbR^2$ qui rel\`eve $\gamma$ (au sens o\`u $p \circ \wt \gamma =
%\gamma$)~; alors $\wt \gamma(1)=\tau(\wt \gamma(0))$).

% On dira qu'un hom\'eomorphisme $\wt h $ du plan
%est un relev\'e de l'hom\'eomorphisme $h$ si $\pi \wt h = h \pi$.
% Dans ce
%cas, $\wt h$ commute avec la translation $\tau$.
%L'indice de $\wt h$ le long d'une courbe allant d'un point
%$x$ au point $\tau(x)$ ne d\'epend pas de la courbe choisie, on le note
%$\indi(\wt h, \tau)$. La proposition suivante est contenu dans le texte 
%\cite{0} (proposition 12.18 et ...). \footnote{Il ne s'agit pas d'un
%r\'esultat tr\`es difficile~: dans \cite{0}, il suffit de lire les
%sections 16 et 17 ???}

%%%%%%%%%%%%%%%%%%%%%%
\begin{prop}
\label{p.releve-canonique-1}
Soit $h$ un hom\'eomorphisme du plan
pr\'eservant l'orientation, fixant le point $0$.
 On suppose que 
$$
\diam(h) \geq 4.
$$ 
Alors il existe un unique relev\'e $\wt h$ de $h$ par $\pi$, appel\'e \emph{relev\'e canonique de $h$}, pour lequel tout relev\'e de tout arc de translation de $h$ est un arc de translation de $\wt h$.
De plus, 
\begin{enumerate}
\item  l'image par $\pi$ de tout arc libre pour $\wt h$ est un arc libre
pour $h$ ;
\item si $\cF$ est une d\'ecomposition en briques pour $h$, alors $\wt \cF := \pi^{-1}(\cF)$ est une d\'ecomposition en briques pour $\wt h$ ;
\item on a la formule d'indice 
$$
\indi_{\tau}(\wt h) = \indi(h,0) -1.
$$
\end{enumerate}
%Enfin, la propri\'et\'e 1 caract\'erise \'egalement le relev\'e $\wt h$.
\end{prop}
%\footnote{Donner la relation d'indice tout de suite.}
%\begin{defi}
%Cet hom\'eomorphisme $\wt h$ est appel\'e \emph{relev\'e canonique} de $h$.
%\end{defi}

Nous prouvons d'abord la premi\`ere partie de l'\'enonc\'e, les propri\'et\'es 1 et 2 seront d\'emontr\'ees dans la section suivante, la propri\'et\'e $3$ sera obtenue seulement \`a la section~\ref{ss.resultats-1}.

\begin{demo}[de la premi\`ere partie de la proposition~\ref{p.releve-canonique-1} (structure)]
%Soit $h$ un \'el\'ement de ${\cal H}_0$ de module sup\'erieur
%ou \'egal \`a $4$.
Nous donnons d'abord la structure de la preuve, fond\'ee sur trois affirmations. Nous d\'emontrons ensuite les trois affirmations.
La premi\`ere \'etape de la preuve consiste \`a montrer qu'il existe des
arcs de translation pour $h$.
\begin{affi}\label{affi.bon-releve-1}
Pour tout point $x \neq 0$, il existe un arc de translation pour $h$
au point $x$.
\end{affi}
On peut maintenant choisir un point $x_{0}$, un arc de translation $\gamma$ au point $x_{0}$, et appeler $\wt h$ le relev\'e associ\'e \`a $\gamma$.
Il s'agit ensuite de voir que tous les arcs de translation induisent
le m\^eme relev\'e. Nous fixons d'abord le point $x_{0}$, et nous
montrons que $\wt h$ ne d\'epend pas du choix de $\gamma$. Ceci revient
\`a s'int\'eresser aux classes d'homotopie des arcs de translation, ce que fait la deuxi\`eme affirmation. 
\begin{affi}\label{affi.bon-releve-2}
Soient $\gamma_1$ et $\gamma_2$ deux arcs de translation tels que
$\gamma_1(0)=\gamma_2(0)$ (et donc aussi
$\gamma_1(1)=\gamma_2(1)$). Alors ces deux arcs sont homotopes, \`a
extr\'emit\'es fix\'ees, dans $\bbR^2 \setminus \{0\}$.
\end{affi}
Il reste enfin \`a v\'erifier que le relev\'e $\wt h$  ne d\'epend pas non plus du choix du point $x_{0}$. Par connexit\'e de l'espace $\bbR^2 \setminus \{0\}$, il suffit de voir ceci  localement. Ceci revient \`a construire, localement,  une famille continue d'arcs de translation, et c'est l'objet de la derni\`ere affirmation.
\begin{affi}\label{affi.bon-releve-3}
Soit $x_0$ un point de $\bbRo$. Il existe alors un voisinage $V$ de $x_0$,
et une famille continue $(\gamma_x)_{x \in V}$ d'arcs tels que, pour
tout point $x$ de $V$, l'arc $\gamma_x$ est un arc de translation
au point $x$.
\end{affi}
La continuit\'e de la famille $(\gamma_x)$ s'entend au sens de la topologie de la
convergence uniforme sur l'espace des arcs.
% Notons maintenant $\wt h$ le relev\'e associ\'e \`a $\gamma_{x_{0}}$~; il s'agit de montrer que $\wt h$est aussi le relev\'e associ\'e \`a chaque arc $\gamma_{x}$ pour $x$ dans $V$.
%Soit  $\wt V$ une composante connexe de $\pi^{-1}(V)$. On peut choisir $V$ simplement connexe, et dans ce cas chaque point $x$ de $V$ poss\`ede un unique relev\'e $\wt x$ dans $\wt V$. La famille d'arcs $(\gamma_x)_{x \in V}$ se rel\`eve en une famille continue $(\wt \gamma_x)_{x \in V}$ d'arcs, chaque arc $\wt \gamma_x$ allant du point $\wt x$ au point $\wt \gamma_{x}(1)$, qui est un relev\'e de $h(x)$, c'est. L'application $x \mapsto \wt \gamma_{x}(1)$ est continue, et $\wt \gamma_{x_{0}}(1)=\wt h(\wt x_{0})$, ce qui implique que, pour tout $x$, $\wt \gamma_{x}(1) = \wt h(\wt x)$. C'est  ce que l'on voulait montrer.
Ceci termine la preuve de la premi\`ere partie de la proposition \`a partir des trois affirmations.
\end{demo}

Nous aurons besoin du lemme suivant (qui est \`a rapprocher des techniques de la  section~\ref{s.module-courbes}).
\begin{lemm}\label{l.degre-grand}
Soit $\gamma :  [0,1] \ra \bbR^2 \setminus \{0\}$ une courbe ferm\'ee. Si $\gamma$ est de degr\'e $> 1$, alors il existe deux r\'eels $t_{1}, t_{2}$, $0 \leq t_{1} <  t_{2} \leq 1$,  tels que $[t_{1} t_{2}]_{\gamma}$ est une courbe ferm\'ee de degr\'e $1$.
En particulier, la $h$-longueur de $\gamma$ est sup\'erieure ou \'egale au module de $h$.
\end{lemm}
En r\'ealit\'e, on pourra montrer (\textit{a posteriori}) que la $h$-longueur de $\gamma$ est minor\'ee par le module de $h$ multipli\'e par la valeur absolue du degr\'e de $\gamma$.

\begin{demo}[du lemme]
Soit $\wt \gamma : [0,1] \ra \bbR^2 $ une courbe du plan qui rel\`eve $\gamma$.
Notons $k>1$ le degr\'e de $\gamma$~; par hypoth\`ese, on a $\wt \gamma (1) = \tau^k(\wt \gamma(0))$.
Or si $\tau(\wt \gamma)$ ne rencontrait pas $\wt \gamma$, alors aucun des it\'er\'es de $\gamma$ par $\tau$ ne rencontrerait $\gamma$ (on peut voir ceci comme une cons\'equence du lemme de Franks appliqu\'e \`a la translation $\tau$). On en d\'eduit le r\'esultat.
\end{demo}

Nous prouvons maintenant les trois affirmations.

\begin{demo}[de l'affirmation~\ref{affi.bon-releve-1}]
Soit $x \neq 0$, on cherche un arc de translation au point $x$.
On va r\'eutiliser une construction classique des arcs de translation pour les hom\'eomorphismes de Brouwer. 
 Quitte \`a conjuguer, on peut supposer que le point $x$ est plus proche de son image par $h$ que du point fixe $0$ (pour la distance euclidienne).
On consid\`ere maintenant un  disque euclidien $D$, centr\'e en $x$, de rayon $r$. Si $r$ est assez petit, le disque est libre (disjoint de son image)~; si $r$ est \'egal \`a  la distance de $x$ \`a $h(x)$, il ne l'est plus ; on peut donc choisir pour $r$ la plus petite valeur pour laquelle $D$ n'est pas libre. Le disque $D$ ne contient pas le point fixe $0$, et  est un disque \emph{critique}, c'est-\`a-dire qu'il v\'erifie
\begin{itemize}
\item $D$ n'est pas libre,
\item l'int\'erieur de $D$ est libre.
\end{itemize}
Il existe donc un point $y$ sur le bord du disque $D$ dont l'image $h(y)$ est aussi sur
le bord de $D$. On d\'efinit les segments euclidiens $\alpha_{1}$ reliant $y$ \`a $x$, et $\alpha_{2}$ reliant $x$ \`a $h(y)$, et on consid\`ere l'arc $\alpha$ obtenu en concat\'enant $\alpha_{1}$ et $\alpha_{2}$, arc qui va de $y$ \`a son image $h(y)$.

\begin{figure}[htbp]
\begin{center}
\def\JPicScale{0.7}
\ifx\JPicScale\undefined\def\JPicScale{1}\fi
\psset{unit=\JPicScale mm}
\psset{linewidth=0.3,dotsep=1,hatchwidth=0.3,hatchsep=1.5,shadowsize=1}
\psset{dotsize=0.7 2.5,dotscale=1 1,fillcolor=black}
\psset{arrowsize=1 2,arrowlength=1,arrowinset=0.25,tbarsize=0.7 5,bracketlength=0.15,rbracketlength=0.15}
\begin{pspicture}(0,0)(107.32,62.5)
\rput{0}(50,39.29){\psellipse[](0,0)(20,20)}
\pscustom[]{\psbezier(70,39.29)(70,44.29)(78.04,43.21)(76.96,48.57)
\psbezier(76.96,48.57)(75.89,53.93)(66.25,51.07)(63.57,53.93)
\psbezier(63.57,53.93)(60.89,56.79)(64.46,59.82)(68.04,59.46)
\psbezier(68.04,59.46)(71.61,59.11)(86.61,62.5)(95.71,57.14)
\psbezier(95.71,57.14)(104.82,51.79)(106.79,42.68)(105.71,39.11)
\psbezier(105.71,39.11)(104.64,35.54)(105.54,23.39)(96.96,23.57)
\psbezier(96.96,23.57)(88.39,23.75)(77.58,25.8)(74.64,27.5)
\psbezier(74.64,27.5)(71.7,29.2)(70,34.29)(70,39.29)
\closepath}
\psdots[](70,39.29)
(70,39.29)
\psdots[](50,39.29)
(50,39.29)
\psdots[](33.75,51.25)
(33.75,51.25)
\psline[linestyle=dashed,dash=1 1](33.75,51.25)(50,39.29)
\psdots[linestyle=dashed,dash=1 1](50,39.29)
(70,39.29)
\psbezier[linestyle=dashed,dash=1 1](95,34.29)(95,34.29)(95,34.29)(95,34.29)
\psline[linestyle=dashed,dash=1 1](50,39.29)(70,39.29)
\psdots[linestyle=dashed,dash=1 1](95,34.29)
(95,34.29)
\pscustom[linestyle=dashed,dash=1 1]{\psbezier(70,39.29)(80,34.29)(80,39.29)(85,44.29)
\psbezier(85,44.29)(90,49.29)(95,34.29)(95,34.29)
}
\rput(49.29,34.64){$x$}
\rput(31.25,53.75){$y$}
\rput(95,29.29){$h(x)$}
\rput(63.93,35){$h(y)$}
\rput(45,47.5){$\alpha_1$}
\rput(45,49.29){}
\rput(60,42.5){$\alpha_2$}
\rput(87,49){$h(\alpha_1)$}
\rput(57.5,46.79){}
\rput(50,14){$\alpha$}
\rput(72.5,3.75){$\beta$}
\pscustom[]{\psbezier(30,18.75)(30,12.5)(48.75,17.5)(48.75,17.5)
\psline(48.75,17.5)(50,16.25)
\psline(50,16.25)(51.25,17.5)
\psbezier(51.25,17.5)(51.25,17.5)(70,12.5)(70,18.75)
}
\pscustom[]{\psbezier(50,11.25)(50,5)(71.09,10)(71.09,10)
\psline(71.09,10)(72.5,8.75)
\psline(72.5,8.75)(73.91,10)
\psbezier(73.91,10)(73.91,10)(95,5)(95,11.25)
}
\rput(25,40){$D$}
\rput(25,40){}
\rput(107.32,55){$h(D)$}
\end{pspicture}
\caption{Construction d'un arc de translation}
\label{f.arc-de-translation}
\end{center}
\end{figure}
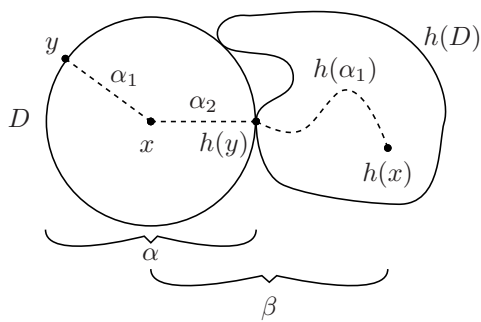

Soit $A$ la courbe infinie obtenue en concat\'enant tous les it\'er\'es de $\alpha$~:
$$
A = \cdots h^{-1}(\alpha) \star \alpha \star h(\alpha) \star \cdots .
$$
Voici l'id\'ee de la fin de la preuve. Supposons que cette courbe infinie $A$ ait un point double, elle contient alors une courbe de Jordan $J$. \`A l'aide de la th\'eorie de Brouwer, on va montrer que la courbe $J$ doit entourer le point fixe $0$ ; mais alors, si le module de $h$ est assez grand, $J$ doit contenir suffisamment d'it\'er\'es de l'arc ``presque libre'' $\alpha$ ; et l'arc $\alpha_{2} \star h(\alpha_{1})$ sera un arc de translation au point $x$.
\footnote{En r\'ealit\'e, on peut montrer que l'hypoth\`ese de module entra\^ine m\^eme que $A$ est une courbe sans point double (nous n'utiliserons pas ceci).}

R\'ealisons cette id\'ee. Soit $\wt h$ l'hom\'eomorphisme relev\'e de $h$ associ\'e \`a l'arc $\alpha$. On consid\`ere un disque topologique $\wt D$, composante connexe de $\pi^{-1}(D)$~; on note $\wt y$ et $\wt \alpha$ les relev\'es respectifs de $y$ et $\alpha$ inclus dans $\wt D$. Le disque $\wt D$ est un disque critique pour $\wt h$.
Comme $\wt h$ n'a pas de point fixe, on dispose du lemme de Franks (voir par exemple~\cite{0}, paragraphe 3.2.B p35)~: notamment, $\wt h$ n'a pas de point p\'eriodique, et en particulier les it\'er\'es du point $\wt y$ sont tous distincts ; et tout disque topologique (ouvert) libre est disjoint de ses it\'er\'es, et en particulier  les it\'er\'es de l'int\'erieur du  disque $\wt D$ sont deux \`a deux disjoints. On en d\'eduit que la courbe infinie $\wt A$ obtenue en concat\'enant les it\'er\'es de $\wt \alpha$ est une courbe sans point double.

On d\'efinit maintenant la courbe $\beta = \alpha_{2} \star h(\alpha_{1})$. On voit facilement que cette courbe est un arc (elle est sans point double), qui va de $x$ \`a $h(x)$.
 Montrons qu'il s'agit d'un arc de translation pour $h$. Ceci revient \`a montrer que la courbe $\beta \star h(\beta)$ est une courbe simple.
Donnons-nous un param\'etrage de cette courbe par l'intervalle $[0,1]$, qui envoie respectivement les intervalles 
$$
\left[0, \frac{1}{4}\right], \left[\frac{1}{4},\frac{2}{4}\right],  
\left[\frac{2}{4},\frac{3}{4}\right], \left[\frac{3}{4},\frac{4}{4}\right],
$$
 sur les courbes $\alpha_{2}$, $h(\alpha_{1})$, $h(\alpha_{2})$, $h^2(\alpha_{1})$.
Les arcs $\alpha_{2}$ et $h^2(\alpha_{1})$ sont libres~; il existe donc deux r\'eels $t_{1}$ et $t_{2}$, v\'erifiant 
$$
\frac{1}{4} < t_{1} < t_{2} < \frac{3}{4}
$$
tels que les courbes $[0, t_{1}]_{\beta \star h(\beta)}$ et $[t_{2}, 1]_{\beta \star h(\beta)}$
sont encore  libres. D'autre part, la courbe $[t_{1}, t_{2}]_{\beta \star h(\beta)}$ est libre car elle est incluse dans l'int\'erieur du disque $h(D)$.
La courbe $\beta \star h(\beta)$ s'\'ecrit dont comme concat\'enation de trois courbes libres~: sa $h$-longueur est sup\'erieure ou \'egal \`a $3$. 

D'autre part, remarquons que les courbes 
$$
 [0, \frac{3}{4}]_{\beta \star h(\beta)} = \alpha_{2} \star h(\alpha_{1}) \star h(\alpha_{2})  = \mbox{ et } 
  [\frac{1}{4},1]_{\beta \star h(\beta)} = h(\alpha_{1}) \star h(\alpha_{2}) \star h^2(\alpha_{1})
$$
sont sans point double (ceci provient de la libert\'e de l'int\'erieur du disque $D$).

Supposons maintenant, par l'absurde, que la courbe $\beta \star h(\beta)$ poss\`ede un point double~:
 il existe deux r\'eels $s_{1}$ et $s_{2}$, v\'erifiant $0 < s_{1} < s_{2} < 1$, tels que la courbe $J=[s_{1}, s_{2}]_{\beta \star h(\beta)}$ est une courbe ferm\'ee. D'apr\`es la remarque pr\'ec\'edente, on a n\'ecessairement $s_{1} < \frac{1}{4}$ et $\frac{3}{4} < s_{2}$. En choisissant  une telle courbe $J$ minimale, on peut supposer que $J$ est une courbe de Jordan.
%\footnote{En g\'en\'eral, une courbe ayant des points doubles ne contient pas n\'ecessairement de courbe de Jordan~; il est donc essentiel de minorer l'\'ecart entre les param\`etres $s_{1}$ et $s_{2}$.}
 Le degr\'e de $J$ dans $\bbR^2 \setminus \{0\}$ vaut alors $-1$, $0$ ou $1$.
 %\pourmoi{PREUVE ??? R\'eutilis\'e section ``equivalence des donnes''.}
Or ce degr\'e ne peut pas \^etre nul, car $J$ est inclus dans la courbe infinie $A$ dont un relev\'e $\wt A$ est sans point double.
Quitte \`a changer le sens de parcourt de $J$, on peut supposer qu'elle est de degr\'e $1$.
D'autre part, la $h$-longueur de $J$ est inf\'erieure ou \'egale \`a trois puisque $J$ est une sous-courbe de la courbe $\beta \star h(\beta)$.
Ceci contredit l'hypoth\`ese que le module de $h$ est sup\'erieur ou \'egal \`a quatre.
\end{demo}

Nous en profitons pour \'enoncer un corollaire utile.
\begin{coro}
\label{c.periode-2}
$h$ n'a pas de point p\'eriodique de p\'eriode $2$.
\end{coro}
\begin{demo}
D'apr\`es la d\'efinition des arcs de translation, 
si $x$ \'etait un point de période $2$, alors il n'existerait pas d'arc de translation au point $x$, contrairement \`a ce qu'on vient de montrer.
\end{demo}

\begin{demo}[de l'affirmation~\ref{affi.bon-releve-2}]
Soit $\gamma$ la courbe obtenue en concat\'enant $\gamma_{1}$ et
l'arc $\gamma_{2}$ parcouru en sens contraire. 
Les arcs $\gamma_{1}$ et $\gamma_{2}$ sont  homotopes si et seulement si la courbe $\gamma$ est de degr\'e nul dans $\bbR^2 \setminus \{0\}$.
Dans le cas contraire, la $h$-longueur de $\gamma$ devrait \^etre sup\'erieure ou \'egale au module de $h$ (lemme~\ref{l.degre-grand}). Or on voit facilement que la $h$-longueur de $\gamma$ est inf\'erieure ou \'egale \`a trois.
Ceci prouve l'affirmation.
\end{demo}

\begin{demo}[de l'affirmation~\ref{affi.bon-releve-3}]
Il s'agit simplement de g\'en\'eraliser la construction de la preuve de la premi\`ere affirmation. On note \`a nouveau $D$ le disque critique centr\'e au point $x_{0}$, et $y$ un point du bord de $D$ dont l'image $h(y)$ est aussi sur le bord de $D$.
Pour tout point $x$ de $D$, on d\'efinit les segment $\alpha_{1}^x$ reliant $y$ \`a $x$, et $\alpha_{2}^x$ reliant $x$ \`a $h(x)$, et on consid\`ere l'arc $\beta^x = \alpha_{2}^x \star
h(\alpha_{1}^x)$. On montre, \`a l'aide des arguments utilis\'es dans la premi\`ere affirmation, que $\beta^x$ est un arc de translation au point $x$. D'autre part, l'application $x \mapsto \beta^x$ est clairement continue. La preuve est termin\'ee.
\end{demo}

%%%%%%%%%%%%
%\pourmoi{\subsection{Propri\'et\'es et corollaires dynamiques}
%UTILE ??

%Il faudrait peut-etre \'enoncer un lemme de Fransk
%restreint : il n'existe pas de chaine de disque de longueur plus
%petite que le module (ou juste de longueur 2) ??? permettrait peut-etre de d\'efinir
%l'invfariant en bas ??

%
%\begin{lemm}\label{lemm.propriete-bon-releve-1}
%On consid\`ere un \'el\'ement $h$ de $\ho$ v\'erifiant $\diam(h) \geq
%4$, et son bon
%relev\'e $\wt h$.
%Si $\gamma$ est un arc libre pour $\wt h$, alors $\pi(\gamma)$ est un arc
%libre pour $h$.
%\end{lemm}

%\begin{demo}[du lemme~\ref{lemm.propriete-bon-releve-1}]
%...
%\end{demo}
%
%Montrer que la conclusion du lemme de Franks est v\'erifi\'ee ? Utile ? Au moins, le mentionner. Preuve ??

%%%%%%%%%%%%%%%%%%%%%%%
\subsection{Propri\'et\'es et cons\'equences dynamiques}
Dans cette section, nous commen\c cons par d\'emontrer la propri\'et\'e 1 et  de la proposition~\ref{p.releve-canonique-1}, puis nous en d\'eduirons quelques cons\'equences dynamiques (notamment, l'errance des points).
Le rel\`evement des d\'ecompositions en briques (propri\'et\'e 2 de la proposition) ne pose pas de probl\`eme, le lecteur peut se r\'ef\'erer \`a ~\cite{0}, preuve du troisi\`eme point de la proposition 4.17, p98.
Pour obtenir la formule d'indice (troisi\`eme point de la proposition), il faut d'abord obtenir une droite de Brouwer pour $h$ ; la preuve de cette formule est donc repouss\'ee \`a la section~\ref{ss.resultats-1}.

%%%%%%%%%%%%%%%%%%%%%%%
\subsubsection*{Arcs libres, d\'ecompositions}
\begin{demo}[de la premi\`ere propri\'et\'e de la proposition~\ref{p.releve-canonique-1}]
On se place sous les hypoth\`eses de la proposition.
Soit $\gamma$ un arc libre pour $\wt h$, il s'agit de voir que $\pi(\gamma)$ est libre, et est sans point double.
Si $\pi(\gamma)$ n'\'etait pas libre, il contiendrait une courbe joignant un point $x$ \`a son image $h(x)$, donc \'egalement un arc $\alpha$ joignant un point \`a son image ; on trouverait facilement un sous-arc de $\alpha$ joignant un point \`a son image, et minimal pour cette propri\'et\'e, autrement dit tel que l'int\'erieur de $\alpha$ est libre.
% (voir par exemple~\cite{} lucien ???). 
Comme $h$ n'a pas de point de p\'eriode $2$ (corollaire~\ref{c.periode-2}), on en d\'eduirait que $\alpha$ est un arc de translation. D'autre part, $\alpha$ a un relev\'e inclus dans $\gamma$, ce relev\'e est donc libre, ce n'est pas un arc de translation, et ceci contredit la d\'efinition du relev\'e canonique.

La courbe $\pi(\gamma)$ est donc libre. Si elle avait un point double, elle contiendrait une courbe de Jordan, qui ne pourrait pas \^etre de degr\'e $0$ puisque son relev\'e n'est pas une courbe ferm\'ee ; on aurait donc une courbe de Jordan de degr\'e $1$, libre, ce qui contredirait l'hypoth\`ese sur le module de $h$.
Finalement, $\pi(\gamma)$ est bien un arc libre.
\end{demo}

%%%%%%%%%%%%%%%%%%%%%%%
%\subsubsection*{Formule d'indice}
%Pour obtenir la formule d'indice (troisi\`eme point de la proposition~\ref{p.releve-canonique-1}), il faut d'abord obtenir une droite de Brouwer pour $h$ ; la preuve de cette formule est donc repouss\'ee \`a la section~\ref{ss.resultats-1}.

%%%%%%%%%%%%%%%%%%%%%%%
\subsubsection*{Libert\'e et errance}
Sous l'hypoth\`ese que le module de $h$ est plus grand que $4$, nous allons montrer que la dynamique de $h$ ressemble \`a celle du cas o\`u le point fixe est d'indice diff\'erent de $1$.

Consid\'erons un hom\'eomorphisme $H$ d'une surface.
Un ensemble $E$ est \emph{libre} pour $H$ si il est
disjoint de son image $H(E)$.
Suivant Morton Brown (\cite{brown}), on dira que  $H$
est \emph{libre} si tout disque topologique ferm\'e
 $D$ qui est libre est
encore disjoint de tous ses it\'er\'es $H^n(D)$, $n \neq 0$. Les
orbites d'un hom\'eomorphisme libre sont particuli\`erement simple~:
en effet, tout point $x$ qui n'est pas fixe poss\`ede un voisinage $D$
qui est un disque topologique ferm\'e libre. Le disque $D$ est alors
disjoint de tous ses it\'er\'es, par cons\'equent le point $x$ est
\emph{errant} (par d\'efinition). Si les points fixes de $H$ sont
isol\'es, on en d\'eduit que l'orbite positive de $x$ tend vers un point fixe
de $H$ ou bien sort de tout compact.

\begin{coro}\label{coro.errent}
Sous les hypoth\`eses de la proposition~\ref{p.releve-canonique-1},
l'hom\'eomorphisme $h$ est libre. En particulier, tout  point est errant, 
et l'orbite positive (ou n\'egative)
d'un point du plan tend vers $0$ ou vers le point \`a l'infini.
\end{coro}

\begin{demo}[du corollaire~\ref{coro.errent}]
Soit $D$ un disque topologique libre pour $h$. Soit $\wt D$ une composante connexe de $\pi^{-1}(D)$ : c'est aussi un disque topologique.
La propri\'et\'e $1$ du relev\'e canonique entra\^ine que ce disque $\wt D$ est libre pour $\wt h$. Supposons  que le disque $D$ rencontre un it\'er\'e $h^n(D)$, 
alors il existe un entier $k$ tel que $\tau^k(\wt D)$ rencontre $\wt h^n(\wt D)$. 
L'hom\'eomorphisme $\wt h$ est un hom\'eomorphisme de Brouwer, il est donc libre, par cons\'equent le disque $\wt D$, qui est libre, est disjoint de tous ses it\'er\'es par $\wt h$ ; donc $k \neq 0$.
La th\'eorie de Brouwer nous dit que $\wt D$ est inclus dans un domaine de translation ; ce domaine est invariant, il rencontre donc $\tau^k(\wt D)$, et en particulier il existe un arc libre $\wt \gamma$ qui va d'un point de $\wt D$ \`a un point de $\tau^k(\wt D)$. On peut prolonger cet arc par un second arc libre $\wt \delta$, inclus dans $\tau^k(\wt D)$, de mani\`ere \`a obtenir une courbe $\wt J:= \wt \gamma * \wt \delta$ qui va d'un point $\wt x $ au point $\tau^k(\wt x)$. Puisque $\wt J$ rencontre $\tau^k(\wt J)$, elle doit aussi rencontrer $\tau(\wt J)$ ($\tau$ est un hom\'eomorphisme de Brouwer, il est donc libre) : quitte \`a raccourcir $\wt J$, on peut alors supposer que $k=1$. La courbe projet\'ee $J=\pi(\wt J)$
est une courbe ferm\'ee de degr\'e $1$ compos\'ee des deux arcs $\pi(\wt \gamma)$  et $\pi(\wt \delta)$, qui sont libres d'apr\`es le premier point. Ceci contredit l'hypoth\`ese sur le module de $h$.
\end{demo}

\`A l'aide des d\'ecompositions en briques, on montrerait facilement que tout point de $\bbR^2 \setminus\{0\}$ est dans un domaine de translation pour $h$ (ce qui permet de r\'ecup\'erer la d\'efinition originelle du module de $h$, d\'efinition~\ref{d.module}).
Ce r\'esultat peut  aussi \^etre vu comme une cons\'equence de la libert\'e de $h$, \emph{via} un th\'eor\`eme de Slaminka (voir ~\cite{slam}, et aussi l'article de Guillou~\cite{guil}).

%%%%%%%%%%%%%%%%%%%%%%%
\subsection{R\'esultats}
\label{ss.resultats-1}
Dans cette section, nous expliquons pourquoi les preuves des r\'esultats \'enonc\'es dans  la section~\ref{s.definitions} fonctionnent encore en indice $1$, en module plus grand que $4$.

\begin{theo}
Dans tous les \'enonc\'es des sections~\ref{s.definitions}, ~\ref{s.ouverts-dynamiques}, ~\ref{s.croissants-petales} et ~\ref{s.autres}, on peut remplacer l'hypoth\`ese
\begin{quote}
\og\ l'indice de Poincar\'e-Lefschetz du point fixe est  diff\'erent de $1$ \fg\
\end{quote}
par l'hypoth\`ese
\begin{quote}
\og\ le module de $h$ est sup\'erieur ou \'egal \`a $4$ \fg.
\end{quote}
Ceci concerne notamment 
les th\'eor\`emes~\ref{t.compo-Reeb}, \ref{t.invariant}, \ref{t.indice-courbes}, \ref{t.lien-indice}, \ref{t.ouverts-dynamiques}, les propositions~\ref{p.dynamique-bord}, \ref{p.bords-adjacents} et le corollaire~\ref{c.conservatif}.
\end{theo}

\begin{demo}
Nous expliquons rapidement comment adapter les preuves d\'ej\`a expos\'ees en d\'etail quand l'indice \'etait suppos\'e diff\'erent de $1$. Le principal changement est le suivant : la formule reliant l'indice de $h$ \`a celui de son relev\'e canonique est obtenue plus tard (troisi\`eme \'etape ci-dessous).

\paragraph{Premi\`ere \'etape}
On montre que les r\'esultats du corollaire~\ref{c.bouts-geodesique} sont v\'erifi\'es : si $\gamma$ est une courbe de Jordan g\'eod\'esique (fournie par la proposition~\ref{p.jordan-geodesique}), alors la droite topologique $\Gamma$ qui rel\`eve $\gamma$ est une droite g\'eod\'esique. En effet, la preuve de la section~\ref{ss.geodesique-canonique} utilisait seulement la correspondance entre les arcs libres de $\wt h$ et ceux de $h$, et la minoration du module : la premi\`ere propri\'et\'e est encore v\'erifi\'ee d'apr\`es la proposition~\ref{p.releve-canonique-1}, et la seconde l'est par hypoth\`ese. On note encore $[\tau]$ la classe d'\'equivalence de la droite $\Gamma$.

\paragraph{Deuxi\`eme \'etape}
Comme \`a la section~\ref{ss.dynamique-equivariants}, on peut alors consid\'erer les composantes de Reeb et le mot cyclique associ\'es \`a $[\tau]$.
On peut recopier mot pour mot la preuve du th\'eor\`eme~\ref{t.lien-indice}-bis (section~\ref{s.lien-indice}), et on obtient la conclusion de ce th\'eor\`eme, \`a savoir  la formule d'indice pour le relev\'e canonique $\wt h$. 

\paragraph{Troisi\`eme \'etape}
Nous pouvons maintenant prouver la formule qui nous manquait, 
qui relie l'indice de $h$ et celui de $\wt h$ (troisi\`eme point de la proposition~\ref{p.releve-canonique-1}). En effet, l'\'etape pr\'ec\'edente a produit des droites de Brouwer $\wt \Delta_{k}$ pour $\wt h$, et on voit facilement que ces droites se projettent en des droites de Brouwer $\Delta_{k}$ pour $h$ qui vont de $0$ \`a $\infty$ (voir par exemple la preuve de la proposition~\ref{p.croissants}). On peut alors appliquer le lemme~\ref{l.lien-indices}, qui donne la formule voulue.

\paragraph{Fin}
On peut maintenant recopier tout le reste des preuves :
la compatibilit\'e entre les composantes de Reeb (section~\ref{s.composantes-de-reeb}), 
puis les preuves des th\'eor\`emes (section~\ref{s.preuves}). 
\end{demo}

%\pagebreak
 %%%%%%%%%%%%%%%%%%%%%%%%
%%%%%%%%%%%%%%%%%%%%%%%%
\section{Ordre cyclique}
\label{as.ordre-cyclique}
%%%%%%%%%%%%%%%%%%%%%%%
\subsection{Parties du plan}
On se donne trois parties connexes $A$, $B$ et $C$ du plan.
Si l'ensemble $A$ est  disjoint de  $B$, il est inclus dans une composante connexe du compl\'ementaire de $B$ ; on notera $O_{A}(B)$ cette composante.
Remarquons que l'adh\'erence de $O_{A}(B)$ rencontre $B$, et, en particulier, l'ensemble $O_{A}(B) \cup B$ est connexe.

\begin{defi}
On dit que \emph{$B$ s\'epare $A$ et $C$}
si les trois ensembles $A$, $B$ et $C$ sont  deux \`a deux disjoints,
 et si les composantes $O_{A}(B)$ et $O_{C}(B)$
sont distinctes (et donc disjointes). On \'ecrit alors
$$
A < B < C.
$$
\end{defi}
On montre facilement une propri\'et\'e d'antisym\'etrie : si $B$ s\'epare $A$ et $C$, alors $A$ ne s\'epare pas $B$ et $C$, autrement dit
$$
\mbox{si } A < B < C \mbox{ alors on n'a pas } B < A < C.
$$
On en d\'eduit une propri\'et\'e de transitivit\'e : si $A,B,C,D$ sont des parties connexes deux \`a deux disjointes,
$$
\mbox{si } A < B < C \mbox{ et } B < C < D \mbox{ alors } A < B < D \mbox{ et } A < C < D.
$$
Dans ce cas, on \'ecrira simplement
$$
A < B < C < D.
$$

\begin{defi}
On \'ecrira 
$$
A \leq B < C
$$
si $C$ est disjoint de $B$, et s'il existe une composante connexe $O'$ du compl\'ementaire de $B$ telle que
\begin{enumerate}
\item $A$ est inclus dans l'adh\'erence de $O'$ ;
\item $O'$ et $O_{C}(B)$ sont distinctes (et donc disjointes).
\end{enumerate}
\end{defi}
Les propri\'et\'es d'antisym\'etrie et de transitivit\'e  se g\'en\'eralisent : 
$$
\mbox{si } A \leq B < C \mbox{ alors on n'a pas } B < A < C \ \  ;
$$
$$
\mbox{si } A \leq B < C \mbox{ et } B < C \leq D \mbox{ alors } A \leq B < D \mbox{ et } A < C \leq D.
$$
On \'ecrira alors
$$
A \leq B < C \leq D.
$$
On peut g\'en\'eraliser ces \'ecritures \`a un nombre quelconque de parties connexes, par r\'ecurrence : l'in\'egalit\'e d'ordre $k$
$$
A_{1} < A_{2} < \cdots < A_{k}
$$
(avec $k \geq 4$) est v\'erifi\'ee si  les $k$  ``sous-in\'egalit\'es'' d'ordre $k-1$ sont v\'erifi\'ees.

\small
\paragraph{Application}
Consid\'erons par exemple une suite de composantes de Reeb $(\wt F_{k}, \wt G_{k})$ associ\'ee \`a une suite g\'eod\'esique $(\wt x_{k})_{k \in \bbZ}$, pour un hom\'eomorphisme de Brouwer $H$ (voir~\cite{I}, section~7.1). Alors l'ordre est donn\'e par 
$$
\cdots \leq \wt F_{k} < \wt G_{k} \leq \wt F_{k+1} < \wt G_{k+1} \leq \cdots
$$
Pour voir ceci, il s'agit de montrer que $\wt F_{k} < \wt G_{k} \leq \wt F_{k+1}$.
Or $\wt G_{k}$ s\'epare $\wt x_{k-1}$ de $\wt x_{k+2}$ (voir par exemple~\cite{I}, section 7.1) ;
l'ensemble $\wt F_{k} $ est inclus dans la composante connexe $O_{\wt x_{k-1}} (\wt G_{k})$  ;
 l'ensemble  $\wt F_{k+1} $ est inclus dans l'adh\'erence de la $H$-boule de centre $\wt x_{k+2}$ et de rayon $1$, et cette boule est disjointe de $\wt G_{k}$ (par in\'egalit\'e triangulaire), donc incluse dans  
 la composante connexe $O_{\wt x_{k+2}} (\wt G_{k})$.
\normalsize

\begin{defi}
\label{d.ofg}
Soient $F$ et $G$ deux parties ferm\'ees et connexes du plan. On appelle \emph{ouvert entre $F$ et $G$} l'ensemble 
$$
O(F,G) = O_{F}(G) \cap O_{G}(F).
$$
\end{defi}

\begin{lemm}
L'ensemble $O(F,G)$ est l'unique composante connexe du compl\'ementaire de $ F \cup G$ dont l'adh\'erence rencontre \`a la fois $F$ et $G$.
\end{lemm}
Contrairement \`a ce qui pr\'ec\`ede, ce lemme est particulier \`a la topologie du plan : la preuve utilise une version du lemme d'Alexander
(voir ~\cite{I}, th\'eor\`eme A.3).

%%%%%%%%%%%%%%%%%%%%%%%
\subsection{Parties du plan trou\'e}
%On consid\`ere maintenant des parties $F$ ferm\'ees et connexes du ``plan trou\'e'' $\bbR^2 \setminus \{0\}$, telles que $\hat F := F \cup \{0,\infty \}$ soit encore connexe. On les appellera des \emph{bonnes} parties.
On consid\`ere des parties connexes du plan trou\'e $\bbR^2 \setminus \{0\}$.
Un \emph{relev\'e} $\wt A$ d'une telle partie $A$ est une composante connexe de $\pi^{-1}(A)$, o\`u $\pi$ est le rev\^etement universel du plan trou\'e.
On note $\tau$ l'automorphisme du rev\^etement universel correspondant \`a un lacet de degr\'e $1$. Tout relev\'e $\wt A$ est ou bien \'egal \`a $\tau(\wt A)$, ou bien disjoint de  $\tau(\wt A)$.
\begin{defi}
Soient $A,B,C$ trois parties connexes du plan trou\'e, deux \`a deux disjointes. On notera
$$
A < B < C < A
$$
s'il existe des relev\'es respectifs  $\wt A, \wt B, \wt C$ de $A,B,C$
tels que 
$$
\wt A < \wt B < \wt C < \tau(\wt A) < \tau(\wt B) < \cdots.
$$
\end{defi}
Comme avant, on g\'en\'eralise ceci \`a un nombre quelconque de bonnes parties :
$A_{1} < A_{2} < \cdots < A_{k} < A_{1}$ si on peut choisir des relev\'es tels que
$\wt A_{1} < \wt A_{2} < \cdots <\wt  A_{k} <\tau(\wt A_{1}) <\tau(\wt A_{2}) < \cdots$.
On d\'efinit les relations $ A < B \leq C < A$ de mani\`ere analogue.
%\footnote{Il suffit d'avoir $\wt A < \wt B < \wt C < \tau(\wt A)$... utile ??}

On consid\`ere maintenant des parties $F$ ferm\'ees et connexes du ``plan trou\'e'' $\bbR^2 \setminus \{0\}$, telles que $\hat F := F \cup \{0,\infty \}$ soit encore connexe. On les appellera des \emph{bonnes} parties.
On peut encore d\'efinir les ouverts entre $F$ et $G$, \emph{via} le lemme suivant.
\begin{lemm}
Soient $F$ et $G$ deux bonnes parties de $\bbR^2 \setminus \{0\}$, disjointes.
Il existe alors exactement deux composantes connexes $O(F,G)$ et $O(G,F)$ du compl\'ementaire de $\hat F \cup \hat G$ dont l'adh\'erence rencontre \`a la fois $F$ et $G$. De plus, on peut choisir les notations pour qu'on ait
$$
F < O(F,G) < G < O(G,F) < F.
$$
\end{lemm}
Voici une esquisse de d\'emonstration. On peut construire un disque topologique ferm\'e $D$ dans la sph\`ere $\bbR^2 \cup \{\infty\}$, qui contient $ F$, et qui est disjoint de $G$ (le bord de $D$ passe par $0$ et $\infty$).
Soit $\wt F$ un relev\'e de $F$. 
%Notons $\tau$ l'automorphisme du rev\^etement. 
Au moyen du disque $D$, on montre que $\wt F$ est disjoint de $\tau(\wt F)$, et que
$$
\tau^{-1}(\wt F) < \wt F < \tau(\wt F).
$$
De plus, il existe un relev\'e $\wt G$ de $G$ tel que
$$
 \wt F < \wt G  < \tau(\wt F) < \tau (\wt G).
$$
On consid\`ere alors les ouverts $O(\wt F, \wt G)$ et $O(\wt G, \tau(\wt F))$ (d\'efinition~\ref{d.ofg}).  En projetant ces deux ouverts sur $\bbR^2 \setminus \{0\}$, on obtient les ensembles  $O(F,G)$ et $O(G,F)$ v\'erifiant les propri\'et\'es requises.

%%%%%%%%%%%%%%%%%%%%%%%%%%%%%%

 %%%%%%%%%%%%
 \clearpage
 \pagenumbering{roman}

%\pagebreak	
%\tableofcontents

\end{document}